\documentclass[11pt,twoside]{article}

\usepackage{amsbsy,amsfonts,amsmath,amssymb,eucal,mathrsfs}
\usepackage[all]{xy}
\usepackage{pstricks,epsfig}
\usepackage{pst-plot}

\newtheorem{definition}{Definition}[section]
\newenvironment{defi}{\begin{definition} \rm}{\end{definition}}

\newtheorem{prop}[definition]{Proposition}

\newtheorem{lemm}[definition]{Lemma}
\newtheorem{fact}[definition]{Fact}

\newtheorem{coro}[definition]{Corollary}
\newtheorem{theo}{Theorem}
\newtheorem{notation}[definition]{Notation}
\newenvironment{nota}{\begin{notation} \rm}{\end{notation}}
\newtheorem{construction}[definition]{Construction}

\newtheorem{remark}[definition]{Remark}
\newenvironment{rema}{\begin{remark} \rm}{\end{remark}}
\newtheorem{remarks}[definition]{Remarks}

\newtheorem{example}[definition]{Example}
\newenvironment{exam}{\begin{example} \rm}{\end{example}}
\newtheorem{examples}[definition]{Examples}

\newtheorem{nothing}[definition]{$\!\!$}

\newenvironment{proo}{{\flushleft \it Proof.}}{\hfill $\square$ \vspace{2mm}}

\newtheorem{conjecture}[definition]{Conjecture}

\newtheorem{definition*}{Definition}[section]
\newenvironment{defi*}{\begin{definition*} \rm}{\end{definition*}}
\newtheorem{definitions*}[definition*]{Definitions}
\newenvironment{defis*}{\begin{definitions*} \rm}{\end{definitions*}}
\newtheorem{prop*}[definition*]{Proposition}
\newtheorem{lemm*}[definition*]{Lemma}
\newtheorem{coro*}[definition*]{Corollary}
\newtheorem{theo*}[definition*]{Theorem}
\newtheorem{remark*}[definition*]{Remark}
\newenvironment{rema*}{\begin{remark*} \rm}{\end{remark*}}
\newtheorem{remarks*}[definition*]{Remarks}
\newenvironment{remas*}{\begin{remarks*} \rm}{\end{remarks*}}
\newtheorem{example*}[definition*]{Example}
\newenvironment{exam*}{\begin{example*} \rm}{\end{example*}}
\newtheorem{examples*}[definition*]{Examples}
\newenvironment{exams*}{\begin{examples*} \rm}{\end{examples*}}
\newtheorem{nothing*}[definition*]{$\!\!$}
\newenvironment{noth*}{\begin{nothing*} \rm}{\end{nothing*}}

\newtheorem{commentaire*}[definition*]{Commentaire}

\begin{document}

\def \pp {{\odot}}
\def \Rh {{\widehat{R}}}
\def \Sh {{\widehat{S}}}
\def \supp {{\rm Supp}}
\def \codim {{\rm Codim}}
\def \b {{\beta}}
\def \T {{\Theta}}
\def \t {{\theta}}
\def \L {{\cal L}}
\def \sca #1#2{\langle #1,#2 \rangle}
\def\pt{\{{\rm pt}\}}
\def\x {{\underline{x}}}
\def\y {{\underline{y}}}
\def\aut{{\rm Aut}}
\def\ra{\rightarrow}
\def\s{\sigma}\def\OO{\mathbb O}\def\PP{\mathbb P}\def\QQ{\mathbb Q}
 \def\CC{\mathbb C} \def\ZZ{\mathbb Z}\def\JO{{\mathcal J}_3(\OO)}
\newcommand{\G}{\mathbb{G}}
\def\proof{\noindent {\it Proof.}\;}
\def\qed{\hfill $\square$}
\def \uh {{\widehat{u}}}
\def \vh {{\widehat{v}}}
\def \uveeh {{\widehat{u^\vee}}}
\def \vveeh {{\widehat{v^\vee}}}
\def \fh {{\widehat{f}}}
\def \wh {{\widehat{w}}}
\def \wt {{\widetilde{w}}}
\def \Wh {{{W_{{\rm aff}}}}}
\def \Wt {{\widetilde{W}_{{\rm aff}}}}
\def \Qt {{\widetilde{Q}}}
\def \Waff {{W_{{\rm aff}}}}
\def \Waffm {{W_{{\rm aff}}^-}}
\def \Wpaff {{{(W^P)}_{{\rm aff}}}}
\def \Wtpaff {{{(\widetilde{W}^P)}_{{\rm aff}}}}
\def \Wtaffm {{\widetilde{W}_{{\rm aff}}^-}}
\def \lh {{\widehat{\lambda}}}
\def \pit {{\widetilde{\pi}}}
\def \lt {{{\lambda}}}
\def \xh {{\widehat{x}}}
\def \yh {{\widehat{y}}}
\def \a {\alpha}
\def \b {\beta}
\def \l {\lambda}
\def \t {\theta}


\newcommand{\expxy}{\exp_{x \to y}}
\newcommand{\drat}{d_{\rm rat}}
\newcommand{\dmax}{d_{\rm max}}
\newcommand{\zl}{Z(x,L_x,y,L_y)}


\newcommand{\N}{\mathbb{N}}
\newcommand{\A}{{\mathbb{A}_{\rm Aff}}}
\newcommand{\Ah}{{\mathbb{A}_{\rm Aff}}}
\newcommand{\At}{{\widetilde{\mathbb{A}}_{\rm Aff}}}
\newcommand{\Ht}{{{H}^T_*(\Omega K^{\ad})}}
\renewcommand{\H}{{\rm Hi}}
\newcommand{\Ih}{{I_{\rm Aff}}}
\newcommand{\psit}{{\widetilde{\psi}}}
\newcommand{\xit}{{\widetilde{\xi}}}
\newcommand{\Jt}{{\widetilde{J}}}
\newcommand{\Zt}{{\widetilde{Z}}}
\newcommand{\Xt}{{\widetilde{X}}}
\newcommand{\at}{{\widetilde{A}}}
\newcommand{\Z}{\mathbb Z}
\renewcommand{\S}{{\mathbb S}}
\newcommand{\fgcoad}{{G/P_\theta}}
\newcommand{\R}{\mathbb{R}}
\newcommand{\Q}{\mathbb{Q}}
\newcommand{\C}{\mathbb{C}}
\renewcommand{\O}{\mathbb{O}}
\newcommand{\F}{\mathbb{F}}
\newcommand{\p}{\mathbb{P}}
\newcommand{\co}{{\cal O}}
\newcommand{\pos}{{\bf P}}

\renewcommand{\a}{{\alpha}}
\newcommand{\az}{\a_\Z}
\newcommand{\ak}{\a_k}

\newcommand{\rc}{\R_\C}
\newcommand{\cc}{\C_\C}
\newcommand{\hc}{\H_\C}
\newcommand{\oc}{\O_\C}

\newcommand{\rk}{\R_k}
\newcommand{\ck}{\C_k}
\newcommand{\hk}{\H_k}
\newcommand{\ok}{\O_k}

\newcommand{\rz}{\R_Z}
\newcommand{\cz}{\C_Z}
\newcommand{\hz}{\H_Z}
\newcommand{\oz}{\O_Z}

\newcommand{\RR}{\R_R}
\newcommand{\CR}{\C_R}
\newcommand{\HR}{\H_R}
\newcommand{\OR}{\O_R}

\newcommand{\re}{\mathtt{Re}}

\newcommand{\matttr}[9]{
\left (
\begin{array}{ccc}
{} \hspace{-.2cm} #1 & {} \hspace{-.2cm} #2 & {} \hspace{-.2cm} #3 \\
{} \hspace{-.2cm} #4 & {} \hspace{-.2cm} #5 & {} \hspace{-.2cm} #6 \\
{} \hspace{-.2cm} #7 & {} \hspace{-.2cm} #8 & {} \hspace{-.2cm} #9
\end{array}
\hspace{-.15cm}
\right )   }


\newcommand{\dual}{{\bf v}}
\newcommand{\com}{\mathtt{Com}}
\newcommand{\rg}{\mathtt{rg}}
\newcommand{\pu}{{\mathbb{P}^1}}
\newcommand{\scal}[1]{\langle #1 \rangle}
\newcommand{\MK}[2]{{\overline{{\rm M}}_{#1}(#2)}}
\newcommand{\mor}[2]{{{\rm Mor}_{#1}(\pu,#2)}}

\newcommand{\fg}{\mathfrak g}
\newcommand{\fgad}{G/P_\T}
\renewcommand{\fh}{\mathfrak h}
\newcommand{\fu}{\mathfrak u}
\newcommand{\fz}{\mathfrak z}
\newcommand{\fn}{\mathfrak n}
\newcommand{\fe}{\mathfrak e}
\newcommand{\fp}{\mathfrak p}
\newcommand{\ft}{\mathfrak t}
\newcommand{\fl}{\mathfrak l}
\newcommand{\fq}{\mathfrak q}
\newcommand{\fsl}{\mathfrak {sl}}
\newcommand{\fgl}{\mathfrak {gl}}
\newcommand{\fso}{\mathfrak {so}}
\newcommand{\fsp}{\mathfrak {sp}}
\newcommand{\ff}{\mathfrak {f}}

\newcommand{\ad}{{\rm ad}}
\newcommand{\jad}{{j^\ad}}
\newcommand{\id}{{\rm id}}


\newcommand{\dynkinadeux}[2]
{
$#1$
\setlength{\unitlength}{1.2pt}
\hspace{-3mm}
\begin{picture}(12,3)
\put(0,3){\line(1,0){10}}
\end{picture}
\hspace{-2.4mm}
$#2$
}

\newcommand{\mdynkinadeux}[2]
{
\mbox{\dynkinadeux{#1}{#2}}
}

\newcommand{\dynkingdeux}[2]
{
$#1$
\setlength{\unitlength}{1.2pt}
\hspace{-3mm}
\begin{picture}(12,3)
\put(1,.8){$<$}
\multiput(0,1.5)(0,1.5){3}{\line(1,0){10}}
\end{picture}
\hspace{-2.4mm}
$#2$
}

\newcommand{\poidsesix}[6]
{
\hspace{-.12cm}
\left (
\begin{array}{ccccc}
{} \hspace{-.2cm} #1 & {} \hspace{-.3cm} #2 & {} \hspace{-.3cm} #3 &
{} \hspace{-.3cm} #4 & {} \hspace{-.3cm} #5 \vspace{-.13cm}\\
\hspace{-.2cm} & \hspace{-.3cm} & {} \hspace{-.3cm} #6 &
{} \hspace{-.3cm} & {} \hspace{-.3cm}
\end{array}
\hspace{-.2cm}
\right )      }

\newcommand{\copoidsesix}[6]{
\hspace{-.12cm}
\left |
\begin{array}{ccccc}
{} \hspace{-.2cm} #1 & {} \hspace{-.3cm} #2 & {} \hspace{-.3cm} #3 &
{} \hspace{-.3cm} #4 & {} \hspace{-.3cm} #5 \vspace{-.13cm}\\
\hspace{-.2cm} & \hspace{-.3cm} & {} \hspace{-.3cm} #6 &
{} \hspace{-.3cm} & {} \hspace{-.3cm}
\end{array}
\hspace{-.2cm}
\right |      }

\newcommand{\poidsesept}[7]{
\hspace{-.12cm}
\left (
\begin{array}{cccccc}
{} \hspace{-.2cm} #1 & {} \hspace{-.3cm} #2 & {} \hspace{-.3cm} #3 &
{} \hspace{-.3cm} #4 & {} \hspace{-.3cm} #5 & {} \hspace{-.3cm} #6
\vspace{-.13cm}\\
\hspace{-.2cm} & \hspace{-.3cm} & {} \hspace{-.3cm} #7 &
{} \hspace{-.3cm} & {} \hspace{-.3cm}
\end{array}
\hspace{-.2cm}
\right )      }

\newcommand{\copoidsesept}[7]{
\hspace{-.12cm}
\left |
\begin{array}{cccccc}
{} \hspace{-.2cm} #1 & {} \hspace{-.3cm} #2 & {} \hspace{-.3cm} #3 &
{} \hspace{-.3cm} #4 & {} \hspace{-.3cm} #5 & {} \hspace{-.3cm} #6
\vspace{-.13cm}\\
\hspace{-.2cm} & \hspace{-.3cm} & {} \hspace{-.3cm} #7 &
{} \hspace{-.3cm} & {} \hspace{-.3cm}
\end{array}
\hspace{-.2cm}
\right |      }

\newcommand{\poidsehuit}[8]{
\hspace{-.12cm}
\left (
\begin{array}{ccccccc}
{} \hspace{-.2cm} #1 & {} \hspace{-.3cm} #2 & {} \hspace{-.3cm} #3 &
{} \hspace{-.3cm} #4 & {} \hspace{-.3cm} #5 & {} \hspace{-.3cm} #6 &
{} \hspace{-.3cm} #7   \vspace{-.13cm}\\
\hspace{-.2cm} & \hspace{-.3cm} & {} \hspace{-.3cm} #8 &
{} \hspace{-.3cm} & {} \hspace{-.3cm}
\end{array}
\hspace{-.2cm}
\right )      }

\newcommand{\copoidsehuit}[8]{
\hspace{-.12cm}
\left |
\begin{array}{cccccc}
{} \hspace{-.2cm} #1 & {} \hspace{-.3cm} #2 & {} \hspace{-.3cm} #3 &
{} \hspace{-.3cm} #4 & {} \hspace{-.3cm} #5 & {} \hspace{-.3cm} #6 &
{} \hspace{-.3cm} #7  \vspace{-.13cm}\\
\hspace{-.2cm} & \hspace{-.3cm} & {} \hspace{-.3cm} #8 &
{} \hspace{-.3cm} & {} \hspace{-.3cm}
\end{array}
\hspace{-.2cm}
\right |      }

\newcommand{\im}{\mathtt{Im}}


\def\cA{{\cal A}} \def\cC{{\cal C}} \def\cD{{\cal D}} \def\cE{{\cal E}}
\def\cF{{\cal F}} \def\cG{{\cal G}} \def\cH{{\cal H}} \def\cI{{\cal I}}
\def\cK{{\cal K}} \def\cL{{\cal L}} \def\cM{{\cal M}} \def\cN{{\cal N}}
\def\cO{{\cal O}}
\def\cP{{\cal P}} \def\cQ{{\cal Q}} \def\cT{{\cal T}} \def\cU{{\cal U}}
\def\cV{{\cal V}} \def\cX{{\cal X}} \def\cY{{\cal Y}} \def\cZ{{\cal Z}}

 \title{On the quantum cohomology of adjoint varieties}
 \author{P.E. Chaput, N. Perrin}

\maketitle

\begin{abstract}
We study the quantum cohomology of quasi-minuscule and quasi-cominuscule
homogeneous spaces. The product of any two Schubert cells does
not involve powers of the quantum parameter higher than 2. With the help of
the quantum to classical principle we give presentations of the quantum
cohomology algebras. These algebras are semi-simple for adjoint
non coadjoint varieties and some properties of the induced strange
duality are shown.
\end{abstract}

 {\def\thefootnote{\relax}
 \footnote{ \hspace{-6.8mm}
 }
 }

\section{Introduction}

In this paper, we study the quantum and classical cohomology rings of some 
rational homogeneous varieties. The papers \cite{cmp}, \cite{CMP1}, 
\cite{cmp3} dealt with the quantum cohomology of minuscule and cominuscule 
rational homogeneous spaces. Here we are interested in quasi-minuscule and 
quasi-cominuscule varieties (see Definition \ref{quasi-min}). 

The choice of 
these varieties comes from the fact that Schubert calculus has a simpler 
and in many cases combinatorial expression: jeu de Taquin and its 
generalisations. This was remarked for (co)minuscule varieties by H. Thomas 
and A. Yong in \cite{TY} and generalised for a large class of Schubert
classes (called $\Lambda$-minuscule, see Definition \ref{defi-minuscule}),
in \cite{littlewood}. 
In particular, for quasi-(co)minuscule varieties, all Schubert classes up to
half the dimension are $\Lambda$-minuscule (see Corollary 
\ref{l-min-dim-moit}). This enables computing
presentations of the classical cohomology ring of these varieties.

Moreover, as for the (co)minuscule case, we will
explain here that the quantum Schubert calculus is also simpler for these 
varieties, and give several technics for computing 
Gromov-Witten invariants. The first one (see Theorem 
\ref{theo-chevalley}) is a simple interpretation of the quantum Chevalley 
formula as given in \cite{FW}. We then prove a result on the power of $q$, 
the quantum parameter, appearing in a product of two Schubert classes. 
In particular for $X$ an adjoint of coadjoint variety (see Definition 
\ref{quasi-min}) we prove a vanishing result for the quantum Schubert calculus
and relate it to the geometry of these varieties:

\begin{theo}       
Let $\s_u$ and $\s_v$ be classical cohomology classes. Then if we write
$$\s_u*\s_v=\sum_{d,w}c_{u,v}^w(d)q^d\s_w\ \ ,$$ 
we have $c_{u,v}^w(d)=0$ unless $d\leq 2$.
\end{theo}

We also explain that the quantum to classical principle as initiated in 
\cite{BKT} can also be applied to compute some Gromov-Witten invariants. With 
all these technics, we are able to give presentations of the quantum 
cohomology rings of all quasi-(co)minuscule homogeneous spaces.
These results are 
presented in Section \ref{section-pres}. We do not give here a presentation in 
the classical cases as they were obtained by A. Buch, A. Kresch and H. 
Tamvakis in \cite{BKTnouveau} but focus on exceptional groups. Some of the 
computations were made using a computer program written in Java
(and available at 
the webpage
\texttt{www.math.sciences.univ-nantes.fr/}\verb|~|\texttt{chaput/quiver-demo.html}),
which may be
of some interest for those who want to make explicit computations in these
quantum cohomology rings.

Finally we focus on strange duality, a known phenomenon for (co)minuscule 
homogeneous spaces $X$: in that case there exists an involution 
$\iota$ of the quantum cohomology algebra $QH^*(X,\C)_{loc}$ with the quantum 
parameter $q$ inverted sending classes of degree $d$ to classes of degree 
$-d$. This involution, in the (co)minuscule cases also sends a Schubert 
class to a multiple of a Schubert class. The existence of such an involution 
is a consequence of the semi-simplicity of the ring $QH^*(X,\Z)_{loc}$. We 
conjectured in \cite{cmp3} that the only homogeneous spaces $X$ with 
$QH^*(X,\Z)_{loc}$ semi-simple are the (co)minuscule homogeneous spaces.

In Section \ref{section-ss} we are able to prove, for many examples, that for 
$X$ a rational homogeneous space with Picard number one the algebra 
$QH^*(X,\Z)_{loc}$ is not semi-simple (see Theorem 
\ref{theo-non-semisimple}). However, we give counterexamples to our 
conjecture and prove that for $X$ an adjoint non coadjoint variety,  
the localised algebra $QH^*(X,\C)_{loc}$ \emph{is} semi-simple (see Theorem 
\ref{ss-adjoint}). This implies the existence of an algebra involution 
$\iota$ as above. We prove that this involution, once specialised
at $q=1$, is the identity on classes of degree multiple of $c_1(X)$ 
(Proposition \ref{prop-inv-c1}) but is more complicated in 
general than the involution in the minuscule and cominuscule cases (see 
Theorem \ref{pas-classe}). In particular, we are not able to compute
it except for very special classes.

\tableofcontents

\newcommand{\Paff}{P^{{\rm aff}}}
\newcommand{\Raff}{R^{{\rm aff}}}
\newcommand{\aad}{\a_{ad}}
\newcommand{\g}{\gamma}

\section{Minuscule, Cominuscule, Adjoint and Coadjoint varieties}
\label{section-hasse}

Let $G$ be a semi-simple algebraic group. Let $T$ be a maximal torus
in $G$ and $B$ a Borel subgroup containing $T$. We denote by $W$ the
associated Weyl group, with $\Delta$ the root system and with $S$ the
set of simple roots. We shall denote by $\T$ the highest root and
with $\t$ the highest short root. For simply-laced groups, we have
$\T=\t$. We will denote by $\a_{ad}$ any simple root such that
$\scal{\a_{ad},\T^\vee}=1$ (such a root is unique except in type $A$).

\subsection{First definitions}

\begin{defi}
\label{quasi-min} Let $\varpi$ be a dominant weight, we shall call $\varpi$

$\bullet$ \emph{minuscule} if for any positive root $\a\in R_+$ we have
$\sca{\a^\vee}{\varpi}\leq1$ ;

$\bullet$ \emph{cominuscule} if it is minuscule for the dual root system ;

$\bullet$ \emph{quasi-minuscule} if for any positive root $\a\in R_+$ we have
$\sca{\a^\vee}{\varpi}\leq2$ with equality iff 
$\a=\varpi$ ;

$\bullet$ \emph{quasi-cominuscule} if it is
quasi-minuscule for the dual root system.
\end{defi}


\begin{rema}
A non minuscule quasi-minuscule weight $\varpi$ is a root. More
precisely we have $\varpi=\t$. A quasi-cominuscule non cominuscule weight
$\varpi$ is also a root and we have $\varpi=\T$. In this
case we shall say that $\varpi$ is \emph{adjoint}. We call a weight
\emph{coadjoint} if it is adjoint for the dual root system, in this case
we have $\varpi = \theta$. 

In simply laced cases, minuscule and cominuscule weights, quasi-minuscule 
and quasi-cominuscule weights and adjoint and coadjoint weights coincide.
\end{rema}

\centerline{\begin{tabular}{c|cc|cc}
\hline {Type} & \multicolumn{2}{c|}{Quasi-minuscule} &
\multicolumn{2}{c}{Quasi-cominuscule} \\
of $G$ & Minuscule & Coadjoint & Cominuscule& Adjoint\\
\hline $A_n$& $\varpi_1,\cdots,\varpi_n$ &$\varpi_1+\varpi_n$&
$\varpi_1,\cdots,\varpi_n$&$\varpi_1+\varpi_n$\\
$B_n$& $\varpi_n$ &$\varpi_1$&
$\varpi_1$&$\varpi_2$\\
$C_n$& $\varpi_1$ & $\varpi_2$ &
$\varpi_n$& $\varpi_1$\\
$D_n$& $\varpi_1,\varpi_{n-1}$, $\varpi_n$ &$\varpi_2$&
$\varpi_1,\varpi_{n-1}$, $\varpi_n$&$\varpi_2$\\
$E_6$& $\varpi_1$, $\varpi_6$ & $\varpi_2$&
$\varpi_1$, $\varpi_6$&$\varpi_2$\\
$E_7$& $\varpi_7$ &$\varpi_1$&
$\varpi_7$&$\varpi_1$\\
$E_8$& None &$\varpi_8$&
None&$\varpi_8$\\
$F_4$& None &$\varpi_4$&
None&$\varpi_1$\\
$G_2$& None &$\varpi_1$& None &$\varpi_2$\\
\hline
\end{tabular}}

\vskip 0.2 cm

\centerline{Minuscule, Quasi-minuscule, Cominuscule, Quasi-cominuscule, 
Adjoint and Coadjoint weights.}

\begin{defi}
Let $P$ be the standard parabolic subgroup associated to a weight
$\varpi$. We will say that the homogeneous variety $G/P$ is
\emph{minuscule}, \emph{cominuscule}, \emph{quasi-minuscule}, 
\emph{quasi-cominuscule},
\emph{adjoint} or \emph{coadjoint} if the weight
$\varpi$ is minuscule, cominuscule,
quasi-minuscule, quasi-cominuscule, adjoint or coadjoint.
\end{defi}

\begin{rema}
(\i) The adjoint homogeneous spaces correspond to the
classical terminology of adjoint varieties.

(\i\i) The geometry of minuscule homogeneous varieties is
better understood than the geometry of general homogeneous spaces.
See for example \cite{seshadri}, \cite{TY} or \cite{small}.
Quasi-minuscule weights were in particular introduced to generalise
the classical theory of standard monomials in minuscule varieties
(see \cite{seshadri}) to other homogeneous spaces (see \cite{LMS}).
\end{rema}

\begin{rema}
For minuscule and cominuscule homogeneous spaces, the quantum
cohomology was described in \cite{BKT} for classical groups and in
\cite{CMP1} in general. So we will mainly focus in the sequel on
adjoint and coadjoint homogeneous spaces except in Subsection
\ref{q-hasse} where we include the minuscule and cominuscule cases in
our simplified version of quantum Chevalley formula.
\end{rema}

The adjoint and coadjoint varieties are tabulated in the following array
(for the simply-laced Lie algebras the notions of adjoint and coadjoint
varieties coincide and therefore we leave the corresponding box empty):

$$\begin{array}{c|cccc|cccc}
&\multicolumn{4}{|c|}{\textrm{Adjoint}}&
\multicolumn{4}{|c}{\textrm{Coadjoint}}
\\
\hline \vspace{-.4cm} &&&&&&&& \\
G & \Theta & \fgad & \dim(\fgad) & c_1(\fgad) & \theta & \fgcoad &
\dim(\fgcoad) & c_1(\fgcoad) \\ 
\hline
A_n & \varpi_1 + \varpi_n & \F(1,n,n+1) & 2n+1 & (n+1,n+1) \\
B_n & \varpi_2 & \G_Q(2,2n+1) & 4n-5 & 2n-2 & \varpi_1 & \Q^{2n-1} & 2n-1 & 2n-1\\
C_n & 2\varpi_1 & \p^{2n-1} & 2n-1 & n & \varpi_2 & \G_\omega(2,2n) & 4n-5 & 2n-1 \\
D_n & \varpi_2 & \G_Q(2,2n) & 4n-7 & 2n-3 \\
E_6 & \varpi_2 & & 21 & 11 \\
E_7 & \varpi_1 & & 33 & 17 \\
E_8 & \varpi_8 & & 57 & 29 \\
F_4 & \varpi_1 & & 15 & 8 & \varpi_4 & & 15 & 11 \\
G_2 & \varpi_2 & & 5 & 3 & \varpi_1 & \Q^5 & 5 & 5 \\
\hline
\end{array}$$
\begin{equation}
\label{array-adjoint}
\textrm{Invariants of Adjoint and Coadjoint varieties.}
\end{equation}

In this array, we denote by $\F(1,n,n+1)$, $\G_Q(2,n)$,
$\G_\omega(2,2n)$ resp. $\Q^n$ the point-hyperplane incidence in
$\p^n$, the Grassmannian of 2 dimensional isotropic subspaces in an
$n$-dimensional subspace endowed with a nondegenerate quadratic form, 
the Grassmannian of 2 dimensional isotropic subspaces in an
$2n$-dimensional subspace endowed with a symplectic
form resp. a $n$-dimensional quadric. We denote by $c_1(X)$ the
index of the variety $X$. In the case of type $C_n$, the adjoint
variety is isomorphic to $\p^{2n-1}$, but the projective embedding is
given by twice the generator of ${\rm Pic}(\p^{2n-1})$. Note that we
always have $\dim(\fgad) = 2 c_1(\fgad) - 1$ (this well-known fact can
be seen as a consequence of Propositions \ref{prop_y2_adjoint} and
\ref{prop_dim_yd}).

\subsection{Roots and Schubert varieties in the adjoint and
  coadjoint cases}

In this subsection, we recall some 
results on Schubert varieties specific to the adjoint and coadjoint
homogeneous spaces. We shall need the following:

\begin{nota}
(\i) We define the affine weight $\rho$ by $\scal{\rho,\a^\vee}=1$ and
the affine coweight $\rho^\vee$ by $\scal{\a,\rho^\vee}=1$ for all
simple root $\a$ of the affine Weyl group.

(\i\i) Recall that there is a natural order on the roots defined by
$\a\leq\b$ if $\b-\a$ is in the monoid generated by positive roots.
\end{nota}

\begin{fact}
Let $X=G/P$ be a coadjoint (respectively an adjoint)
homogeneous space. Then, the Schubert varieties are in one to one
correspondence with short (respectively long) roots. The
corresponding map from $W^P$ to this set of roots is given by
$u\mapsto \a=u(\varpi)$ and we set $\s_\a=\s_u$.
\end{fact}

\begin{proo}
This is a reformulation of the fact that the orbit of the point
$\varpi$ in $X$ under the action of the Weyl group is the set of
$T$-fixed points in $X$.
\end{proo}

Let $u\in W^P$ and $\a=u(\varpi)$. We denote by $X(\a)$ the Schubert
variety such that $[X(\a)]=\s_\a$ in $H^{2l(u)}(X,\Z)$.

\begin{prop}
\label{chev}
Let $X=G/P$ be a coadjoint (respectively an adjoint)
homogeneous space and let $\a$ be a short (respectively long) root.
Denote by $X(\a)$ the associated Schubert variety. Then, the
following propositions hold:

(\i) The integer $\dim(X)$ is odd, we define $r$ by the equality 
$\dim(X)=2r+1$.


(\i\i) There is an equivalence 
$$X(\a)\subset X(\b)\Leftrightarrow \left\{
\begin{array}{ll}
\a\leq \b & \textrm{for $\a$ and $\b$ of the same sign,}\\ 
\supp(\a)\cup\supp(\b) \textrm{ is connected} & \textrm{for $\a$
  negative and $\b$ positive.} \end{array}\right.$$ 




(\i\i\i) Denote by $X(\a)^\vee$ the Poincar{\'e} dual of $X(\a)$ and by $i$ the
Weyl involution, then we have the equality $X(\a)^\vee=X(-i(\a))$.


\end{prop}

\begin{proo}
All these results can be easily deduced from results in \cite[Section
3.1]{LMS} in the coadjoint case. Indeed, $(\imath)$ follows from
Theorem 3.1 in \emph{loc. cit.}, $(\imath\imath)$ follows
from the conjonction of Theorem 3.1 and Theorem 3.10 in
\emph{loc. cit.}. For $(\imath\imath\imath)$, recall that Poincar{\'e}
duality is given on 
$W^P$ by $u\mapsto w_0uw_0w_X$ where $w_0$ in the longest element in
$W$ and $w_X$ the longest element in $W^P$. In terms of roots, if
$\a=u(-\varpi)$, then the root corresponding to $X(\a)^\vee$ is
$$w_0uw_0w_X(-\varpi)=w_0uw_0(\varpi)=w_0u(-\varpi)=w_0(\a)=-i(\a).$$

To prove the corresponding statements in the adjoint case, little more
is needed. Indeed, all the results only involve Weyl groups and are
therefore true simply because the Weyl group of a root system and the
Weyl group of its dual are isomorphic. 
\end{proo}

\subsection{Description of $\varpi$-minuscule elements in the adjoint
  and coadjoint cases}

Let us recall, according to Dale Peterson
\cite[p.273]{proctor_minuscule} (see also \cite{littlewood}), the
following definition:

\begin{defi}
\label{defi-minuscule}
Let $\Lambda = \sum_{i} \Lambda_i \varpi_i$ 
be a dominant weight (where $(\varpi_i)_i$ are the fundamental weights).
\begin{itemize}
\item
An element $w \in W$ is $\Lambda$-minuscule if there exists a
reduced decomposition
$w = s_{i_1} \cdots s_{i_l}$ such that for any $k \in [1,l]$ we have
$s_{i_k} s_{i_{k+1}}\cdots s_{i_l}(\Lambda) 
= s_{i_{k+1}}\cdots s_{i_l}(\Lambda) - \alpha_{i_k}$.
\item
$w$ is $\Lambda$-cominuscule if $w$ is $(\sum \Lambda_i 
\varpi_i^\vee)$-minuscule.
\end{itemize}
\end{defi}

\begin{rema}
  For $P$ a parabolic subgroup associated to a (co)minuscule weight
  $\varpi$, then any element in $W^P$ is $\varpi$-(co)minuscule.
\end{rema}

Let $X=G/P$ be a coadjoint or an adjoint variety associated to the
fundamental weight $\varpi$. Recall that there
exists a non negative integer $r$ such that $\dim(X)=2r+1$.

\begin{prop}
\label{lambda-min-adjoint}
If $X$ is coadjoint (resp. adjoint), then any element $w$ in $W^P$
with $l(w)\leq r$ is $\varpi$-minuscule (resp. $\varpi$-cominuscule).
%
%
\end{prop}

\begin{proo}
The proof of Theorem 3.1 in \cite{LMS} proves that, in the
coadjoint case, such a $u$ is $\varpi$-minuscule and that it is
maximal for this property.
The adjoint case is obtained from the previous one by taking the dual
root system.
%
\end{proo}

We refer to \cite{littlewood} for a definition of jeu de taquin. As a
consequence of the main result in \emph{loc. cit.} 
we have: 

\begin{coro}
\label{l-min-dim-moit}
  Let $u$, $v$ and $w$ be elements in $W^P$, then the
  Littlewood-Richardson coefficient $c_{u,v}^w$ can be computed using
  jeu de taquin for $l(w)\leq r$.
\end{coro}

\subsection{Quantum Hasse diagram}
\label{q-hasse}

In this subsection, we shall give a simplified statement for the
quantum Chevalley formula as formulated in \cite{FW}. We shall also
see that our parametrisation of Schubert classes by roots as in the
previous subsection extends to the quantum monomials, which are by
definition of the form $q^d\s_u$ for some integer $d$ and some $u\in
W^P$: these monomials are parametrised by the roots 
of the affine root system $\Rh$ associated to $R$ and this 
parametrisation is compatible with the quantum Chevalley formula.

For the moment we only assume that $X$ has Picard number one.
Let us first define an injection of the set of
quantum monomials in the set of affine weights $\Paff$.

\begin{defi}\
\label{defi-eta}

\begin{itemize}
\item
We denote $\cM$ the set of quantum monomials (with $q$ inverted),
namely the set of pairs $(w,d)$ in $W^P \times \Z$.
\item
Let $w \in W^P$ and $d \in \Z$.
To the quantum monomial $\s_w \cdot q^d$ we associate the affine weight
$\eta(w,d) =  u(\varpi) - d\delta$, 
where $\delta$ is the minimal
positive imaginary root in $\Rh$ (see for example \cite{kac}).
\item
For $\pi = \eta(w,d)$ we denote $\s_\pi = \s_w \cdot q^d$ the corresponding
quantum monomial.
\item
For $\pi \in \eta(\cM)$ we denote $l(\pi)$ the integer $l(w)+dc_1$, if
$(w,d)$ is the element such that $\eta (\s_w,d) = \pi$.
\end{itemize}
\end{defi}

\begin{rema}
  In the coadjoint case, $\eta(\cM)$ is the set of all short
real roots while in the adjoint case it is the set of all long real
roots.
\end{rema}

\begin{prop}
\label{prop-fulton}
Let $\pi \in \eta(\cM)$; we have
$$
h * \s_\pi = \sum_{\g>0,l(s_\g(\pi)) = l(\pi)+1} \scal{\pi,\g^\vee}
\s_{s_\g(\pi)}.
$$
\end{prop}
\begin{proo}
This is the particular case, when the Picard number of $X$ is one, of
\cite[Theorem 10.1]{fulton}. In fact, let us assume that
$\pi = \eta(w,0)$ with $w \in W^P$;
according to this result, the
quantum product of $h$ and $\s_w$ is
\begin{equation}
\label{equa-fulton}
\sum \scal{\varpi,\a^\vee} \s_{us_\a} +
\sum q^{d(\a)} \scal{\varpi,\a^\vee} \s_{us_\a}
\end{equation}
where $\a$ is a positive root, $d(\a) = \scal{\varpi,\a^\vee}$, and
in the first resp. second summand we have $l(ws_\a) = l(w) + 1$ resp.
$l(ws_\a) = l(u)+1 - c_1d(\a)$. We denote $A$ resp. $A_1,A_2$ the set
of these roots resp. the roots in the first, second case. We consider
an injection $\varphi : A \to \Rh$ by mapping $\a \in A_1$ to
$w(\a)$ and $\a \in A_2$ to $w(\a)+\delta$.

Let $\a \in A$ and let us denote $\g = \varphi(\a)$.
We claim that for $\a \in A_1$ resp. $\a \in A_2$,
$s_\g w(\varpi)$ is equal to $ws_\a (\varpi)$ resp. 
$ws_\a (\varpi) - d(\a)\delta$. In fact the first case follows from the
equality $ws_\a = ws_\a w^{-1}w = s_{w(\a)}w$. The second case may be computed
as follows:
$$
\begin{array}{rcl}
s_\g w(\varpi) & = & s_{w(\a)+\delta} ( \varpi ) =
w(\varpi) - \scal{w(\a)^\vee,w(\varpi)}(w(\a) + \delta) \\
& = & w(\varpi) - \scal{\a^\vee,\varpi}w(\a)-d(\a)\delta =
ws_\a (\varpi) - d(\a)\delta\ ,
\end{array}
$$
where the last equality follows from the first case. Thus in particular
we get $s_\g w(\varpi) = \eta ( ws_\a , d(\a) )$.

Moreover it follows that $\g$ is a positive root. In fact, in the second
case it has positive coefficient on $\delta$, and in the first case we have
$s_\g w(\varpi) = w(\varpi) - d(\a) \g$ and on the other hand
$s_\g w(\varpi) = ws_\a (\pi) < w(\a)$, so that $\g > 0$.

Finally $\varphi$ is a bijection from $A$ to the set of positive roots $\g$
such that $l(s_\g w(\pi)) = l(w(\pi)) + 1$, and the proposition is proved.
\end{proo}

Our next goal is to restrict the set of root $\g$ in the former
proposition in the quasi-(co)minuscule cases. 
We express the index $c_1$ of our varieties in terms of the highest roots.
We start with a:

\begin{lemm}
\label{lemm-c1-general}
Let $X=G/P$ be a homogeneous space corresponding to a fundamental
weight $\varpi$ and let $\beta$ be the simple root with
$\scal{\varpi,\beta^\vee}=1$. Let $\a^\vee$ be the biggest coroot with
the following properties:
\begin{itemize}
\item
$\scal{\varpi,\a^\vee} = 1$.
\item
$|\a^\vee| = |\beta^\vee|$.
\end{itemize}
Then we have $c_1(X) = \scal{\rho,\a^\vee}+1$.
\end{lemm}
\begin{proo}
Let us recall (see \cite[Page 85]{BK}) that
$c_1(X)=\scal{\delta_P,\beta^\vee}$
where $\delta_P=\rho+w^P\rho$ and
$w^P$ is the longest element in $W_P$ the Weyl group of $P$.
To prove the lemma it suffices to explain why $w^P(\beta^\vee)=\a^\vee$. The
coroot $w^P(\beta^\vee)$ must have the same length as $\beta^\vee$.

First let us prove
that the set $B$ of coroots $\a^\vee$ of the same length as $\beta^\vee$ and with
$\scal{\varpi,\a^\vee}=1$ is equal to the orbit of $\beta^\vee$ under
$W_P$. Clearly the orbit is contained in $B$. Now let $\beta_0^\vee\in B$. We
proceed by induction on $\scal{\beta_0^\vee,\rho}$. If
$\scal{\beta_0^\vee,\rho}=1$, then $\beta_0^\vee=\beta^\vee$
and we are done. Otherwise,
there exists $\gamma^\vee$ a simple coroot of $P$ with
$\scal{\gamma,\beta_0^\vee}>0$. Indeed, otherwise, we would have
$\scal{\rho_P,\beta_0^\vee}\leq 0$ where $\rho_P$ is half the sum of the
roots of the Levi part $L(P)$ of $P$. But because $\beta_0^\vee$ is
positive this
would imply $\scal{\rho_P,\beta_0^\vee}=0$ thus $\beta_0^\vee=\beta^\vee$. Now
$s_{\gamma^\vee}(\beta_0^\vee)$ is again in $B$ with
$\scal{\rho,s_{\gamma^\vee}(\beta_0^\vee)}<\scal{\rho,\beta_0^\vee}$
and the result on $B$ follows by induction.

We end by proving that $w^P(\beta^\vee)$ is the
highest coroot of $B$. Indeed, if $u$ and $v$ are in $W_P$ with $u\leq
v$, and if $\varpi'$ is a fundamental weight of $L(P)$, then
$u(\varpi')\geq v(\varpi')$ i.e. $u(\varpi')-v(\varpi')$ is a sum of
positive roots of $L(P)$. We get
$\scal{u(\varpi')-v(\varpi'),\beta^\vee}\leq0$ thus
$\scal{u(\beta),{\varpi'}^\vee}\leq\scal{v(\a),{\varpi'}^\vee}$ and
$u(\beta^\vee)\leq v(\beta^\vee)$. 
\end{proo}

If $R$ is a root system, remark that the coroot $\T^\vee$ is the
highest short root of the dual root system while $\t^\vee$ is the
highest root of the dual root system.
From now on we assume that $X$ is quasi-(co)minuscule and has Picard
rank one. We have:

\begin{lemm}
\label{lemm-c1}
The following formulas yield the first Chern 
class of minuscule, cominuscule, coadjoint or adjoint varieties:

\vskip 0.4 cm

\centerline{
\begin{tabular}{c|cccc}
  & Minuscule & Coadjoint & Cominuscule & Adjoint\\
\hline
$c_1$ & $\sca{\rho}{\t^\vee}+1$ & $\sca{\rho}{\t^\vee}$ & 
$\sca{\rho}{\T^\vee}+1$ & $\sca{\rho}{\T^\vee}$ \\
$c_1$ & $\scal{\T,\rho^\vee}+1$ & $\scal{\T,\rho^\vee}$ \\
\end{tabular}}
\end{lemm}
\begin{proo}
If $G$ is of type $A$, then $X$ is a Grassmannian and the lemma is
true. If not, let $\aad^\vee$ denote the adjoint simple coroot, that
it the only coroot such that $\scal{\aad,\T^\vee}=1$. 
To prove the first line it is enough to check that the root
$\a$ of Lemma \ref{lemm-c1-general}
is $\t^\vee,\t^\vee-\aad^\vee,\T^\vee,\T^\vee-\aad^\vee$, according to
the four cases.

In the minuscule case, the simple root $\beta$ with
$\scal{\varpi,\beta^\vee}=1$ is short and therefore $\beta^\vee$ is long; 
moreover by definition of a minuscule weight we have
$\scal{\varpi,\g^\vee} \leq 1$ for any coroot $\g^\vee$.
It thus follows that $\a^\vee$ is the highest coroot.

In the cominuscule case, $\beta^\vee$ must be short. Moreover for any short
coroot $\g^\vee$ we have $\scal{\varpi,\g^\vee} \leq 1$. Thus in this case
$\a^\vee = \T^\vee$.

The coadjoint and adjoint cases are similar, except that the highest
coroot $\g^\vee$
of the same length as $\beta^\vee$ has coefficient 2 on $\beta^\vee$,
so that we must have $\a^\vee = \g^\vee - \aad^\vee$.

\vskip .2cm
To prove the second line of the table, we only need to apply the
quasi-cominuscule case to the dual root system. We may also observe
that the sum of the coefficients of the highest root is the same in a
root system and its dual (for the minuscule and coadjoint cases).
\end{proo}

\begin{lemm}
\label{lemm-length}
If $[a/b]$ denotes the integral part of $a/b$, we have the following
equalities:

\vskip .4cm

\centerline{
\begin{tabular}{c|cc}
 & Minuscule & Coadjoint \\
\hline
$\scal{\varpi-\pi,\rho^\vee}$ & $l(\pi)$ & $l(\pi)+[l(\pi)/c_1]$
\vspace{.3cm} \\
 & Cominuscule & Adjoint \\
\hline
$\scal{\rho,(\varpi-\pi)^\vee}$ & $l(\pi)$ & $l(\pi)+[l(\pi)/c_1]$ \\
\end{tabular}
}
\end{lemm}
\begin{proo}
Note that we have $\scal{\delta,\rho^\vee} = \scal{\T,\rho^\vee}+1$
since $\delta = \T + \a_0$.
Thus using Lemma \ref{lemm-c1} and the definition of $l(\pi)$ given in
Definition \ref{defi-eta},
it is enough to prove the lemma when
$\pi = \eta(w,0)$.

In this case, assuming we are in the minuscule case, if we write
$w = s_{i_l} \cdots s_{i_1}$, we have
$w(\varpi) = \varpi - \sum \a_{i_k}$, and so
$\scal{\varpi-w(\varpi),\rho^\vee} = l = l(w)$. Thus the lemma is proved in
this case. The cominuscule case follows by duality.

In the coadjoint case, if $l(w)<c_1$, then any reduced expression
of $w$ will have the same property as in the minuscule case and the same
argument gives $\scal{\varpi-w(\varpi),\rho^\vee} = l(w)$. If $l(w) \geq c_1$,
then we can find a reduced expression
$w = s_{i_l} \cdots s_{i_1}$ such that 
$w(\varpi) = \varpi - \sum d_k \a_{i_k}$ where $d_k=1$ for all $k$ but $k=c_1$,
and $d_{c_1} = 2$. Thus $\scal{\varpi-w(\varpi),\rho^\vee} = l(w) + 1$ and the
lemma is again true. The adjoint case follows by duality.
\end{proo}

\vskip .2cm

\begin{defi}
Let $\pi \in \eta(\cM)$ and $\g$ be a root. We say that $\g$ interacts with
$\pi$ if $\scal{\pi,\g^\vee} > 0$ and either:
\begin{itemize}
\item
We are in the minuscule or cominuscule cases and $\g$ is a simple root.
\item
We are in the coadjoint or adjoint cases, $c_1$ does not
divide $l(\pi)+1$, and $\g$ is a simple root.
\item
We are in the coadjoint or adjoint cases and $c_1$
divides $l(\pi)+1$; this implies that there exist a simple root $\a$ and an
integer $d$ such that $\pi = \a + d\delta$. In this case $\g$
interacts with $\pi$ if either $\g = \a$ or $\g = \a+ \beta$,
where $\beta$ is another simple root.
\end{itemize}
\end{defi}

We are now in a position to give a more explicit formula for the quantum
Chevalley formula of Proposition \ref{prop-fulton}:
\begin{theo}
\label{theo-chevalley}
We have a formula
$$
h * \s_\pi = \sum_{\g} \scal{\pi,\g^\vee}
\s_{s_\g(\pi)}\ ,
$$
where the sum runs over all positive roots $\g$ which interact with $\pi$.
\end{theo}
\begin{proo}
Assume first we are in the minuscule case.
Let $\g>0$ be such that $l(s_\g(\pi)) = l(\pi)+1$. Then by Lemma
\ref{lemm-length} we have
$\scal{\varpi-s_\g(\pi),\rho^\vee} =
\scal{\varpi-\pi+\g,\rho^\vee} = l(\pi)+1 = \scal{\varpi-\pi,\rho^\vee} + 1$,
thus $\scal{\g,\rho^\vee}=1$ and $\g$ is a simple root. In the cominuscule
case the same argument gives that $\g^\vee$ is a simple root, so $\g$
is simple. The coadjoint and adjoint cases are similar, except that
one has to take care whether $l(s_\g(\pi))$ is divisible by $c_1$ or not.
\end{proo}

\subsection{Affine symmetries in the adjoint and coadjoint cases}
\label{subsection-affine}

The description of the quantum monomials in terms of roots gives a
nice interpretation of the affine symmetries described in
\cite{affine}. Let us recall the description of these symmetries, we
refer to \cite{affine} for more details. Let $c$ be a cominuscule
vertex of the Dynkin diagram and let $\a_c$ be the corresponding
cominuscule simple root. Denote by $\tau_c$ the orientation preserving
automorphism of the extended Dynkin diagram sending the added vertex
to the vertex $c$ (there is a unique such automorphism). Denote by
$\varpi_c^\vee$ the cominuscule coweight associated to $c$ and by
$v_c$ the shortest element such that
$v_c\varpi_c^\vee=w_0\varpi_c^\vee$ with $w_0$ the longest element of
the Weyl group. We proved in \cite{affine} the following formula
$$\s_{v_c}*\s_u=q^{\eta_P(\varpi^\vee_c-u^{-1}(\varpi^\vee_c))}\s_{v_cu}$$ 
where $\eta_P:Q^\vee\to Q^\vee/Q_P^\vee$ is the projection from the coroot
lattice to its quotient by the coroot lattice of $P$ and where we
identify this quotient with $\Z$.

\begin{prop}
(\i) We have the equality 
$$\s_{v_c}=\left\{
\begin{array}{ll}
\s_{-\a_c-\Theta+\t}&\textrm{in the coadjoint case}\\  
\s_{-\a_c}&\textrm{in the adjoint case.}\\
\end{array}\right.$$
In particular
$\deg(\s_{v_c})=c_1$.

(\i\i) For any real long root $\a$ in the adjoint case (resp. short
root $\a$ in the coadjoint case), we have the formula
$$\s_{v_c}*\s_\a=\s_{\tau_c(\a)-\delta}=q\s_{\tau_c(\a)}.$$
\end{prop}

\begin{proo}
Point ($\imath$) is an easy computation. For ($\imath\imath$), let $u\in W^P$
such that $u(\varpi)=\a$. Using the fact that
$\tau_c=v_ct_{-\varpi_c^\vee}$ we get
$\tau_c(\a)=v_c(\a)+\sca{\a}{\varpi_c^\vee}\delta$. But we check the
equality
$\eta_P(\varpi^\vee_c-u^{-1}(\varpi^\vee_c))=1-\sca{\a}{\varpi_c^\vee}$
and the result follows.
\end{proo}

\section{Maximal degree for quantum cohomology}

\subsection{Dimension argument}

In this section we want to bound the possible degrees of the quantum
monomials in the quantum parameter $q$ that may appear in a quantum
product $\s_u*\s_v$. The first easy restriction is given by the
following:

\begin{lemm}
\label{lemm-dmax3}
Let $X$ be a adjoint or coadjoint homogeneous space. Let
$u$ and $v$ in $W^P$, then if we write
$$\s_u*\s_v=\sum_{d,w}c_{u,v}^w(d)q^d\s_w$$ 
we have $c_{u,v}^w(d)=0$ unless $d\leq 3$.
\end{lemm}

\begin{proo}
  This is a simple dimension count. Recall that $c_{u,v}^w(d)=0$
  unless $l(u)+l(v)+\dim X-l(w)=\dim X+c_1(X)d$ thus unless
  $c_1(X)d\leq 2\dim X$. But in the adjoint case, we have $\dim
  X=2c_1(X)-1$ giving $c_{u,v}^w(d)=0$ unless $d\leq 3$. 
  In the coadjoint non adjoint case,
  there are three cases. 
For the Lagrangian Grassmannian, we have
  $\dim \G_\omega(2,2n)=4n-5=2c_1(\G_\omega(2,2n))-3$.
  Thus we deduce that $c_{u,v}^w(d) = 0$
  if $d>3$. For $F_4/P_4$ we have $\dim(F_4/P_4) = 15$ and
  $c_1(F_4/P_4) = 11$, thus  
  $c_{u,v}^w(d)=0$ unless $d\leq 2$. For $G_2/P_1$, this variety being
  a quadric, the result is well known (see for example \cite{cmp}).
\end{proo}

In the next subsections we show that $c_{u,v}^w(d) = 0$ if $d>2$.
In the coadjoint non adjoint case, we are in fact more interested in
understanding the geometric reason of this vanishing, since the proof
of the above Lemma almost 
proves this point (the only case left is that of $\G_\omega(2,2n)$). 
Rather what we will do is to describe in a uniform way the locus $Y_3$
covered by degree 3 curves passing through 2 general points (see Definition
\ref{yd} below) for adjoint and coadjoint varieties and explain that the
vanishing of the quantum coefficients is related to a porism. 
We believe that it is worth it describing the rich geometry of the
subvarieties $Y_3$.

\subsection{Rational curves on homogeneous spaces}

For the following definition, we shall assume that $X$ is a
homogeneous space with Picard group ${\rm Pic}(X)=\Z$ and we pick the only 
ample line bundle ${\cal L}$ which is a generator of the Picard group (we may 
also include the adjoint variety for $SL_n$ by taking as ample line bundle 
${\cal L}$ the line bundle given by its minimal embedding). All degrees are 
considered with respect to ${\cal L}$.

\begin{defi}
\label{yd}
Let $d$ be an integer and $x$ and $y$ be points in $X$. 

(\i) Let $X_d(x)$ be the locus of points in $X$ connected to $x$ 
by a rational curve of degree $d$. In symbols:
$X_d(x)=\{z\in X\ /\ \textrm{there exists $C$ a rational degree $d$ 
curve through $x$ with $z\in C$}\}.$

(\i\i) Let $Y_d(x,y)$ be the locus of points in $X$ connected to $x$ 
and $y$ by the same rational curve of degree $d$. In symbols:
$Y_d(x,y)=\{z\in X\ /\ \textrm{there exists $C$ a rational degree $d$ 
curve 
with $x,y,z\in C$}\}.$
%
\end{defi}

We observed in \cite{CMP1} that, for $X$ a minuscule or cominuscule rational 
homogeneous space, the geometry of rational curves and especially the 
geometry of the varieties $X_d(x)$ and $Y_d(x,y)$ is intimately related 
to the quantum cohomology of $X$. We shall see that this is still the
case for adjoint and coadjoint varieties.

Let us denote by $V=V(\varpi)$ the irreducible representation of
highest weight $\varpi$. For $\pi$ a weight of the representation
$V$, we denote by $x_\pi$ a vector of weight $\pi$. All points in
$X$ are equivalent and we therefore will often take $x_\varpi$, a
highest weight vector of $V$, as general point in $X$ (recall that $X$
is the orbit of $[x_\varpi]$, the class of $x_\varpi$ in $\p V$).

We introduce the notation $M_{d,3}(X)$ to denote the moduli 
space of degree $d$ morphism from $\pu$ to $X$ with three marked points 
on $\pu$. Recall from \cite{thomsen} or \cite{perrin} that $M_{d,3}(X)$ is 
irreducible smooth and has the expected dimension $\dim(X)+dc_1(X)$. In 
particular, we may speak of general degree $d$ rational curves on $X$.
We shall denote by $ev_i:M_{d,3}(X)\to X$ the evaluation map defined by the 
$i$-th marked point.

\begin{prop}      \label{prop_irreductibilite}
Let $x,y$ be two generic points on a generic curve of degree $d$,
then $X_d(x)$ and $Y_d(x,y)$ are irreducible.
\end{prop}

\begin{proo}
We know that $M_{d,3}(X)$ is irreducible. We have a $G$-equivariant 
evaluation $ev_1 : M_{d,3}(X) \to X$. As $X$ is covered by rational 
curves of any degree, its image is $X$ itself. Assume that $x=x_\varpi$.
Let $U$ be a unipotent subgroup of $G$ such that $U \to X$ is an isomorphism 
from $U$ to an open neighbourhood of $x$, which we denote by $O$.
We have an isomorphism
$$
\begin{array}{rcl}
ev_1^{-1}(x_\varpi) \times U  & \rightarrow & ev_1^{-1}(O) \\
(f,u)  &  \mapsto  &  u \circ f
\end{array}
$$
Since $ev_1^{-1}(O)$, an open subset of $M_{d,3}(X)$, is irreducible,
$ev_1^{-1}(x)$
must be also irreducible. Therefore its image under
$ev_2$ which is $X_d(x_\varpi)$, is irreducible.

For $Y_d(x,y)$ we argue essentially in the same way. We may assume that 
$x=x_\varpi$. We note that $X_d(x_\varpi)$ is a Schubert variety since 
it is irreducible from the above and $B$-stable since $x_\varpi$ is fixed by 
$B$. Because $y$ and the degree $d$ curve through $x_\varpi$ are supposed to 
be general, we may assume that $y$ is in the open $B$-orbit of $X_d(x_\varpi)$.
Thus there is again a unipotent subgroup $U$ of $B$ such that $U \to X$ 
induces an isomorphism between $U$ and an open neighbourhood of $y$ in 
$X_d(x_\varpi)$. Since $U$ is included in $B$ it preserves 
$ev_1^{-1}(x_\varpi)$, and thus induces an isomorphism
$ev_1^{-1}(x_\varpi) \simeq (ev_1 \times ev_2)^{-1}(x_\varpi,y) \times U$.
Therefore $(ev_1 \times ev_2)^{-1}(x_\varpi,y)$ is irreducible and we 
conclude as before.
\end{proo}

Note that the proof shows that $Y_d(x,y)$ does not depend, up to isomorphism,
on the choice of the two points $x,y$ on a generic curve of degree $d$. Thus
we can define the following:

\begin{defi}
\label{defi_x(2)} Let $X_d(X)$ (resp. $Y_d(X)$) be the isomorphism
class of the varieties $X_d(x)$ (resp. $Y_d(x,y)$) for one (resp.
two) generic point(s) $x$ (resp. $x$ and $y$) on a generic degree $d$ curve.
When no confusion on the variety $X$ can be made we simplify $X_d(X)$ and 
$Y_d(X)$ in $X_d$ and $Y_d$.
\end{defi}


\begin{coro}
The variety $X_d$ in $X$ is a Schubert variety stable under
$P$.
\end{coro}

\begin{defi}
(\i) We shall denote by $x_\varpi(d)$ the $T$-stable point in the open 
$B$-orbit  of the Schubert variety $X_d(x_\varpi)$.

(\i\i) For $d,n$ integers we denote by $\delta_{X}(d,n)$ the dimension of 
the space of degree $d$ curves passing through $n$ generic points of a 
generic curve of degree $d$ in $X$.
\end{defi}

\begin{exam}
  The value of $\delta_{\p^k}(1,n)$ is $2(k-1)$ if $n=0$, it is 
$k-1$ if $n=1$, and it is 0 otherwise. 
\end{exam}

We now prove a formula involving 
the dimension of $X_d(x)$ and of $Y_d(x,y)$. 

\begin{prop}      \label{prop_dim_yd}
Let $x,y$ be two generic points on a generic degree $d$ curve, then we 
have 
$$\dim(X_d(x)) + \dim(Y_d(x,y)) + \delta(d,3) = c_1d.$$
\end{prop}


\begin{proo}
The moduli space $M_{d,3}(G/P)$ has dimension
$\dim(X) + dc_1$. Therefore, as $X$ is covered by rational curves of any 
degree, the inverse image
$ev_1^{-1}(x_\varpi)$ has dimension
$dc_1$. Considering the morphism
$ev_2 : ev_1^{-1}(x_\varpi) \rightarrow X$, with image $X_d(x_\varpi)$, 
we deduce the following equality of dimensions: 
$\dim((ev_1\times ev_2)^{-1}(x_\varpi,x_\varpi(d)) = c_1 d - \dim(X_d)$.

Now the image of $ev_3 : (ev_1\times ev_2)^{-1}(x_\varpi,x_\varpi(d))
\rightarrow X$ is by definition $Y_d(x,y)$, and the fibers of this 
morphism are $\delta(d,3)$-dimensional.
\end{proo}

We recall here how to describe chains of $T$-invariant rational curves 
(see for example \cite{FW}): 

\begin{prop}      \label{prop_chaine}
The chains of $k$ irreducible $T$-invariant rational curves starting from 
$x_\varpi$ of degrees $(d_i)_{i\in[1,k]}$ are in bijection
with sequences $(\alpha_i)_{i\in[1,k]}$
of positive roots with $\scal{s_{\alpha_{m-1}} \cdots 
s_{\alpha_1}(\varpi) , \alpha_m^\vee}
= d_m $ for all $m$.
The end-point of the $m$-th component of such a chain is the projectivisation
of the weight line of weight
$s_{\alpha_m} \cdots s_{\alpha_1}(\varpi) =
\varpi - d_1 \alpha_1 - \cdots - d_m \alpha_m$.
\end{prop}

We shall need the following:

\begin{lemm}      \label{lemm_une_chaine}
Let $x,y \in X$ be two $T$-invariant points and $d$ an integer. Assume
that there is exactly one
$T$-invariant chain of degree $d$ through $x$ and $y$.
Then there is only one chain of degree $d$ through $x$ and $y$.
\end{lemm}
\begin{proo}
Assume that there is a chain which is not $T$-invariant. Then the $T$-orbit
of this chain provides a positive-dimensional family of chains of degree
$d$ through $x$ and $y$. In the 
moduli space of curves, the
corresponding positive dimensional variety must have at least 2 fixed points,
a contraction.
\end{proo}

\subsection{Adjoint case}
\label{adjointes}

We first investigate the geometry of conics in adjoint homogeneous
spaces so we fix an adjoint homogeneous space $X$.

\begin{prop}       \label{prop_y2_adjoint}
We have $X_2 = X$ and $Y_2$ is a conic.
\end{prop}
\begin{proo}
Recall that $\Theta$ denotes the highest root. Since
$\langle \Theta , \Theta \rangle = 2$, the curve through $x_\Theta$
with tangent
direction $-\Theta$ has degree 2; its other $T$-fixed point has weight
$s_\Theta ( \Theta ) = - \Theta$, thus it is $x_{- \Theta}$, that is, the
$T$-fixed point representing the dense orbit. Therefore $X_2 = X$.

To understand $Y_2$, we determine the degree 2 curves through
$x_\Theta$ and $x_{-\Theta}$. By Proposition \ref{prop_chaine}, such a
curve is either given by a positive root $\alpha$ such that
$\scal{\Theta , \alpha} = 2$ and $\Theta - 2 \alpha = - \Theta$ (in which
case $\alpha = \Theta$: this corresponds to the previous conic), or by a
sequence $(\beta , \alpha)$ of two positive roots such that
$\scal{\Theta , \alpha} = \scal {s_\alpha ( \Theta ) , \beta} = 1$
and $\Theta - \alpha - \beta = - \Theta$. Since $\Theta$ is the highest
root, this last equation implies $\alpha = \beta = \Theta$, contradicting
the degrees.

Thus there is only one $T$-invariant degree 2 curve through
$x_\Theta$ and $x_{-\Theta}$. By Lemma \ref{lemm_une_chaine} there is only
one degree 2 curve through
$x_\Theta$ and $x_{-\Theta}$ and this implies that $Y_2$ is a conic.
\end{proo}

\begin{defi}
\label{defi-yp3}
We denote $Y'_3 \subset X$ the union of all the lines in $X$ which
meet $Y_2$. It is also the union of all the \emph{reducible} curves of
degree 3 through two general points $x$ and $y$.
\end{defi}

\begin{prop}
\label{prop_degre3_adjoint}
Let $x,y,z$ be three points in $X$. If there
is a curve of degree 3 through $x,y,z$, then there are infinitely many.
Thus $Y_3 = Y'_3$.
\end{prop}
\begin{proo}
Since the dimension of the space of curves through 3 points is an
upper-continuous function on the 3 points, we may assume that $x$ and $y$
are generic, thus that $x=x_\Theta$ and $y=x_{-\Theta}$.

The proposition amounts to the fact that $\delta(3,3) = 1$, or, using
proposition \ref{prop_dim_yd}, that $\dim Y_3 = c_1$. Let $C$ be the conic
through $x$ and $y$ and let $Y'_3$ be the variety covered by the lines 
meeting $C$. Since the dimension of the variety of lines through
any point in $X$ is $c_1 - 2$, the dimension of $Y'_3$ is at most
$c_1$. It is in fact $c_1$ since $T_x X$ and $T_y X$ only meet along
the projectivisation of the line $[\fg_\Theta , \fg_{-\Theta}]$ in $\fg$.

To prove that $Y_3(x_\Theta,x_{-\Theta}) = Y'_3$, we use a suitable
decomposition of $\fg$. Let $L \subset G$ be the pointwise
stabiliser of $C$ and $\fl$ its Lie algebra. If $\a$ is a root
of $\fl$ then $\scal{\Theta,\a^\vee}=0$. Let $\fsl_2(\Theta)$ be
generated by $\fg_\Theta , \fg_{-\Theta}$ and $[\fg_\Theta ,
\fg_{-\Theta}]$. We consider $\fg$ as a module over $\fsl_2 \times
\fl$. This module certainly contains $\fsl_2 \times \fl$. The tangent space
$T_xX$ of $X$ at $x$ (resp. $y$) can be identified with the $T$-stable 
subspace of $\fg$ with weights the negative roots $\a$ such that 
$\scal{\Theta,\a^\vee}<0$ (resp. the positive roots $\a$ such that 
$\scal{\Theta,\a^\vee}>0$). We define the subspace $T_x^1X$ of $T_xX$ 
(resp. $T_x^1X$ of $T_xX$) as the $T$-stable subspace with weights the 
negative roots $\a$ such that $\scal{\Theta,\a^\vee}=-1$ (resp. the 
positive roots $\a$ such that $\scal{\Theta,\a^\vee}=1$). Since $L
\subset P$, the space $V_L := T_x^1 X$ is a representation of $L$,
included in $\fg$. Similarly $T_y^1 X$ is also an
$\fl$-representation, and moreover $T_x^1 X$ and $T_y^1 X$
are disjoint. Counting dimensions, we see that we have

$$
\fg = \fsl_2 \times \fl \oplus \C^2 \otimes V_L
$$

Concretely, here is the decomposition in each case~:

$$
\begin{array}{|c|c|}

\hline

\fg & \mbox{Decomposition of }\fg \\

\hline

\fsl_n & (\fgl_2 \times \fgl_{n-2})_0 \oplus
\C^2 \otimes ( \C^{n-2} \oplus (\C^{n-2})^* ) \\

\fso_{2n+1} & \fsl_2 \times (\fsl_2 \times \fso_{2n-3}) \oplus \C^2 \otimes
(\C^2 \otimes \C^{2n-3}) \\

\fsp_{2n} & \fsl_2 \times \fsp_{2n-2} \oplus \C^2 \otimes \C^{2n-2} \\

\fso_{2n} & \fsl_2 \times (\fsl_2 \times \fso_{2n-4}) \oplus \C^2 \otimes
(\C^2 \otimes \C^{2n-4}) \\

\fe_6 & \fsl_2 \times \fsl_6 \oplus \C^2 \otimes \wedge^3 \C^6  \\

\fe_7 & \fsl_2 \times \fso_{12} \oplus \C^2 \otimes V^{32}  \\

\fe_8 & \fsl_2 \times \fe_7 \oplus \C^2 \otimes V^{56}  \\

\ff_4 & \fsl_2 \times \fsp_6 \oplus \C^2 \otimes \wedge^3_\omega \C^6 \\

\fg_2 & \fsl_2 \times \fsl_2 \oplus \C^2 \otimes S^3\C^2  \\

\hline

\end{array}
$$


\medskip

Using this description we can show that $Y_3 = Y'_3$. Let $\nu :
\C^2 \rightarrow \fsl_2$ be the natural map of degree two with image
the affine cone over $\fsl_2^\ad$, the adjoint variety for $\fsl_2$. 
Let $X_L \subset V_L$ be the
affine cone over the closed $L$-orbit in $\p V_L$. Then $Y'_3$ is
described as the set of classes of elements $(s,0,v\otimes \alpha)
\in \fsl_2 \times \fl \oplus \C^2 \otimes V_L$, with $\alpha \in
X_L$, and $s$ and $\nu(v)$ proportional. The points $x,y$ are
respectively $[\nu(1,0):0:0],[\nu(0,1):0:0]$. 

If
$z=[\nu(v),0,v\otimes \alpha] \in Y'_3 \subset X$,
$q:\C^2 \to \C$ is the quadratic form vanishing on $(1,0)$ and $(0,1)$
and such that $q(v)=1$,
$l$ is a linear form such that $l(v)=1$, then the map
$$
\begin{array}{rcl}
\p^1 = \p \C^2 & \rightarrow & Y'_3 \\
\ [w] & \mapsto & [l(w).\nu(w):0:q(w).w\otimes \alpha]
\end{array}
$$
is a rational curve of degree 3 through $x,y$ and $z$. Thus there are 
infinitely
many rational curves of degree 3 through $x,y,z$, so
$Y'_3 = Y_3$ and the proposition is proved.
\end{proo}

\subsection{Coadjoint varieties as hyperplane sections of Scorza varieties}

Scorza varieties were first defined and classified
by F.L. Zak as varieties satisfying an extremal property with respect to 
their higher secants (see \cite{zak}). For our purpose it will be enough 
to use the following classification theorem as a definition:

\begin{defi}
\label{defi-scorza}
A Scorza variety is defined by two integers $\delta$ and $m$, with
$\delta \in \{2,4,8\}$ and $m \geq 3$. If $\delta = 8$ then we must have
$m=3$. The corresponding homogeneous varieties are as follows 
(a generic hyperplane section of a Scorza variety turns out to be itself 
homogeneous; this homogeneous
space is described in the third column):
$$
\begin{array}{c|c|c}
\delta & X(\delta,m) & X(\delta) \cap H \\
\hline
2 & \p^{m-1} \times \p^{m-1} & \F(1,m-1,m) \\
\hline
4 & \G(2,2m) & \G_\omega(2,2m) \\
\hline
8 & E_6/P_1 & F_4/P_4 \\
\end{array}
$$
\end{defi}

The reason why we consider Scorza varieties is that we recover the coadjoint 
non adjoint varieties as hyperplane sections of some Scorza variety (except 
for $G_2/P_1$). Actually F.L. Zak's definition of Scorza varieties allows 
$\delta=1$, but we don't need the varieties $X(1,m)$ and therefore will not 
speak about them.

Those varieties are naturally embedded via the smallest ample divisor $\cO(1)$.
The numbers $\delta,m$ have some geometrical and algebraic meanings: $\delta$
is the secant defect of $X(\delta,m)$, meaning that
$\dim(Sec(X)) = 2\dim(X)+1-\delta$ where $Sec(X)$ denotes the union of
the secant lines to $X$.
More precisely
through two generic points of $X$ there passes a unique $\delta$-dimensional 
quadric. Moreover $m$ is the smallest integer such that the
the $m$-th higher secant of $X$ is the whole projective space \cite{zak}.

On the other hand, the vector space in the projectivisation of which 
$X(\delta,m)$ is embedded has a natural structure of a Jordan algebra, 
the integer $m$ is the generic rank of this algebra and $\delta$ the 
dimension of the composition algebra which coordinatises 
this algebra 
\cite{jacobson}.

%

\begin{lemm}
\label{lemm-c1-scorza}
We have
$c_1(X) = \frac{\delta}{2}m$ and
$\dim(X) = \delta (m-1)$.
\end{lemm}
\begin{proo}
The following numbers are well-known 
$$
\begin{array}{c|c|c|c}
X & \p^{m-1} \times \p^{m-1} & \G(2,2m) & E_6/P_1 \\ \hline
\delta & 2 & 4 & 8 \\ \hline
\dim(X) & 2(m-1) & 2(2m-2) & 16 \\ \hline
c_1(X) & m & 2m & 12  \\
\end{array}
$$
and yield the lemma.\end{proo}

\subsection{A description of $Y'_3$ for Scorza varieties}

In this subsection, we fix $X$ a Scorza variety with parameters
$(\delta,m)$.
We first describe $Y_2$.

\begin{prop}
\label{prop-y2-scorza}
Through 2 generic points in $X$ there passes a conic, and $Y_2$ is
a smooth quadric of dimension $\delta$.
\end{prop}
\begin{proo}
Let $x,y \in X$ be generic. By \cite{zak}, through $x$ and $y$
there passes a unique quadric
of dimension $\delta$, that we denote
$Q_{x,y}$. This proves the first claim of the lemma. Moreover we
have $Q_{x,y} \subset Y_2(x,y)$. Since $Y_2(x,y)$ is irreducible
by Proposition \ref{prop_irreductibilite}, and has dimension at most
$2c_1(X) - \dim(X) = \delta$ by Lemma \ref{lemm-c1-scorza} and Proposition 
\ref{prop_dim_yd}, the proposition is proved.
\end{proo}

\begin{defi}
\label{defi-sd}
(\i) Let $d>0$, $e\geq 0$ be integers, $Spin_{2d}$ the corresponding spinor 
group, $A$ resp. $\S$ the vector resp. a spinor representation, with highest
weight vector $a_0$ resp. $s_0$. We denote $S \subset \p \S$ the spinor 
variety. Let $W$ denote the natural representation
of $SL_e$ with highest weight vector $w_0$.

(\i\i) We denote $S(d,e)$ the horospherical variety which is
the closure of the $(Spin_{2d} \times SL_e)$-orbit
of $[a_0 + s_0 \otimes w_0]$
in $\p(A \oplus (S \otimes W))$ (for more on horospherical varieties,
see for example \cite{boris}).
The closed orbit of $[a_0]$, which is a smooth quadric
of dimension $2d-2$, is a closed subvariety of $S(d,e)$ denoted $Q(d)$.
\end{defi}

\begin{lemm}
\label{lemm-sd}
The variety $S(d,e)$ has dimension $\frac{(d+2)(d-1)}2 + e$.
Any point in $S(d,e)-Q(d)$ is contained in a linear space of dimension
$d$ included in $S(d,e)$, and intersecting $Q(d)$ along a maximal
isotropic subspace (of dimension $d-1$).
\end{lemm}
\begin{proo} 
Let $x = [a + s \otimes w]$ be a generic point in $S(d,e)$, 
with $a \in A,s \in \S, w \in W$
and $s \otimes w \not = 0$. Since $x$ is in the $Spin_{2d}\times SL_e$-orbit of
$[a_0+s_0 \otimes w_0]$, the class $[s]$ of $s$ in $\p \S$
it is a point of the spinor variety, thus it corresponds to a
maximal isotropic linear subspace $L_d(s) \subset A$, and $a$ is isotropic.
Moreover $a \in L_d(s)$. Since the stabilisor in $Spin_{2d}$ of an element
$[s] \in S$ acts transitively on
the projectivisation of $L_d(s)$, $S(d,e)$ is the variety of all
classes $[a'+s'\otimes w']$ satisfying
$[s'] \in S$ and $a' \in L_d(s')$. It has dimension $d(d-1)/2 + (d-1) + (e-1) + 1$.
Since $x \in S(d,e)$, $a \in L_d(s)$.
Thus the linear space
$\p(L_d(s) \oplus \C \cdot s \otimes w)$
satisfies the conditions of the lemma.
\end{proo}

Recall that we denote by $\varpi$ the dominant weight corresponding to $X$. 
For $X$ a Scorza variety, remark that $\varpi$ is minuscule. We define two 
subgroups of $G$ as follows. Let $P$ be the 
parabolic subsgroupf of $G$ such that $G/P=X$ and denote by $L$ its Levi 
subgroup. The Dynkin diagram of $L$ can be recovered from the Dynkin diagram 
of $G$ by removing the vertex $v$ corresponding to $P$ and the vertex 
$i(v)$ where $i$ is the Weyl involution.

\begin{defi}
(\i) We denote by $L_1$ the subgroup of $L$ whose Dynkin diagram 
is the union of the connected components of the Dynkin diagram of $L$ non 
orthogonal to $\T$ the highest root of $G$ --- or equivalently connected to 
the added vertex in the extended Dynkin diagram of $G$. 

(\i\i) We denote by $L_2$ the subgroup of $L$ whose Dynkin diagram 
is the union of the connected components of the Dynkin diagram of $L$  
orthogonal to $\T$ the highest root of $G$ --- or equivalently non 
connected to the added vertex in the extended Dynkin diagram of $G$. 

(\i\i\i) We denote by $\tilde L_1$ the subgroup of $G$ generated 
by $L_1$ and $SL_2(\T)$ (the subgroup of $G$ whose Lie algebra is 
$\fg_{-\T}+[\fg_{-\T},\fg_\T]+\fg_\T$). Its Dynkin diagram is the Dynkin 
diagram of $L_1$ together with the added vertex in the extended Dynkin 
diagram of $G$.
\end{defi}

Pictures of the Dynkin diagrams of $G$, $L$, $\tilde L_1$ and $L_2$
for Scorza varieties are given in the proof of the following:

\begin{prop}
\label{prop-y3-scorza}
The variety $Y'_3(E_6/P_1)$ resp. $Y'_3(\G(2,2m))$
is projectively equivalent to $S(5,1)$ resp. $S(3,2m-4)$.
The variety $Y'_3(\p^{m-1} \times \p^{m-1})$ is a union of two
irreducible components equivalent to
$S(2,m-2) \simeq \p^1 \times \p^{m-1}$ and meeting
along a smooth quadric surface.

In all cases $Y'_3(X)$ is the linear section $X \cap \p(T_PSec(X))$,
where $P$ is a generic point in the span of $Y_2$.
\end{prop}
\begin{proo}
Let 
$\a$ be a maximal root with vanishing pairing with $i(\varpi)$.
In the case of $E_6/P_1$ resp. $\G(2,2m),$ resp. $\p^{m-1} \times \p^{m-1}$,
we have
$$\a=\poidsesix 122101 \textrm{ resp. } \a=(1,1,\cdots,1,0,0), 
\textrm{ resp. } \a=\left\{
\begin{array}{l}
  ((1,\cdots,1,0),(0,\cdots,0)) \textrm{ or } \\
( (0,\cdots,0) , (1,\cdots,1,0) ).\\
\end{array}\right.$$ 
Note that we have
$ \scal { \Theta^\vee , \a } = 1 $ thus
$ \scal { -\Theta^\vee , \varpi-\a } = 0 $.
We let $v_1 \in V$ be a weight vector of weight $\varpi$, and
$v_d \in V$ of weight $s_\a(\varpi) = \varpi - \alpha$. It follows
that the line through $v_1$ and $v_d$ is included in $X$, and moreover
we see that $v_1$ is the lowest weight of the vector representation of
$\tilde L_1 \simeq Spin_{\delta+2}$ and $v_d$ is the lowest weight
of the representation $\S \otimes W_e$ where $\S$ is a spin representation
of $\tilde L_1$ and $W_e$ is the natural representation of a simple factor
$SL_e$ of the group $L_2$ (see the following array: the black circles represent
the Dynkin diagram, one has to remove the crossed ones to find $L_1$;
the blue square represent the longest root and the green disks the weight
$\varpi-\alpha$).

$$
\begin{array}{|c|c|c|c|}
\hline
X & \p^{m-1} \times \p^{m-1} & \G(2,2m) & E_6/P_1 \\  \hline
& \begin{pspicture*}(0,0)(3.6,2.4)

\psset{unit=.5}

\psellipse(0.6,3.0)(0.13,0.13)
\psline(0.508,2.908)(0.692,3.092)
\psline(0.508,3.092)(0.692,2.908)
\psellipse(0.6,0.6)(0.13,0.13)
\psline(0.508,0.508)(0.692,0.692)
\psline(0.508,0.692)(0.692,0.508)
\psellipse(1.2,3.0)(0.12,0.12)
\psellipse[linecolor=green,fillstyle=solid,fillcolor=green](1.2,0.6)(0.13,0.13)
\pspolygon[linecolor=green,fillstyle=solid,fillcolor=green](1.919,4.319)(1.681,4.319)(1.681,4.081)(1.919,4.081)
\pspolygon[linecolor=blue,fillstyle=solid,fillcolor=blue](1.919,1.919)(1.681,1.919)(1.681,1.681)(1.919,1.681)
\psellipse(2.4,3.0)(0.13,0.13)
\psellipse(2.4,0.6)(0.13,0.13)
\psellipse(3.0,3.0)(0.13,0.13)
\psline(2.908,2.908)(3.092,3.092)
\psline(2.908,3.092)(3.092,2.908)
\psellipse(3.0,0.6)(0.13,0.13)
\psline(2.908,0.508)(3.092,0.692)
\psline(2.908,0.692)(3.092,0.508)
\psellipse(4.2,3.0)(0.13,0.13)
\psline(4.108,2.908)(4.292,3.092)
\psline(4.108,3.092)(4.292,2.908)
\psellipse(4.2,0.6)(0.13,0.13)
\psline(4.108,0.508)(4.292,0.692)
\psline(4.108,0.692)(4.292,0.508)
\psellipse[linecolor=green,fillstyle=solid,fillcolor=green](4.8,3.0)(0.13,0.13)
\psellipse(4.8,0.6)(0.13,0.13)
\pspolygon[linecolor=blue,fillstyle=solid,fillcolor=blue](5.519,4.319)(5.281,4.319)(5.281,4.081)(5.519,4.081)
\pspolygon[linecolor=green,fillstyle=solid,fillcolor=green](5.519,1.919)(5.281,1.919)(5.281,1.681)(5.519,1.681)
\psellipse(6.0,3.0)(0.13,0.13)
\psellipse(6.0,0.6)(0.13,0.13)
\psellipse(6.6,3.0)(0.13,0.13)
\psline(6.508,2.908)(6.692,3.092)
\psline(6.508,3.092)(6.692,2.908)
\psellipse(6.6,0.6)(0.13,0.13)
\psline(6.508,0.508)(6.692,0.692)
\psline(6.508,0.692)(6.692,0.508)
\psline(0.692,3.092)(1.681,4.081)
\psline(1.919,4.081)(2.908,3.092)
\psline(4.292,3.092)(5.281,4.081)
\psline(5.519,4.081)(6.508,3.092)
\psline(6.47,3.0)(6.13,3.0)
\psline(4.67,3.0)(4.33,3.0)
\psline(2.87,3.0)(2.53,3.0)
\psline(1.08,3.0)(0.73,3.0)
\psline(0.692,0.692)(1.681,1.681)
\psline(1.919,1.681)(2.908,0.692)
\psline(2.87,0.6)(2.53,0.6)
\psline(1.07,0.6)(0.73,0.6)
\psline(4.67,0.6)(4.33,0.6)
\psline(6.13,0.6)(6.47,0.6)
\psline(6.508,0.692)(5.519,1.681)
\psline(5.281,1.681)(4.292,0.692)
\psline[linestyle=dashed](1.32,3.0)(2.27,3.0)
\psline[linestyle=dashed](4.93,3.0)(5.87,3.0)
\psline[linestyle=dashed](1.33,0.6)(2.27,0.6)
\psline[linestyle=dashed](4.93,0.6)(5.87,0.6)

\psset{unit=1}

\end{pspicture*}
 & \begin{pspicture*}(0,0)(7.2000003,2.4)

\psellipse[linecolor=green,fillstyle=solid,fillcolor=green](0.6,0.6)(0.12,0.12)
\psellipse(1.2,0.6)(0.12,0.12)
\psline(1.115,0.515)(1.285,0.685)
\psline(1.115,0.685)(1.285,0.515)
\psellipse(1.8,0.6)(0.12,0.12)
\psellipse(2.4,0.6)(0.12,0.12)
\pspolygon[linecolor=blue,fillstyle=solid,fillcolor=blue](3.719,1.919)(3.481,1.919)(3.481,1.681)(3.719,1.681)
\psellipse(4.8,0.6)(0.12,0.12)
\psellipse[linecolor=green,fillstyle=solid,fillcolor=green](5.4,0.6)(0.12,0.12)
\psellipse(6.0,0.6)(0.12,0.12)
\psline(5.915,0.515)(6.085,0.685)
\psline(5.915,0.685)(6.085,0.515)
\psellipse(6.6,0.6)(0.12,0.12)
\psline(0.711,0.645)(3.481,1.752)
\psline(3.719,1.752)(6.489,0.645)
\psline(6.48,0.6)(6.12,0.6)
\psline(5.88,0.6)(5.52,0.6)
\psline(5.28,0.6)(4.92,0.6)
\psline(2.28,0.6)(1.92,0.6)
\psline(1.68,0.6)(1.32,0.6)
\psline(1.08,0.6)(0.72,0.6)
\psline[linestyle=dashed](2.52,0.6)(4.68,0.6)

\end{pspicture*}
 
& \begin{pspicture*}(0,0)(3.6000001,2.4)

\psellipse(0.6,1.8)(0.122,0.122)
\psline(0.514,1.714)(0.686,1.886)
\psline(0.514,1.886)(0.686,1.714)
\psellipse[linecolor=green,fillstyle=solid,fillcolor=green](1.2,1.8)(0.122,0.122)
\psellipse(1.8,1.8)(0.12,0.12)
\psellipse(1.8,1.2)(0.122,0.122)
\pspolygon[linecolor=blue,fillstyle=solid,fillcolor=blue](1.915,0.715)(1.685,0.715)(1.685,0.485)(1.915,0.485)
\psellipse(2.4,1.8)(0.12,0.12)
\psellipse(3.0,1.8)(0.122,0.122)
\psline(2.914,1.714)(3.086,1.886)
\psline(2.914,1.886)(3.086,1.714)
\psline(0.722,1.8)(1.078,1.8)
\psline(1.322,1.8)(1.68,1.8)
\psline(1.92,1.8)(2.28,1.8)
\psline(2.52,1.8)(2.878,1.8)
\psline(1.8,1.68)(1.8,1.322)
\psline(1.8,1.078)(1.8,0.715)

\end{pspicture*}
 \\  \hline
d & 2 & 3 & 5 \\  \hline
e & m-2 & 2m-4 & 1 \\  \hline
\end{array}
$$

It follows that $Y'_3$ contains the closure $S$ of the
$(\tilde L_1 \times L_2)$-orbit of $[v_1 + v_d]$, which is by definition
equivalent to a $S(d,e)$. The equality $Y'_3 = S$
follows counting dimensions. In fact, observe from the above array that
we always have $e = \frac{d-1}2(2m-d)-1$ (we have no better way to
check this than making the computation in the three cases
$d=2,3,5$). This may be rewritten as
\begin{equation}
\label{equa-dim}
\frac{(d+4)(d-1)}2 + e-1 = (d-1)m + 2d-4.
\end{equation}

We may find lines in $S$ which meet $Y_2$ by choosing
any point in $S$ and then a line in the $\p^d$ which contains it and
was described in Lemma \ref{lemm-sd}, taking into account the fact that in this
way each line is counted $\infty^1$ times. 
The dimension of the set of these
lines is thus the left hand-side of (\ref{equa-dim}) since
$\dim S = \frac{(d+2)(d-1)}2 + e$ by Lemma
\ref{lemm-sd}. On the other hand the variety of lines in $X$ passing
through a point in $Y_2$ has dimension
$\dim\ Y_2 + c_1(X) - 2$, which is the right hand-side of (\ref{equa-dim})
by Lemma \ref{lemm-c1-scorza} since $\delta = 2d-2$.
Thus it follows that any line meeting
$Y_2$ is included in $S$, so that $S=Y'_3$.

\vskip .5cm

Moreover we show that this is also $X \cap \p(T_P Sec(X))$ by decomposing
$V$ as an $(\tilde L_1 \times L_2)$-module (we include a proof for completeness,
but this result will not be used in the sequel). Note that we have a decomposition
$V = V_1 \oplus V_{1/2} \oplus V_0$ with $V_1,V_{1/2},V_0$ as in the following
array:

$$
\begin{array}{c|ccc}
\delta & 2 & 4 & 8 \\ \hline
V & A \otimes B\ (\dim A = \dim B = m) &
\wedge^2 A\ ( \dim A = 2m ) & V_{27} \\
\mbox{with} & A = A_2 \oplus A_{2m-2} , B = B_2 \oplus B_{2m-2} &
A = A_4 \oplus A_{2m-4} &  \\
\tilde L_1 & SL(A_2) \times SL(B_2) & SL(A_4) & Spin_{10} \\
L_2 & SL(A_{m-2}) \times SL(B_{m-2}) & SL(A_{2m-4}) & \{1\} \\
V_1 & A_2 \otimes B_2 & \wedge^2 A_4 & V_{10} \\
V_{1/2} & A_2 \otimes B_{2m-2} \oplus B_2 \otimes A_{2m-2} &
A_4 \otimes A_{2m-4} & \S_{16} \\
V_0 & A_{m-2} \otimes B_{m-2} & \wedge^2 A_{2m-4} & \C \\
\end{array} 
$$

Assuming $P$ is an idempotent
in the Jordan algebra $V$ (which we can always do up to an action of $G$),
the above decomposition is the Pierce decomposition
$V = V_1 \oplus V_{1/2} \oplus V_0$, where $V_x$ is the eigenspace for 
multiplication
by $P$ with eigenvalue $x$.

Note that the linear span of $Y_2$ is $V_1$ and that
$T_PSec(x) = V_1 \oplus V_{1/2}$.
Let $y \in Y'_3$: by definition this means that there exists $x \in Y_2$ 
such that
the line $(xy)$ is included in $X$. Therefore $y \in T_xX$.
By Terracini lemma we have $Y_xX \subset T_PSec(X)$, so that we deduce
$Y'_3 \subset T_PSec(X)$.
In the case $\delta = 2$ resp. $\delta = 4$, $Y'_3$ is clearly
the set of couples $(x,y) \in \p A \times \p B$ with $x \in \p A_2$
or $y \in \p B_2$ resp. the set of subspaces $V_2$ meeting non trivially
$A_4$. Thus this is indeed $X \cap \p(T_P Sec(X))$.
In the case of $E_6/P_1$, $Y'_3$ has dimension 15 and must therefore
equal the hyperplane section $\p(T_PSec(X)) \cap X$.
\end{proo}

\subsection{$Y'_3$ for a hyperplane section of a Scorza variety}

Let $X$ be a Scorza variety and $W$ a general hyperplane section of $X$. 
We will reduce the computation
of $Y_3(W)$ to the case when $X = \p^1 \times \p^2$
(which is not a Scorza variety), and therefore we start studying this
case.

\begin{lemm}
\label{lemm-p1p2}
Let $X = \p^1 \times \p^2$ and $W$ a generic hyperplane section of $X$.
Given 3 points $x_1,x_2,x_3$ in $W$ and a reducible curve $R$ of degree 3
through these
points, $R$ is a limit of some irreducible curves passing through
$x_1,x_2,x_3$ and contained in $W$. 
\end{lemm}
Thus the same result holds for a generic curve in $X$.
\begin{proo}
First of all, let us describe $W$ as a scroll.
Restricting the linear form $H$ on the $\p^1$'s of
$\p^1 \times \p^2$ yields a
linear map $\C^3 \to \C^2$,  with 1-dimensional kernel. Let $x \in \p^2$ correspond
to this kernel: it is the only element in $\p^2$ such that
$\p^1 \times \{x\}$ is contained in $H$.
Take a line $l$ in $\p^2$ not containing $x$ and let
$C := (\p^1 \times l) \cap H$. Then $C$ is a conic which may be described as the
set of elements $(a,f(a))$ for some isomorphism $f:\p^1 \to l$. Then
$W$ is the union of all the lines joining $(a,f(a)) \in C$
and $(a,x) \in \p^1 \times \{x\}$.

Now, given three generic
points $x_1,x_2,x_3$ on
$W$, we may choose $C$ such that $x_1,x_2 \in C$, so that
$x_i = (a_i,f(a_i)) \in \p^1 \times l$
for some $a_i \in \p^1$ and $i=1,2$. Then the third point is
on the line joining a point $(a_3,f(a_3))$ in $C$ and $(a_3,x) \in \p^1 \times \{x\}$.
Thus we may assume that $R = C \cup l$, where $l$ is the line between
$(a_3,f(a_3))$ and $(a_3,x)$.

If
$h$ resp. $q$ is a linear resp. quadratic form on $\p^1$ then the
morphism 
$$t \mapsto (t,h(t)f(t)+q(t)x) \in W \subset \p^1 \times \p^2$$
will pass through $x_1,x_2$ provided that
$q$ vanishes at $a_1,a_2$, and imposing moreover that
it passes through $x_3$ yields one more relation between $h$ and $q$. Thus
there are indeed infinitely many irreducible
curves of degree 3 through $x_1,x_2,x_3$.

If $q$ tends to 0 with the vanishing point of $h$ tending accordingly to
$a_3$, we see that the limit of those curves must contain $C$. Since it is 
connected and contains
$x_3$, this limit can only be $R$.
\end{proo}

Let $X$ be a Scorza variety and $W = X \cap H$ a generic hyperplane
section of $X$.

\begin{prop}
\label{deg3-coadj}
We have $Y_3(X) = Y'_3(X)$ and $Y_3(W) = Y_3(X) \cap H = Y'_3(W)$,
and if there passes a curve of degree
3 through 3 points of $W$, then there passes infinitely many.
\end{prop}
\begin{proo}
We show that any reducible curve $R$ in $X$ or $W$ of degree 3
through $x_1,x_2,x_3$ is a limit of some irreducible curves passing
through the same points, which implies the proposition.

We first consider the case of $X$.
To see this, we use the description of $Y'_3(X)$ given in Proposition
\ref{prop-y3-scorza}. By an action of $G$, we may assume that $x_1$ resp.
$x_2$ is the projectivisation of the highest resp. lowest weight line in $V$, so that
$x_1=a_1 , x_2=a_2$ are some points
on the quadric $Q(d)$, and $x_3$ is on the line joining
a point $a_3$ in $Q(d)$ and a point $b_3 = s_3 \otimes w_3$ in
$S \times \p W$, with $a_3 \in L_d(s_3)$, and finally that $R$ is the
union of this line
and a conic $C$ included in $Q(d)$ and passing through $a_1$ and $a_2$. 
Let $V_{d-2} = L_d(s_3) \cap \scal{a_1,a_2,a_3}^\bot$, an isotropic
subspace of dimension $d-2$. We choose a non-degenerate subspace $V_4$
of dimension 4 included in $V_{d-2}^\bot$,
containing $a_1,a_2,a_3$, and disjoint from $V_{d-2}$.

Since $\p V_4 \cap Q(d)$ is a smooth quadric surface,
to each point $a \in \p V_4 \cap Q(d)$ corresponds an isotropic plane
that we denote $V_a$, such that $\tilde V_a := V_{d-2} \oplus V_a$ is a maximal
isotropic subspace in $V_{2d}$ parametrised by a point $s_a$ in $S$.
Remark that any two distinct subspaces $\tilde V_a , \tilde V_b$
meet in codimension 2 (namely,
along $V_{d-2}$), so that the curve traced out by the $s_a$'s in the 
spinor variety is in fact a line.

It follows that $X$ contains a $\p^1 \times \p^2$, namely the union of
all the planes $\p (V_a \oplus \C.s_a \otimes w_3)$ for
$a$ in the quadric surface. Moreover, since $\tilde V_{a_3} = L_d(s_3)$,
this
$\p^1 \times \p^2$ contains the line joining $a_3$ and $b_3$.
Since $\p V_4 \cap Q(d)$ contains $a_1,a_2,a_3$, it contains $C$ and thus
our $\p^1 \times \p^2$ contains our initial curve $R$. Thus by Lemma
\ref{lemm-p1p2}, $R$ can be deformed to an irreducible curve passing
through $x_1,x_2,x_3$.

If we assume from the beginning that $R$ is included in $W = X \cap H$,
then again
from Lemma \ref{lemm-p1p2}, the same conclusion holds.
\end{proo}

\begin{theo}         \label{coro_dmax_adjoint}
Let $X$ be an adjoint or a coadjoint variety and let $\s_u$ and $\s_v$ be
classical cohomology classes. Then if we write
$$\s_u*\s_v=\sum_{d,w}c_{u,v}^w(d)q^d\s_w$$ 
we have $c_{u,v}^w(d)=0$ unless $d\leq 2$.
\end{theo}
\begin{proo}
First, by Proposition \ref{prop_degre3_adjoint} and Proposition
\ref{deg3-coadj}, we know that $c_{u,v}^w(3) = 0$. Lemma
\ref{lemm-dmax3} concludes the proof for $c_{u,v}^w(d)$ with $d>3$.
\end{proo}

\begin{coro}
\label{gw-point-au-carre}
  Let $X$ be an adjoint variety, then in the quantum cohomology
  $QH^*(X,\C)$ we have the formula
$${\rm pt}*{\rm pt}=2q^2\ell$$ 
where ${\rm pt}$ resp. $\ell$ denotes the class of a point resp. the
class of a line in $X$.
\end{coro}

\begin{proo}
  By degree counting the first quantum terms are of degree 2 in $q$. 
Furthermore by the previous result, there are no higher powers of $q$. 
We thus have to prove the equality $I_2({\rm pt},{\rm pt},h)=2$ where $h$ 
is the hyperplane class. As we have seen in Proposition 
\ref{prop_y2_adjoint}, there is a unique conic through two general points 
in $X$. As this conic meets a hyperplane in two points, the result follows.
\end{proo}

\section{Quantum to classical principle for degree one Gromov-Witten 
invariants}  
\label{subsection_computing}

In this section, we explain how to apply by now classical technics 
initiated by A. Buch, A. Kresch and H. Tamvakis \cite{BKT} to compute 
degree one Gromov-Witten invariants. We shall here use the presentation 
of \cite{CMP1}. Most of the arguments are very similar to those of 
\cite[Section 3]{CMP1} and we will therefore not reproduce them. 

Let $F$ be the variety of lines on $X$. This variety is well described in 
\cite{LM}. In particular, note that this variety is homogeneous for $X$ 
adjoint but is not homogeneous for $X$ coadjoint and not adjoint. 
In the last case, the comparison principle is less useful to compute 
Gromov-Witten invariants on $X$ than to compute intersections on $F$.

Let us recall some notation from \cite{CMP1}. Points in $F$ will be 
denoted by $\ell$ and we shall denote by $L$ the associated line 
in $X$. We have an incidence $I=\{(x,\ell)\in X\times F\ /\ x\in L\}$ 
between $X$ and $F$. We denote by $p$ resp. $q$ the projection from $I$ 
to $X$ resp. $F$. For $w\in W^P$, we denote the variety $q(p^{-1}(X(w)))$ 
by $F(\wh)$.

The following lemma is an adaptation of Lemma 3.7 in \cite{CMP1} to the 
non minuscule case. For $w\in W^P$, we denote by $w^\vee$ the element 
of $W^P$ such that $X(w)$ and $X(w^\vee)$ are Poincar{\'e} dual.

\begin{lemm}
\label{intersection}
  Let $X(w)$ and $X(v^\vee)$ be two Schubert varieties such that
  $X(w)\subset X(v^\vee)$. Then there exists an element $g\in G$ such that
  the intersection $X(w)\cap g\cdot X(v)$ is a reduced point.
\end{lemm}

\begin{proo}
There always exists a sequence of inclusions of Schubert varieties
$$X(w)=X(w_0)\subset X(w_1)\subset\cdots\subset X(w_r)=X(v^\vee)$$
where $X(w_i)$ is a moving divisor in $X(w_{i+1})$ for all $i$. We prove, 
by descending induction on $i\in[0,r]$, that there exists an element $g_i$ 
in $G$ such that the intersection $X(w_i)\cap g_i\cdot X(v)$ is a reduced 
point. For $i=r$ this is Poincar{\'e} duality. Assume that the result holds 
for $w_{i+1}$ i.e. there exists an element $g_{i+1}\in G$ such that 
$X(w_{i+1})\cap g_{i+1}\cdot X(v)$ is a reduced point say $x$. Recall that 
because $X(w_i)$ is a moving divisor in $X(w_{i+1})$, we have the equality:
$$\displaystyle{X(w_{i+1})=\bigcup_{h\in {\rm Stab}(X(w_{i+1}))}h\cdot
  X(w_i).}$$
We deduce from this that there exists an element $h\in{\rm
  Stab}(X(w_{i+1}))$ such that $x\in h\cdot X(w_i)$. But now $h\cdot X(w_i)$
  meets $g_{i+1}\cdot X(v)$ in $x$ at least. Since $h\cdot X(w_i)$ is 
contained in $X(w_{i+1})$ and $X(w_{i+1})\cap g_{i+1}\cdot X(v)$ is the 
reduced point $x$, we have that $h\cdot X(w_i)$ and $g_{i+1}\cdot
  X(v)$ meet transversely and only in $x$.
\end{proo}


For the next two results, the proofs of 
\cite[Lemma 3.10 and Corollary 3.11]{CMP1} apply verbatim.

\begin{lemm}
Let $X(u)$, $X(v)$ and $X(w)$ be three proper Schubert subvarieties 
of $X$,  such that 
$$\codim(X(u))+\codim(X(v))+\codim(X(w)) = \dim(X)+c_1(X).$$
Then for $g$, $g'$ and $g''$ three general elements in $G$, the
intersection $g\cdot F(\uh)\cap g'\cdot F(\vh)\cap g''\cdot F(\wh)$ is
a finite set of reduced points. 

Let $\ell$ be a point in this
intersection, then the line $L$ meets each of $g\cdot X(u)$, $g'\cdot X(v)$ and
$g''\cdot X(w)$ in a unique point and these points are in general
position in $L$.
\end{lemm}

\begin{coro}
\label{XF_degre1}
Let $X(u)$, $X(v)$ and $X(w)$ be three proper Schubert subvarieties of
$X$. Suppose that the sum of their codimensions is $\dim(X)+c_1(X)$.  
Then
$$I_1([X(u)],[X(v)],[X(w)])=I_0([F(\uh)],[F(\vh)],[F(\wh)]).$$
\end{coro}

As in \cite{CMP1}, we may use this result to compute combinatorially some
Gromov-Witten invariants. Let us simply give an example of this.

\begin{exam}
\label{calcul-de-q}
   We assume that $X$ is of adjoint type and homogeneous under a group of 
type different from $A_n$. Let $u$ and $v$ be elements in $W^P$ such that 
$l(u)+l(v)=c_1(X)=\deg(q)$. We want to compute the quantum product 
$\s^u*\s^v$. 

Let us first remark that in this situation, if $c_1=c_1(X)$, then we have 
$\dim X=2c_1-1$. In particular, the only quantum non classical term 
appearing in $\s^u*\s^v$ is of the form $aq$ for some $a\in\Z$. In symbols:
$$\s^u*\s^v=\s^u\cup\s^v+aq.$$
The integer $a$ is equal to the Gromov-Witten invariant 
$I_1([X(u^\vee)],[X(v^\vee)],[{\rm pt}])$. By the previous discussion, 
this invariant is equal to $I_0([q(p^{-1}(X(u^\vee)))],
[q(p^{-1}(X(u^\vee)))],[q(p^{-1}({\rm pt}))])$.

In our situation, the variety $F$ is homogeneous of the form $G/Q$ for 
$Q$ a parabolic subgroup. By  the same discussion as in 
\cite[Subsection 3.3]{CMP1}, the variety $F(\wh)=q(p^{-1}(X(w)))$ for 
$X(w)$ a Schubert variety in $X$ is a Schubert variety in $F$. The 
notation $\wh$ will in that situation stand for the corresponding
element in $W^Q$.  
Furthermore, the incidence variety $I$ is also homogeneous of the form 
$G/R$ with $R=P\cap Q$ and the inverse image $p^{-1}(X(w))$ of a Schubert 
variety is again a Schubert variety. We may obtain the element $\wh$ 
easily from $w$ by the following procedure. 
\begin{itemize}
\item First remark that $P$ is maximal associated to a simple root, say 
$\a$, and that $w_X$ the longest element in $W^P$ has the form
$s_\a\wt_X$ with $l(\wt_X)=l(w_X)-1$. In other words $\wt_X=s_\a w_X$.
\item Let $Z$ denote the fiber $p^{-1}({\rm pt})$ of the zero dimensional 
Schubert variety ${\rm pt}$ in $X$. This is a Schubert variety in $I$ and 
let $w_Z\in W^R$ be the corresponding element. Then $p^{-1}(X(w))$ is the 
Schubert variety in $I$ associated to $ww_Z\in W^R$.
\item If $w\neq w_X$, then we have $\wh=w w_Z\in W^Q$. We have
  $\wh_X=s_\a w_Xw_Z$. This last element is the longest element in
  $W^R$, we denote it by $w_F$. This element is self-inverse thus
  $w_F=w_Zw_X s_\a$.
\end{itemize}

\noindent
We want to compute the Littlewood-Richardson coefficient 
$$I_0([F(\uveeh)],[F(\vveeh)],[F(w_Z)])=
c_{\uveeh^\star,\vveeh^\star}^{w_Z}$$ 
in $F$ with $\star$ denoting Poincar{\'e} duality in $F$. Now compute for
$l(u)>0$: 
$$\uveeh^\star = w_0\uveeh w_0w_F = w_0(u^\vee w_Z) w_0w_F =
w_0(w_0uw_0w_X) w_Z w_0w_F =u w_X w_Z w_F = us_\a.$$
The same is true for $\vveeh^\star$. We thus get
$$a=c_{us_\a,vs_\a}^{w_Z}.$$
Now remark that for $X$ adjoint the fiber $Z$ is always 
a product if minuscule and cominuscule homogeneous spaces. In particular 
of $\Lambda$ is the weight associated to $Q$, then $w_Z$ is 
$\Lambda$-minuscule and we way 
compute the above Littlewood-Richardson coefficient using jeu de
taquin by \cite[Theorem 1.3]{littlewood}.

Let us do an explicit example. Let $X=E_6/P_2$, the heap associated to
$w_Z$ is the heap of a Grassmannian as follows:

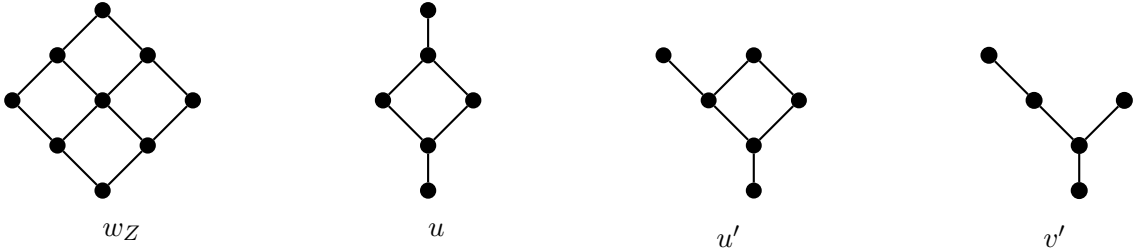
\centerline{\begin{pspicture*}(0,0)(4.8,4.8)

\psellipse[fillstyle=solid,fillcolor=black](1.2,2.4)(0.108,0.108)
\psellipse[fillstyle=solid,fillcolor=black](1.8,3.0)(0.108,0.108)
\psellipse[fillstyle=solid,fillcolor=black](1.8,1.8)(0.108,0.108)
\psellipse[fillstyle=solid,fillcolor=black](2.4,3.6)(0.108,0.108)
\psellipse[fillstyle=solid,fillcolor=black](2.4,2.4)(0.108,0.108)
\psellipse[fillstyle=solid,fillcolor=black](2.4,1.2)(0.108,0.108)
\psellipse[fillstyle=solid,fillcolor=black](3.0,3.0)(0.108,0.108)
\psellipse[fillstyle=solid,fillcolor=black](3.0,1.8)(0.108,0.108)
\psellipse[fillstyle=solid,fillcolor=black](3.6,2.4)(0.108,0.108)
\psline(2.324,3.524)(1.876,3.076)
\psline(1.724,2.924)(1.276,2.476)
\psline(1.276,2.324)(1.724,1.876)
\psline(1.876,1.724)(2.324,1.276)
\psline(2.476,3.524)(2.924,3.076)
\psline(3.076,2.924)(3.524,2.476)
\psline(2.924,2.924)(2.476,2.476)
\psline(2.324,2.324)(1.876,1.876)
\psline(3.524,2.324)(3.076,1.876)
\psline(2.924,1.724)(2.476,1.276)
\psline(1.876,2.924)(2.324,2.476)
\psline(2.476,2.324)(2.924,1.876)
\put(2.4,0.6){$w_Z$}

\end{pspicture*}\begin{pspicture*}(0,0)(3.6,4.8)

\psellipse[fillstyle=solid,fillcolor=black](1.2,2.4)(0.108,0.108)
\psellipse[fillstyle=solid,fillcolor=black](1.8,3.6)(0.108,0.108)
\psellipse[fillstyle=solid,fillcolor=black](1.8,3.0)(0.108,0.108)
\psellipse[fillstyle=solid,fillcolor=black](1.8,1.8)(0.108,0.108)
\psellipse[fillstyle=solid,fillcolor=black](1.8,1.2)(0.108,0.108)
\psellipse[fillstyle=solid,fillcolor=black](2.4,2.4)(0.108,0.108)
\psline(1.8,3.492)(1.8,3.108)
\psline(1.724,2.924)(1.276,2.476)
\psline(1.276,2.324)(1.724,1.876)
\psline(1.8,1.692)(1.8,1.308)
\psline(1.876,2.924)(2.324,2.476)
\psline(2.324,2.324)(1.876,1.876)
\put(1.8,0.6){$u$}

\end{pspicture*}\begin{pspicture*}(0,0)(4.2,4.2)

\psellipse[fillstyle=solid,fillcolor=black](1.2,3.0)(0.106,0.106)
\psellipse[fillstyle=solid,fillcolor=black](1.8,2.4)(0.106,0.106)
\psellipse[fillstyle=solid,fillcolor=black](2.4,3.0)(0.106,0.106)
\psellipse[fillstyle=solid,fillcolor=black](2.4,1.8)(0.106,0.106)
\psellipse[fillstyle=solid,fillcolor=black](2.4,1.2)(0.106,0.106)
\psellipse[fillstyle=solid,fillcolor=black](3.0,2.4)(0.106,0.106)
\psline(1.275,2.925)(1.725,2.475)
\psline(1.875,2.325)(2.325,1.875)
\psline(2.4,1.694)(2.4,1.306)
\psline(2.325,2.925)(1.875,2.475)
\psline(2.475,2.925)(2.925,2.475)
\psline(2.925,2.325)(2.475,1.875)
\uput[l](2.4,0.6){$u'$}

\end{pspicture*}\begin{pspicture*}(0,0)(4.2,4.2)

\psellipse[fillstyle=solid,fillcolor=black](1.2,3.0)(0.113,0.113)
\psellipse[fillstyle=solid,fillcolor=black](1.8,2.4)(0.113,0.113)
\psellipse[fillstyle=solid,fillcolor=black](2.4,1.8)(0.113,0.113)
\psellipse[fillstyle=solid,fillcolor=black](2.4,1.2)(0.113,0.113)
\psellipse[fillstyle=solid,fillcolor=black](3.0,2.4)(0.113,0.113)
\psline(1.28,2.92)(1.72,2.48)
\psline(1.88,2.32)(2.32,1.88)
\psline(2.92,2.32)(2.48,1.88)
\psline(2.4,1.687)(2.4,1.313)
\uput[l](2.4,0.6){$v'$}

\end{pspicture*}}

Let $u=s_{\a_2} s_{\a_4} s_{\a_5} s_{\a_3} s_{\a_4} s_{\a_2}\in
W^{P_2}$ and $v$ be any element in $W^{P_2}$ with $l(v)=5$. We
have $\a=\a_2$ and $us_\a=s_{\a_2} s_{\a_4} s_{\a_5} s_{\a_3}
s_{\a_4}$, thus its heap is not contained in the heap of $w_Z$ and we
get
$$\s^u*\s^v=\s^u\cup\s^v.$$
All the coefficients in this classical product can be computed using
jeu de taquin.

Let us now set $u'=s_{\a_1} s_{\a_4} s_{\a_5} s_{\a_3} s_{\a_4} s_{\a_2}\in
W^{P_2}$ and $v's_{\a_1} s_{\a_5} s_{\a_3} s_{\a_4} s_{\a_2}\in
W^{P_2}$ be elements in $W^{P_2}$. We
have $u's_\a=s_{\a_1} s_{\a_4} s_{\a_5} s_{\a_3} s_{\a_4}$ and
$v's_\a=s_{\a_1} s_{\a_5} s_{\a_3} s_{\a_4}$, thus by jeu de
taquin we get
$$\s^u*\s^v=\s^u\cup\s^v+q.$$
Again, all the coefficients in the classical product $\s^u\cup\s^v$
can be computed using jeu de taquin.
\end{exam}

\section{Presentation of the quantum cohomology rings}
\label{section-pres}


Recall that for the minuscule and cominuscule varieties, presentation of the
quantum cohomology rings are already known: in type $A$ this is done in 
\cite{ST}, in type $B$, $C$ and $D$ in \cite{KT1} and \cite{KT2}. In 
\cite{CMP1} the authors gave a uniform presentation of the quantum
cohomology algebras of the minuscule and cominuscule homogeneous spaces.
In this section we give a presentation of the quantum cohomology algebras
of adjoint and coadjoint homogeneous spaces. One more time for classical 
groups this is already known: in type $A$, a presentation will follow from 
the quantum Chevalley formula (see Section \ref{section-incidence}), in type 
$B$, $C$ and $D$, such a presentation is contained in the more general
results of \cite{BKTnouveau}. Therefore in this section we are
concerned with exceptional types.

Most of the computations here made use of a software written in Java, that
the reader may upload or test at the webpage
\texttt{www.math.sciences.univ-nantes.fr/}\verb|~|\texttt{chaput/quiver-demo.html},
in order
to visualise the quantum Hasse diagrams, apply our version of the
Littlewood-Richardson rule \cite{littlewood}, perform quantum multiplications,
and in general check the assertions in this section.

In the proofs in this section, the Hasse diagrams are colored: the
orange vertices correspond to $\varpi$-minuscule elements. The
blue vertices correspond to quantum Schubert classes with a factor $q$
while grey vertices correspond to quantum Schubert classes with a
factor $q^2$. In the paper \cite{giambelli}, we gathered the quantum
Giambelli formulas for all exceptional adjoint and coadjoint
varieties. We also give in \cite{giambelli} some tableaux used in the
proof of the presentation for $E_8/P_8$ (we did 
not include these tableaux in this paper to save place).

\subsection{Adjoint and coadjoint varieties in type $G_2$}

The homogeneous space $G_2/P_1$ is a smooth 5-dimensional quadric, thus
its quantum cohomology ring is
$QH^*(G_2/P_1,\Q) = \Q[h,q]/(4hq-h^6)$. For $G_2/P_2$ we have the
following:

\begin{prop}
\label{pres-g2p2}
Let $h$ denote the hyperplane class of $G_2/P_2$. We have
$$
QH^*(G_2/P_2) = \Q[h,q]/(h^6 - 18h^3q - 27q^2).
$$
\end{prop}
\begin{proo}
Using Theorem \ref{theo-chevalley}, it is immediate to check that
the quantum Hasse diagram of $G_2/P_2$ begins as follows (and then it
is periodic):\\
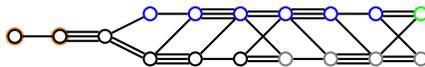
\centerline{

\begin{pspicture*}(0,0)(6.4,1.2)

\psellipse[linecolor=orange](0.2,0.5)(0.12,0.12)
\psellipse(0.2,0.5)(0.1,0.1)
\psellipse[linecolor=orange](0.8,0.5)(0.12,0.12)
\psellipse(0.8,0.5)(0.1,0.1)
\psellipse(1.4,0.5)(0.1,0.1)
\psellipse(2.0,0.2)(0.1,0.1)
\psellipse[linecolor=blue](2.0,0.8)(0.1,0.1)
\psellipse(2.6,0.2)(0.1,0.1)
\psellipse[linecolor=blue](2.6,0.8)(0.1,0.1)
\psellipse(3.2,0.2)(0.1,0.1)
\psellipse[linecolor=blue](3.2,0.8)(0.1,0.1)
\psellipse[linecolor=gray](3.8,0.2)(0.1,0.1)
\psellipse[linecolor=blue](3.8,0.8)(0.1,0.1)
\psellipse[linecolor=gray](4.4,0.2)(0.1,0.1)
\psellipse[linecolor=blue](4.4,0.8)(0.1,0.1)
\psellipse[linecolor=gray](5.0,0.2)(0.1,0.1)
\psellipse[linecolor=blue](5.0,0.8)(0.1,0.1)
\psellipse[linecolor=gray](5.6,0.2)(0.1,0.1)
\psellipse[linecolor=green](5.6,0.8)(0.1,0.1)
\psline(0.3,0.5)(0.7,0.5)
\psline(0.88,0.56)(1.32,0.56)
\psline(0.9,0.5)(1.3,0.5)
\psline(0.88,0.44)(1.32,0.44)
\psline(1.499,0.484)(1.928,0.269)
\psline(1.472,0.431)(1.901,0.216)
\psline(1.489,0.545)(1.911,0.755)
\psline(2.08,0.26)(2.52,0.26)
\psline(2.1,0.2)(2.5,0.2)
\psline(2.08,0.14)(2.52,0.14)
\psline(2.071,0.271)(2.529,0.729)
\psline(2.1,0.8)(2.5,0.8)
\psline(2.7,0.2)(3.1,0.2)
\psline(2.671,0.271)(3.129,0.729)
\psline(2.68,0.86)(3.12,0.86)
\psline(2.7,0.8)(3.1,0.8)
\psline(2.68,0.74)(3.12,0.74)
\psline(3.295,0.23)(3.705,0.23)
\psline(3.295,0.17)(3.705,0.17)
\psline(3.271,0.271)(3.729,0.729)
\psline(3.271,0.729)(3.729,0.271)
\psline(3.295,0.83)(3.705,0.83)
\psline(3.295,0.77)(3.705,0.77)
\psline(3.9,0.2)(4.3,0.2)
\psline(3.871,0.729)(4.329,0.271)
\psline(3.88,0.86)(4.32,0.86)
\psline(3.9,0.8)(4.3,0.8)
\psline(3.88,0.74)(4.32,0.74)
\psline(4.48,0.26)(4.92,0.26)
\psline(4.5,0.2)(4.9,0.2)
\psline(4.48,0.14)(4.92,0.14)
\psline(4.471,0.729)(4.929,0.271)
\psline(4.5,0.8)(4.9,0.8)
\psline(5.095,0.23)(5.505,0.23)
\psline(5.095,0.17)(5.505,0.17)
\psline(5.071,0.271)(5.529,0.729)
\psline(5.071,0.729)(5.529,0.271)
\psline(5.095,0.83)(5.505,0.83)
\psline(5.095,0.77)(5.505,0.77)

\end{pspicture*}}
From this one computes that the degree 2 class is $\s_2=1/3h^2$ and that
the class of a point is $\s_5=1/18h^5-5/6h^2q$. We deduce from this that
$q^2 = (2h\s_5-hq\s_2)/3 = 1/27h^6 - 2/3h^3q$.
\end{proo}

\subsection{Adjoint and coadjoint varieties in type $F_4$}

\begin{prop}
\label{pres-f4p1}
Let $h$ denote the hyperplane class of $F_4/P_1$ and $s \in H^*(F_4/P_1)$
be the Schubert class corresponding to $w=s_2s_3s_2s_1$. We have
$$
QH^*(F_4/P_1,\Q) = 
\Q[h,s,q]/(h^8 - 12s^2 - 16q , 3h^{12} - 18h^8s + 24h^4s^2 + 8s^3).
$$
\end{prop}
\begin{proo}
The quantum Hasse diagram of $F_4/P_1$ begins as follows:\\
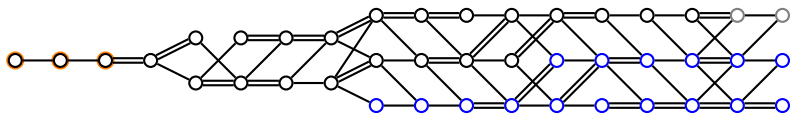
\centerline{

\begin{pspicture*}(0,0)(11.2,1.6)

\psellipse[linecolor=orange](0.2,0.8)(0.12,0.12)
\psellipse(0.2,0.8)(0.1,0.1)
\psellipse[linecolor=orange](0.8,0.8)(0.12,0.12)
\psellipse(0.8,0.8)(0.1,0.1)
\psellipse[linecolor=orange](1.4,0.8)(0.12,0.12)
\psellipse(1.4,0.8)(0.1,0.1)
\psellipse(2.0,0.8)(0.1,0.1)
\psellipse(2.6,0.5)(0.1,0.1)
\psellipse(2.6,1.1)(0.1,0.1)
\psellipse(3.2,0.5)(0.1,0.1)
\psellipse(3.2,1.1)(0.1,0.1)
\psellipse(3.8,0.5)(0.1,0.1)
\psellipse(3.8,1.1)(0.1,0.1)
\psellipse(4.4,0.5)(0.1,0.1)
\psellipse(4.4,1.1)(0.1,0.1)
\psellipse[linecolor=blue](5.0,0.2)(0.1,0.1)
\psellipse(5.0,0.8)(0.1,0.1)
\psellipse(5.0,1.4)(0.1,0.1)
\psellipse[linecolor=blue](5.6,0.2)(0.1,0.1)
\psellipse(5.6,0.8)(0.1,0.1)
\psellipse(5.6,1.4)(0.1,0.1)
\psellipse[linecolor=blue](6.2,0.2)(0.1,0.1)
\psellipse(6.2,0.8)(0.1,0.1)
\psellipse(6.2,1.4)(0.1,0.1)
\psellipse[linecolor=blue](6.8,0.2)(0.1,0.1)
\psellipse(6.8,0.8)(0.1,0.1)
\psellipse(6.8,1.4)(0.1,0.1)
\psellipse[linecolor=blue](7.4,0.2)(0.1,0.1)
\psellipse[linecolor=blue](7.4,0.8)(0.1,0.1)
\psellipse(7.4,1.4)(0.1,0.1)
\psellipse[linecolor=blue](8.0,0.2)(0.1,0.1)
\psellipse[linecolor=blue](8.0,0.8)(0.1,0.1)
\psellipse(8.0,1.4)(0.1,0.1)
\psellipse[linecolor=blue](8.6,0.2)(0.1,0.1)
\psellipse[linecolor=blue](8.6,0.8)(0.1,0.1)
\psellipse(8.6,1.4)(0.1,0.1)
\psellipse[linecolor=blue](9.2,0.2)(0.1,0.1)
\psellipse[linecolor=blue](9.2,0.8)(0.1,0.1)
\psellipse(9.2,1.4)(0.1,0.1)
\psellipse[linecolor=blue](9.8,0.2)(0.1,0.1)
\psellipse[linecolor=blue](9.8,0.8)(0.1,0.1)
\psellipse[linecolor=gray](9.8,1.4)(0.1,0.1)
\psellipse[linecolor=blue](10.4,0.2)(0.1,0.1)
\psellipse[linecolor=blue](10.4,0.8)(0.1,0.1)
\psellipse[linecolor=gray](10.4,1.4)(0.1,0.1)
\psline(0.3,0.8)(0.7,0.8)
\psline(0.9,0.8)(1.3,0.8)
\psline(1.495,0.83)(1.905,0.83)
\psline(1.495,0.77)(1.905,0.77)
\psline(2.089,0.755)(2.511,0.545)
\psline(2.072,0.869)(2.501,1.084)
\psline(2.099,0.816)(2.528,1.031)
\psline(2.695,0.53)(3.105,0.53)
\psline(2.695,0.47)(3.105,0.47)
\psline(2.671,0.571)(3.129,1.029)
\psline(2.671,1.029)(3.129,0.571)
\psline(3.295,0.53)(3.705,0.53)
\psline(3.295,0.47)(3.705,0.47)
\psline(3.271,0.571)(3.729,1.029)
\psline(3.295,1.13)(3.705,1.13)
\psline(3.295,1.07)(3.705,1.07)
\psline(3.9,0.5)(4.3,0.5)
\psline(3.871,0.571)(4.329,1.029)
\psline(3.895,1.13)(4.305,1.13)
\psline(3.895,1.07)(4.305,1.07)
\psline(4.489,0.455)(4.911,0.245)
\psline(4.472,0.569)(4.901,0.784)
\psline(4.499,0.516)(4.928,0.731)
\psline(4.455,0.583)(4.945,1.317)
\psline(4.489,1.055)(4.911,0.845)
\psline(4.472,1.169)(4.901,1.384)
\psline(4.499,1.116)(4.928,1.331)
\psline(5.1,0.2)(5.5,0.2)
\psline(5.071,0.729)(5.529,0.271)
\psline(5.1,0.8)(5.5,0.8)
\psline(5.071,1.329)(5.529,0.871)
\psline(5.095,1.43)(5.505,1.43)
\psline(5.095,1.37)(5.505,1.37)
\psline(5.7,0.2)(6.1,0.2)
\psline(5.671,0.729)(6.129,0.271)
\psline(5.695,0.83)(6.105,0.83)
\psline(5.695,0.77)(6.105,0.77)
\psline(5.671,1.329)(6.129,0.871)
\psline(5.695,1.43)(6.105,1.43)
\psline(5.695,1.37)(6.105,1.37)
\psline(6.295,0.23)(6.705,0.23)
\psline(6.295,0.17)(6.705,0.17)
\psline(6.271,0.729)(6.729,0.271)
\psline(6.3,0.8)(6.7,0.8)
\psline(6.246,0.889)(6.711,1.354)
\psline(6.289,0.846)(6.754,1.311)
\psline(6.3,1.4)(6.7,1.4)
\psline(6.9,0.2)(7.3,0.2)
\psline(6.846,0.289)(7.311,0.754)
\psline(6.889,0.246)(7.354,0.711)
\psline(6.871,0.729)(7.329,0.271)
\psline(6.846,0.889)(7.311,1.354)
\psline(6.889,0.846)(7.354,1.311)
\psline(6.871,1.329)(7.329,0.871)
\psline(6.9,1.4)(7.3,1.4)
\psline(7.5,0.2)(7.9,0.2)
\psline(7.446,0.289)(7.911,0.754)
\psline(7.489,0.246)(7.954,0.711)
\psline(7.5,0.8)(7.9,0.8)
\psline(7.471,1.329)(7.929,0.871)
\psline(7.495,1.43)(7.905,1.43)
\psline(7.495,1.37)(7.905,1.37)
\psline(8.095,0.23)(8.505,0.23)
\psline(8.095,0.17)(8.505,0.17)
\psline(8.071,0.729)(8.529,0.271)
\psline(8.095,0.83)(8.505,0.83)
\psline(8.095,0.77)(8.505,0.77)
\psline(8.071,1.329)(8.529,0.871)
\psline(8.1,1.4)(8.5,1.4)
\psline(8.695,0.23)(9.105,0.23)
\psline(8.695,0.17)(9.105,0.17)
\psline(8.671,0.729)(9.129,0.271)
\psline(8.7,0.8)(9.1,0.8)
\psline(8.671,1.329)(9.129,0.871)
\psline(8.7,1.4)(9.1,1.4)
\psline(9.295,0.23)(9.705,0.23)
\psline(9.295,0.17)(9.705,0.17)
\psline(9.271,0.271)(9.729,0.729)
\psline(9.271,0.729)(9.729,0.271)
\psline(9.295,0.83)(9.705,0.83)
\psline(9.295,0.77)(9.705,0.77)
\psline(9.271,0.871)(9.729,1.329)
\psline(9.271,1.329)(9.729,0.871)
\psline(9.295,1.43)(9.705,1.43)
\psline(9.295,1.37)(9.705,1.37)
\psline(9.895,0.23)(10.305,0.23)
\psline(9.895,0.17)(10.305,0.17)
\psline(9.871,0.271)(10.329,0.729)
\psline(9.9,0.8)(10.3,0.8)
\psline(9.871,0.871)(10.329,1.329)
\psline(9.9,1.4)(10.3,1.4)

\end{pspicture*}}
To compute $s^2$ in $QH^*(F_4/P_1)$,
we first compute the classical product.
Let $\s_{4,2}$ correspond to $s_4s_3s_2s_1$:
the Hasse diagram shows that $\s_{4,2}^2$ has degree 56.
Let $\s_{8,1}$ resp. $\s_{8,2}$ correspond to
$s_1s_2s_3s_2s_4s_3s_2s_1$ resp. $s_2s_1s_3s_2s_4s_3s_2s_1$.
Since $\deg(\s_{8,1})=16$ and $\deg(\s_{8,2})=40$, we deduce
that $\s_{4,2} \cup \s_{4,2} = \s_{8,1} + \s_{8,2}$. Using the Chevalley 
formula it follows that $s \cup s = 6 \s_{8,1} + 8 \s_{8,2}$. The quantum 
product $s*s=s^2$ is then equal to the same class by Example 
\ref{calcul-de-q}. By quantum Chevalley formula, we get the first relation. 
The second relation (of degree 12) is then just a consequence of the 
Chevalley formula.
This is enough in order
to express all the Schubert classes as polynomials in $h,s$ and
$q$, as displayed in \cite{giambelli}. 
\end{proo}

\begin{prop}
\label{pres-f4p4}
Let $h$ denote the hyperplane class of $F_4/P_4$ and $s \in H^*(F_4/P_1)$
be the Schubert class corresponding to $w=s_1s_2s_3s_4$. We have
$$
QH^*(F_4/P_4,\Q) = 
\Q[h,s,q]/(2\ h^{8} - 6\ h^{4}s + 3\ s^{2},
-11\ h^{12} + 26\ h^{8}s + 3\ hq).
$$
\end{prop}
\begin{proo}
The quantum Hasse diagram of $F_4/P_4$ begins as follows:\\
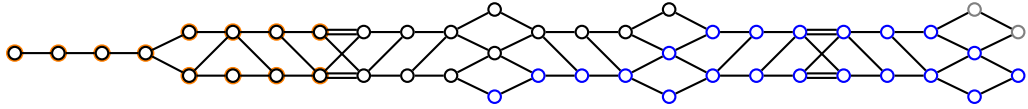
\centerline{

\begin{pspicture*}(0,0)(14.32,1.56)

\psellipse[linecolor=orange](0.2,0.78)(0.116,0.116)
\psellipse(0.2,0.78)(0.097,0.097)
\psellipse[linecolor=orange](0.78,0.78)(0.116,0.116)
\psellipse(0.78,0.78)(0.097,0.097)
\psellipse[linecolor=orange](1.36,0.78)(0.116,0.116)
\psellipse(1.36,0.78)(0.097,0.097)
\psellipse[linecolor=orange](1.94,0.78)(0.116,0.116)
\psellipse(1.94,0.78)(0.097,0.097)
\psellipse[linecolor=orange](2.52,0.48)(0.116,0.116)
\psellipse(2.52,0.48)(0.097,0.097)
\psellipse[linecolor=orange](2.52,1.06)(0.116,0.116)
\psellipse(2.52,1.06)(0.097,0.097)
\psellipse[linecolor=orange](3.1,0.48)(0.116,0.116)
\psellipse(3.1,0.48)(0.097,0.097)
\psellipse[linecolor=orange](3.1,1.06)(0.116,0.116)
\psellipse(3.1,1.06)(0.097,0.097)
\psellipse[linecolor=orange](3.68,0.48)(0.116,0.116)
\psellipse(3.68,0.48)(0.097,0.097)
\psellipse[linecolor=orange](3.68,1.06)(0.116,0.116)
\psellipse(3.68,1.06)(0.097,0.097)
\psellipse[linecolor=orange](4.26,0.48)(0.116,0.116)
\psellipse(4.26,0.48)(0.097,0.097)
\psellipse[linecolor=orange](4.26,1.06)(0.116,0.116)
\psellipse(4.26,1.06)(0.097,0.097)
\psellipse(4.84,0.48)(0.097,0.097)
\psellipse(4.84,1.06)(0.097,0.097)
\psellipse(5.42,0.48)(0.097,0.097)
\psellipse(5.42,1.06)(0.097,0.097)
\psellipse(6.0,0.48)(0.097,0.097)
\psellipse(6.0,1.06)(0.097,0.097)
\psellipse[linecolor=blue](6.58,0.2)(0.097,0.097)
\psellipse(6.58,0.78)(0.097,0.097)
\psellipse(6.58,1.36)(0.097,0.097)
\psellipse[linecolor=blue](7.16,0.48)(0.097,0.097)
\psellipse(7.16,1.06)(0.097,0.097)
\psellipse[linecolor=blue](7.74,0.48)(0.097,0.097)
\psellipse(7.74,1.06)(0.097,0.097)
\psellipse[linecolor=blue](8.32,0.48)(0.097,0.097)
\psellipse(8.32,1.06)(0.097,0.097)
\psellipse[linecolor=blue](8.9,0.2)(0.097,0.097)
\psellipse[linecolor=blue](8.9,0.78)(0.097,0.097)
\psellipse(8.9,1.36)(0.097,0.097)
\psellipse[linecolor=blue](9.48,0.48)(0.097,0.097)
\psellipse[linecolor=blue](9.48,1.06)(0.097,0.097)
\psellipse[linecolor=blue](10.06,0.48)(0.097,0.097)
\psellipse[linecolor=blue](10.06,1.06)(0.097,0.097)
\psellipse[linecolor=blue](10.64,0.48)(0.097,0.097)
\psellipse[linecolor=blue](10.64,1.06)(0.097,0.097)
\psellipse[linecolor=blue](11.22,0.48)(0.097,0.097)
\psellipse[linecolor=blue](11.22,1.06)(0.097,0.097)
\psellipse[linecolor=blue](11.8,0.48)(0.097,0.097)
\psellipse[linecolor=blue](11.8,1.06)(0.097,0.097)
\psellipse[linecolor=blue](12.38,0.48)(0.097,0.097)
\psellipse[linecolor=blue](12.38,1.06)(0.097,0.097)
\psellipse[linecolor=blue](12.96,0.2)(0.097,0.097)
\psellipse[linecolor=blue](12.96,0.78)(0.097,0.097)
\psellipse[linecolor=gray](12.96,1.36)(0.097,0.097)
\psellipse[linecolor=blue](13.54,0.48)(0.097,0.097)
\psellipse[linecolor=gray](13.54,1.06)(0.097,0.097)
\psline(0.297,0.78)(0.683,0.78)
\psline(0.877,0.78)(1.263,0.78)
\psline(1.457,0.78)(1.843,0.78)
\psline(2.026,0.736)(2.434,0.524)
\psline(2.027,0.822)(2.433,1.018)
\psline(2.617,0.48)(3.003,0.48)
\psline(2.588,0.548)(3.032,0.992)
\psline(2.617,1.06)(3.003,1.06)
\psline(3.197,0.48)(3.583,0.48)
\psline(3.168,0.992)(3.612,0.548)
\psline(3.197,1.06)(3.583,1.06)
\psline(3.777,0.48)(4.163,0.48)
\psline(3.748,0.992)(4.192,0.548)
\psline(3.777,1.06)(4.163,1.06)
\psline(4.352,0.51)(4.748,0.51)
\psline(4.352,0.45)(4.748,0.45)
\psline(4.328,0.548)(4.772,0.992)
\psline(4.328,0.992)(4.772,0.548)
\psline(4.352,1.09)(4.748,1.09)
\psline(4.352,1.03)(4.748,1.03)
\psline(4.937,0.48)(5.323,0.48)
\psline(4.908,0.548)(5.352,0.992)
\psline(4.937,1.06)(5.323,1.06)
\psline(5.517,0.48)(5.903,0.48)
\psline(5.488,0.548)(5.932,0.992)
\psline(5.517,1.06)(5.903,1.06)
\psline(6.087,0.438)(6.493,0.242)
\psline(6.086,0.524)(6.494,0.736)
\psline(6.087,1.018)(6.493,0.822)
\psline(6.086,1.104)(6.494,1.316)
\psline(6.667,0.242)(7.073,0.438)
\psline(6.666,0.736)(7.074,0.524)
\psline(6.667,0.822)(7.073,1.018)
\psline(6.666,1.316)(7.074,1.104)
\psline(7.257,0.48)(7.643,0.48)
\psline(7.228,0.992)(7.672,0.548)
\psline(7.257,1.06)(7.643,1.06)
\psline(7.837,0.48)(8.223,0.48)
\psline(7.808,0.992)(8.252,0.548)
\psline(7.837,1.06)(8.223,1.06)
\psline(8.407,0.438)(8.813,0.242)
\psline(8.406,0.524)(8.814,0.736)
\psline(8.407,1.018)(8.813,0.822)
\psline(8.406,1.104)(8.814,1.316)
\psline(8.987,0.242)(9.393,0.438)
\psline(8.986,0.736)(9.394,0.524)
\psline(8.987,0.822)(9.393,1.018)
\psline(8.986,1.316)(9.394,1.104)
\psline(9.577,0.48)(9.963,0.48)
\psline(9.548,0.548)(9.992,0.992)
\psline(9.577,1.06)(9.963,1.06)
\psline(10.157,0.48)(10.543,0.48)
\psline(10.128,0.548)(10.572,0.992)
\psline(10.157,1.06)(10.543,1.06)
\psline(10.732,0.51)(11.128,0.51)
\psline(10.732,0.45)(11.128,0.45)
\psline(10.708,0.548)(11.152,0.992)
\psline(10.708,0.992)(11.152,0.548)
\psline(10.732,1.09)(11.128,1.09)
\psline(10.732,1.03)(11.128,1.03)
\psline(11.317,0.48)(11.703,0.48)
\psline(11.288,0.992)(11.732,0.548)
\psline(11.317,1.06)(11.703,1.06)
\psline(11.897,0.48)(12.283,0.48)
\psline(11.868,0.992)(12.312,0.548)
\psline(11.897,1.06)(12.283,1.06)
\psline(12.467,0.438)(12.873,0.242)
\psline(12.466,0.524)(12.874,0.736)
\psline(12.467,1.018)(12.873,0.822)
\psline(12.466,1.104)(12.874,1.316)
\psline(13.047,0.242)(13.453,0.438)
\psline(13.046,0.736)(13.454,0.524)
\psline(13.047,0.822)(13.453,1.018)
\psline(13.046,1.316)(13.454,1.104)

\end{pspicture*}}\\
%
%
This Hasse diagram shows that the classical (or quantum because its 
cohomological degree is $8<c_1$)
product $s\cup s = s*s$ has degree
$14$. Let us denote $\s_{8,1}$ resp. $\s_{8,2}$ the Schubert classes
corresponding to
$s_3s_2s_4s_1s_3s_2s_3s_4$ resp. $s_4s_3s_2s_1s_3s_2s_3s_4$, of degree
5 resp. 2. We have either $s^2 = 2\s_{8,1}+2\s_{8,2}$ or
$s^2 = 7\s_{8,2}$. Since $s*h^4 \geq s^2$ and
$s*h^4 = 5\s_{8,1}+4\s_{8,2}$, the second case is excluded. Thus we have
$s^2 = 2\s_{8,1}+2\s_{8,2}$. The relations in degree 8 and 12 follow.
\end{proo}

\subsection{Adjoint and coadjoint varieties in type $E_6$}

We now consider the technically more involved cases of groups of type $E$.
In the cohomology of $E_6/P_2$, let $h$ denote the hyperplane class, $s$ resp. $t$
the classes corresponding to the elements $s_3s_4s_2$ resp. $s_1s_3s_4s_2$.

\begin{prop}
\label{pres-e6p2}
The quantum cohomology algebra
$QH^*(E_6/P_2,\Q)$ is the quotient of the polynomial algebra $\Q[h,s,t,q]$
by the relations $h^{8} - 6\ h^{5}s + 3\ h^{4}t + 9\ h^{2}s^{2} - 12\ hst + 6\ t^{2}\ ,\ 
h^{9} - 4\ h^{6}s + 3\ h^{5}t + 3\ h^{3}s^{2} - 6\ h^{2}st + 2\ s^{3}\ ,\ $ and
$-97\ h^{12} + 442\ h^{9}s - 247\ h^{8}t - 507\ h^{6}s^{2} + 624\ h^{5}st - 156\ h^{2}s^{2}t + 48\ hq$.
\end{prop}
\begin{proo}
We know that $QH^*(E_6/P_2,\Q)$ is generated by $h,s$ and $t$ and that there are relations
in degree $8,9,12$.
The quantum Hasse diagram of $E_6/P_2$ begins as follows:

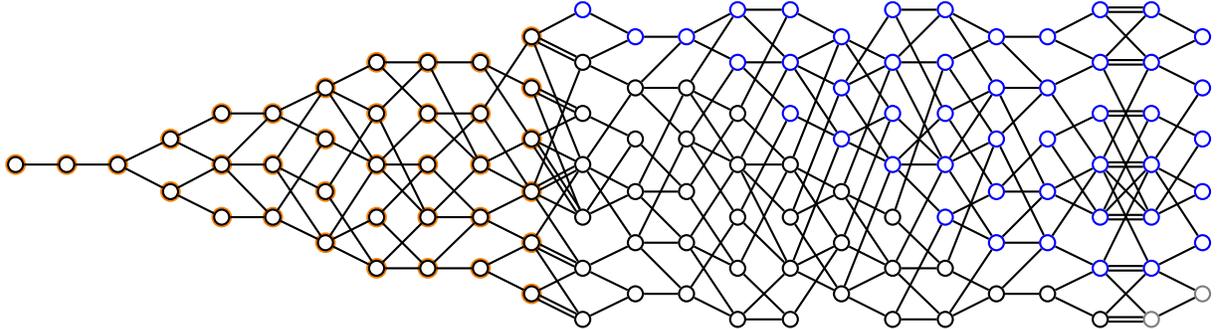
\centerline{\input{hasse-e6p2}}

From this it follows that the degree of $s^2$ is $37752$.
From now on and until the end of this section, we will write Schubert 
classes according to the following notation in the proofs:

\begin{nota}
\label{nota-classes}
To $w \in W/W_P$ we associate the element $\varphi(w)=\a = \varpi - 
w(\varpi)$ of the root
lattice. If $\a = \varphi(w)$ with $w \in W/W_P$, we denote $\sigma(\a)$ the
cohomology class $\s_w$.
\end{nota}

Let us define the following Schubert classes
$$\s_{6,1}=\s\poidsesix112102,\ \s_{6,2}=\s\poidsesix012102,\ 
\s_{6,3}=\s\poidsesix111111,\ \textrm{and}\ \s_{6,4}=\s\poidsesix012111,$$ 
they have respective degrees 10920, 6006, 4992, and 10920.
Considering the following tableaux we deduce from \cite{littlewood}
that $s^2 \geq 2\s_{6,1} + \s_{6,3} + \s_{6,4}$:

\vskip 0.2 cm

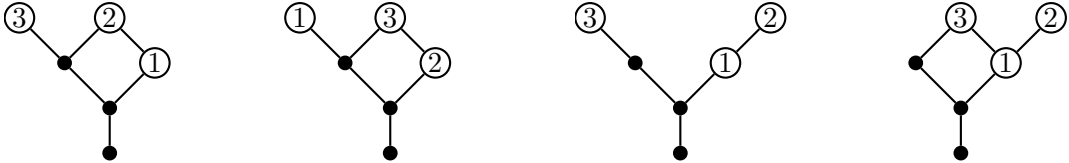
\centerline{

\begin{pspicture*}(0,0)(3.6,2.2)

\psline(0.348,1.851)(0.731,1.469)
\psellipse(0.2,2.0)(0.21,0.21)
\put(0.1,1.876){3}
\psline(1.549,1.851)(1.851,1.549)
\psline(1.252,1.851)(0.869,1.469)
\psellipse(1.4,2.0)(0.21,0.21)
\put(1.3,1.876){2}
\psline(1.851,1.252)(1.469,0.869)
\psellipse(2.0,1.4)(0.21,0.21)
\put(1.9,1.276){1}
\psline(0.869,1.331)(1.331,0.869)
\psellipse[fillstyle=solid,fillcolor=black](0.8,1.4)(0.098,0.098)
\psline(1.4,0.702)(1.4,0.298)
\psellipse[fillstyle=solid,fillcolor=black](1.4,0.8)(0.098,0.098)
\psellipse[fillstyle=solid,fillcolor=black](1.4,0.2)(0.098,0.098)

\end{pspicture*}

\begin{pspicture*}(0,0)(3.6,2.2)

\psline(0.348,1.851)(0.731,1.469)
\psellipse(0.2,2.0)(0.21,0.21)
\put(0.1,1.876){1}
\psline(1.549,1.851)(1.851,1.549)
\psline(1.252,1.851)(0.869,1.469)
\psellipse(1.4,2.0)(0.21,0.21)
\put(1.3,1.876){3}
\psline(1.851,1.252)(1.469,0.869)
\psellipse(2.0,1.4)(0.21,0.21)
\put(1.9,1.276){2}
\psline(0.869,1.331)(1.331,0.869)
\psellipse[fillstyle=solid,fillcolor=black](0.8,1.4)(0.098,0.098)
\psline(1.4,0.702)(1.4,0.298)
\psellipse[fillstyle=solid,fillcolor=black](1.4,0.8)(0.098,0.098)
\psellipse[fillstyle=solid,fillcolor=black](1.4,0.2)(0.098,0.098)

\end{pspicture*}

\begin{pspicture*}(0,0)(3.6,2.2)

\psline(0.348,1.851)(0.731,1.469)
\psellipse(0.2,2.0)(0.21,0.21)
\put(0.1,1.876){3}
\psline(0.869,1.331)(1.331,0.869)
\psellipse[fillstyle=solid,fillcolor=black](0.8,1.4)(0.098,0.098)
\psline(2.451,1.851)(2.148,1.549)
\psellipse(2.6,2.0)(0.21,0.21)
\put(2.5,1.876){2}
\psline(1.851,1.252)(1.469,0.869)
\psellipse(2.0,1.4)(0.21,0.21)
\put(1.9,1.276){1}
\psline(1.4,0.702)(1.4,0.298)
\psellipse[fillstyle=solid,fillcolor=black](1.4,0.8)(0.098,0.098)
\psellipse[fillstyle=solid,fillcolor=black](1.4,0.2)(0.098,0.098)

\end{pspicture*}

\begin{pspicture*}(0,0)(3.6,2.2)

\psline(1.252,1.851)(0.869,1.469)
\psline(1.549,1.851)(1.851,1.549)
\psellipse(1.4,2.0)(0.21,0.21)
\put(1.3,1.876){3}
\psline(2.451,1.851)(2.148,1.549)
\psellipse(2.6,2.0)(0.21,0.21)
\put(2.5,1.876){2}
\psline(0.869,1.331)(1.331,0.869)
\psellipse[fillstyle=solid,fillcolor=black](0.8,1.4)(0.098,0.098)
\psline(1.851,1.252)(1.469,0.869)
\psellipse(2.0,1.4)(0.21,0.21)
\put(1.9,1.276){1}
\psline(1.4,0.702)(1.4,0.298)
\psellipse[fillstyle=solid,fillcolor=black](1.4,0.8)(0.098,0.098)
\psellipse[fillstyle=solid,fillcolor=black](1.4,0.2)(0.098,0.098)

\end{pspicture*}
}

\noindent
On the other hand both $s^2$ and $2\s_{6,1} + \s_{6,3} + \s_{6,4}$ have 
degree 37752: thus we have $s^2 = 2\s_{6,1} + \s_{6,3} + \s_{6,4}$. This 
allows to compute an expression of
all Schubert classes up to degree 6 as polynomials in $h,s,t$, as displayed in
\cite{giambelli}.

Let $\s_{7,i},1\leq i\leq 5$ be the Schubert classes in $H^7(E_6/P_2)$, 
corresponding respectively to the roots
$$\poidsesix112102,\ \poidsesix122101,\ \poidsesix112111,\ \poidsesix012112,
\textrm{ and } \poidsesix012211.$$ 
They have respective degrees 3003, 2925, 4992, 3003, and 2925. The product 
$st$ has degree $7917$. In view of the following
two tableaux it follows that $st = \s_{7,2} + \s_{7,3}$.

\vskip 0.2 cm

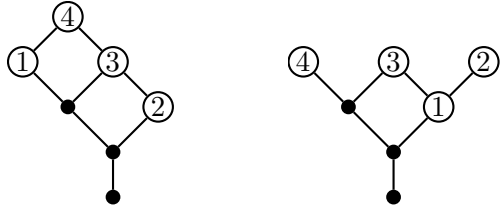
\centerline{

\begin{pspicture*}(0,0)(3.6,2.8)

\psline(0.948,2.451)(1.252,2.148)
\psline(0.651,2.451)(0.348,2.148)
\psellipse(0.8,2.6)(0.21,0.21)
\put(0.7,2.476){4}
\psline(1.549,1.851)(1.851,1.549)
\psline(1.252,1.851)(0.869,1.469)
\psellipse(1.4,2.0)(0.21,0.21)
\put(1.3,1.876){3}
\psline(1.851,1.251)(1.469,0.869)
\psellipse(2.0,1.4)(0.21,0.21)
\put(1.9,1.276){2}
\psline(0.348,1.851)(0.731,1.469)
\psellipse(0.2,2.0)(0.21,0.21)
\put(0.1,1.876){1}
\psline(0.869,1.331)(1.331,0.869)
\psellipse[fillstyle=solid,fillcolor=black](0.8,1.4)(0.098,0.098)
\psline(1.4,0.702)(1.4,0.298)
\psellipse[fillstyle=solid,fillcolor=black](1.4,0.8)(0.098,0.098)
\psellipse[fillstyle=solid,fillcolor=black](1.4,0.2)(0.098,0.098)

\end{pspicture*}

\begin{pspicture*}(0,0)(3.6,2.2)

\psline(0.348,1.851)(0.731,1.469)
\psellipse(0.2,2.0)(0.21,0.21)
\put(0.1,1.876){4}
\psline(1.252,1.851)(0.869,1.469)
\psline(1.549,1.851)(1.851,1.549)
\psellipse(1.4,2.0)(0.21,0.21)
\put(1.3,1.876){3}
\psline(2.451,1.851)(2.148,1.549)
\psellipse(2.6,2.0)(0.21,0.21)
\put(2.5,1.876){2}
\psline(0.869,1.331)(1.331,0.869)
\psellipse[fillstyle=solid,fillcolor=black](0.8,1.4)(0.098,0.098)
\psline(1.851,1.252)(1.469,0.869)
\psellipse(2.0,1.4)(0.21,0.21)
\put(1.9,1.276){1}
\psline(1.4,0.702)(1.4,0.298)
\psellipse[fillstyle=solid,fillcolor=black](1.4,0.8)(0.098,0.098)
\psellipse[fillstyle=solid,fillcolor=black](1.4,0.2)(0.098,0.098)

\end{pspicture*}
}

\noindent
This is enough to express some Giambelli formulas up to degree 9 but to get 
the relation
in degree 8 resp. 9 we must compute $t^2$ resp. $s^3$. To this end we compute
$t^2$ and $s\cdot \s_{6,1}$. The degree of $t^2$ and 
$$\s \poidsesix122111$$ 
are $1638$, and we have the following tableau:

\vskip 0.2 cm

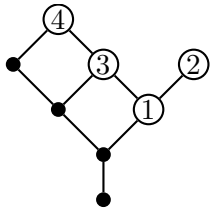
\centerline{

\begin{pspicture*}(0,0)(3.6,2.8)

\psline(0.948,2.451)(1.252,2.148)
\psline(0.651,2.451)(0.269,2.069)
\psellipse(0.8,2.6)(0.21,0.21)
\put(0.7,2.476){4}
\psline(1.252,1.851)(0.869,1.469)
\psline(1.549,1.851)(1.851,1.549)
\psellipse(1.4,2.0)(0.21,0.21)
\put(1.3,1.876){3}
\psline(2.451,1.851)(2.148,1.549)
\psellipse(2.6,2.0)(0.21,0.21)
\put(2.5,1.876){2}
\psline(0.269,1.931)(0.731,1.469)
\psellipse[fillstyle=solid,fillcolor=black](0.2,2.0)(0.098,0.098)
\psline(0.869,1.331)(1.331,0.869)
\psellipse[fillstyle=solid,fillcolor=black](0.8,1.4)(0.098,0.098)
\psline(1.851,1.251)(1.469,0.869)
\psellipse(2.0,1.4)(0.21,0.21)
\put(1.9,1.276){1}
\psline(1.4,0.702)(1.4,0.298)
\psellipse[fillstyle=solid,fillcolor=black](1.4,0.8)(0.098,0.098)
\psellipse[fillstyle=solid,fillcolor=black](1.4,0.2)(0.098,0.098)

\end{pspicture*}
}
\noindent
thus 
$$t^2 = \s \poidsesix122111.$$

The following seven tableaux prove the inequality
$$
s\cdot \s_{6,1} \geq \s \poidsesix123102 + 3 \s \poidsesix122112 + \s
\poidsesix112211 + 2 \s \poidsesix122211:
$$

\centerline{

\begin{pspicture*}(0,0)(3.6,3.4)

\psline(1.4,2.99)(1.4,2.81)
\psline(1.252,3.051)(0.948,2.749)
\psellipse(1.4,3.2)(0.21,0.21)
\put(1.3,3.076){6}
\psline(1.4,2.39)(1.4,2.21)
\psellipse(1.4,2.6)(0.21,0.21)
\put(1.3,2.476){5}
\psline(0.948,2.451)(1.252,2.148)
\psline(0.651,2.451)(0.348,2.148)
\psellipse(0.8,2.6)(0.21,0.21)
\put(0.7,2.476){4}
\psline(1.549,1.851)(1.851,1.549)
\psline(1.252,1.851)(0.869,1.469)
\psellipse(1.4,2.0)(0.21,0.21)
\put(1.3,1.876){3}
\psline(1.851,1.251)(1.469,0.869)
\psellipse(2.0,1.4)(0.21,0.21)
\put(1.9,1.276){2}
\psline(0.348,1.851)(0.731,1.469)
\psellipse(0.2,2.0)(0.21,0.21)
\put(0.1,1.876){1}
\psline(0.869,1.331)(1.331,0.869)
\psellipse[fillstyle=solid,fillcolor=black](0.8,1.4)(0.098,0.098)
\psline(1.4,0.702)(1.4,0.298)
\psellipse[fillstyle=solid,fillcolor=black](1.4,0.8)(0.098,0.098)
\psellipse[fillstyle=solid,fillcolor=black](1.4,0.2)(0.098,0.098)

\end{pspicture*}

\begin{pspicture*}(0,0)(3.6,2.8)

\psline(0.651,2.451)(0.348,2.148)
\psline(0.948,2.451)(1.252,2.148)
\psellipse(0.8,2.6)(0.21,0.21)
\put(0.7,2.476){6}
\psline(1.4,2.39)(1.4,2.21)
\psellipse(1.4,2.6)(0.21,0.21)
\put(1.3,2.476){5}
\psline(2.451,1.851)(2.148,1.549)
\psellipse(2.6,2.0)(0.21,0.21)
\put(2.5,1.876){4}
\psline(0.348,1.851)(0.731,1.469)
\psellipse(0.2,2.0)(0.21,0.21)
\put(0.1,1.876){3}
\psline(1.549,1.851)(1.851,1.549)
\psline(1.252,1.851)(0.869,1.469)
\psellipse(1.4,2.0)(0.21,0.21)
\put(1.3,1.876){2}
\psline(1.851,1.251)(1.469,0.869)
\psellipse(2.0,1.4)(0.21,0.21)
\put(1.9,1.276){1}
\psline(0.869,1.331)(1.331,0.869)
\psellipse[fillstyle=solid,fillcolor=black](0.8,1.4)(0.098,0.098)
\psline(1.4,0.702)(1.4,0.298)
\psellipse[fillstyle=solid,fillcolor=black](1.4,0.8)(0.098,0.098)
\psellipse[fillstyle=solid,fillcolor=black](1.4,0.2)(0.098,0.098)

\end{pspicture*}

\begin{pspicture*}(0,0)(3.6,2.8)

\psline(0.651,2.451)(0.348,2.148)
\psline(0.948,2.451)(1.252,2.148)
\psellipse(0.8,2.6)(0.21,0.21)
\put(0.7,2.476){6}
\psline(1.4,2.39)(1.4,2.21)
\psellipse(1.4,2.6)(0.21,0.21)
\put(1.3,2.476){5}
\psline(2.451,1.851)(2.148,1.549)
\psellipse(2.6,2.0)(0.21,0.21)
\put(2.5,1.876){2}
\psline(0.348,1.851)(0.731,1.469)
\psellipse(0.2,2.0)(0.21,0.21)
\put(0.1,1.876){3}
\psline(1.549,1.851)(1.851,1.549)
\psline(1.252,1.851)(0.869,1.469)
\psellipse(1.4,2.0)(0.21,0.21)
\put(1.3,1.876){4}
\psline(1.851,1.251)(1.469,0.869)
\psellipse(2.0,1.4)(0.21,0.21)
\put(1.9,1.276){1}
\psline(0.869,1.331)(1.331,0.869)
\psellipse[fillstyle=solid,fillcolor=black](0.8,1.4)(0.098,0.098)
\psline(1.4,0.702)(1.4,0.298)
\psellipse[fillstyle=solid,fillcolor=black](1.4,0.8)(0.098,0.098)
\psellipse[fillstyle=solid,fillcolor=black](1.4,0.2)(0.098,0.098)

\end{pspicture*}

\begin{pspicture*}(0,0)(3.6,2.8)

\psline(0.651,2.451)(0.348,2.148)
\psline(0.948,2.451)(1.252,2.148)
\psellipse(0.8,2.6)(0.21,0.21)
\put(0.7,2.476){6}
\psline(1.4,2.39)(1.4,2.21)
\psellipse(1.4,2.6)(0.21,0.21)
\put(1.3,2.476){5}
\psline(2.451,1.851)(2.148,1.549)
\psellipse(2.6,2.0)(0.21,0.21)
\put(2.5,1.876){4}
\psline(0.348,1.851)(0.731,1.469)
\psellipse(0.2,2.0)(0.21,0.21)
\put(0.1,1.876){1}
\psline(1.549,1.851)(1.851,1.549)
\psline(1.252,1.851)(0.869,1.469)
\psellipse(1.4,2.0)(0.21,0.21)
\put(1.3,1.876){3}
\psline(1.851,1.251)(1.469,0.869)
\psellipse(2.0,1.4)(0.21,0.21)
\put(1.9,1.276){2}
\psline(0.869,1.331)(1.331,0.869)
\psellipse[fillstyle=solid,fillcolor=black](0.8,1.4)(0.098,0.098)
\psline(1.4,0.702)(1.4,0.298)
\psellipse[fillstyle=solid,fillcolor=black](1.4,0.8)(0.098,0.098)
\psellipse[fillstyle=solid,fillcolor=black](1.4,0.2)(0.098,0.098)

\end{pspicture*}
}
\centerline{

\begin{pspicture*}(0,0)(3.6,2.8)

\psline(0.348,1.851)(0.731,1.469)
\psellipse(0.2,2.0)(0.21,0.21)
\put(0.1,1.876){6}
\psline(1.4,2.39)(1.4,2.21)
\psellipse(1.4,2.6)(0.21,0.21)
\put(1.3,2.476){5}
\psline(2.148,2.451)(2.451,2.148)
\psline(1.851,2.451)(1.549,2.148)
\psellipse(2.0,2.6)(0.21,0.21)
\put(1.9,2.476){4}
\psline(2.451,1.851)(2.148,1.549)
\psellipse(2.6,2.0)(0.21,0.21)
\put(2.5,1.876){2}
\psline(1.549,1.851)(1.851,1.549)
\psline(1.252,1.851)(0.869,1.469)
\psellipse(1.4,2.0)(0.21,0.21)
\put(1.3,1.876){3}
\psline(1.851,1.251)(1.469,0.869)
\psellipse(2.0,1.4)(0.21,0.21)
\put(1.9,1.276){1}
\psline(0.869,1.331)(1.331,0.869)
\psellipse[fillstyle=solid,fillcolor=black](0.8,1.4)(0.098,0.098)
\psline(1.4,0.702)(1.4,0.298)
\psellipse[fillstyle=solid,fillcolor=black](1.4,0.8)(0.098,0.098)
\psellipse[fillstyle=solid,fillcolor=black](1.4,0.2)(0.098,0.098)

\end{pspicture*}

\begin{pspicture*}(0,0)(3.6,2.8)

\psline(0.651,2.451)(0.348,2.148)
\psline(0.948,2.451)(1.252,2.148)
\psellipse(0.8,2.6)(0.21,0.21)
\put(0.7,2.476){6}
\psline(2.148,2.451)(2.451,2.148)
\psline(1.851,2.451)(1.549,2.148)
\psellipse(2.0,2.6)(0.21,0.21)
\put(1.9,2.476){5}
\psline(2.451,1.851)(2.148,1.549)
\psellipse(2.6,2.0)(0.21,0.21)
\put(2.5,1.876){4}
\psline(0.348,1.851)(0.731,1.469)
\psellipse(0.2,2.0)(0.21,0.21)
\put(0.1,1.876){3}
\psline(1.549,1.851)(1.851,1.549)
\psline(1.252,1.851)(0.869,1.469)
\psellipse(1.4,2.0)(0.21,0.21)
\put(1.3,1.876){2}
\psline(1.851,1.251)(1.469,0.869)
\psellipse(2.0,1.4)(0.21,0.21)
\put(1.9,1.276){1}
\psline(0.869,1.331)(1.331,0.869)
\psellipse[fillstyle=solid,fillcolor=black](0.8,1.4)(0.098,0.098)
\psline(1.4,0.702)(1.4,0.298)
\psellipse[fillstyle=solid,fillcolor=black](1.4,0.8)(0.098,0.098)
\psellipse[fillstyle=solid,fillcolor=black](1.4,0.2)(0.098,0.098)

\end{pspicture*}

\begin{pspicture*}(0,0)(3.6,2.8)

\psline(0.651,2.451)(0.348,2.148)
\psline(0.948,2.451)(1.252,2.148)
\psellipse(0.8,2.6)(0.21,0.21)
\put(0.7,2.476){6}
\psline(2.148,2.451)(2.451,2.148)
\psline(1.851,2.451)(1.549,2.148)
\psellipse(2.0,2.6)(0.21,0.21)
\put(1.9,2.476){5}
\psline(2.451,1.851)(2.148,1.549)
\psellipse(2.6,2.0)(0.21,0.21)
\put(2.5,1.876){4}
\psline(0.348,1.851)(0.731,1.469)
\psellipse(0.2,2.0)(0.21,0.21)
\put(0.1,1.876){1}
\psline(1.549,1.851)(1.851,1.549)
\psline(1.252,1.851)(0.869,1.469)
\psellipse(1.4,2.0)(0.21,0.21)
\put(1.3,1.876){3}
\psline(1.851,1.251)(1.469,0.869)
\psellipse(2.0,1.4)(0.21,0.21)
\put(1.9,1.276){2}
\psline(0.869,1.331)(1.331,0.869)
\psellipse[fillstyle=solid,fillcolor=black](0.8,1.4)(0.098,0.098)
\psline(1.4,0.702)(1.4,0.298)
\psellipse[fillstyle=solid,fillcolor=black](1.4,0.8)(0.098,0.098)
\psellipse[fillstyle=solid,fillcolor=black](1.4,0.2)(0.098,0.098)

\end{pspicture*}
}
\noindent
We leave it to the reader to check that the degrees coincide, so that we 
have equality. 

To express the Schubert classes in $H^{10}(E_6/P_2)$ as polynomials in 
$h,s,t$, 
we must compute $s^2t$, and to this end we compute the product $t\s_{6,1}$. 
According to the following three tableaux and using the same degree 
argument as above this product is 
$$\s \poidsesix123112 + \s \poidsesix122212 + \poidsesix123211.$$

\centerline{

\begin{pspicture*}(0,0)(3.6,3.4)

\psline(1.4,2.99)(1.4,2.81)
\psline(1.252,3.051)(0.948,2.749)
\psellipse(1.4,3.2)(0.21,0.21)
\put(1.3,3.076){4}
\psline(1.4,2.39)(1.4,2.098)
\psellipse(1.4,2.6)(0.21,0.21)
\put(1.3,2.476){3}
\psline(2.451,1.851)(2.069,1.469)
\psellipse(2.6,2.0)(0.21,0.21)
\put(2.5,1.876){2}
\psline(0.948,2.451)(1.331,2.069)
\psline(0.651,2.451)(0.269,2.069)
\psellipse(0.8,2.6)(0.21,0.21)
\put(0.7,2.476){1}
\psline(1.469,1.931)(1.931,1.469)
\psline(1.331,1.931)(0.869,1.469)
\psellipse[fillstyle=solid,fillcolor=black](1.4,2.0)(0.098,0.098)
\psline(1.931,1.331)(1.469,0.869)
\psellipse[fillstyle=solid,fillcolor=black](2.0,1.4)(0.098,0.098)
\psline(0.269,1.931)(0.731,1.469)
\psellipse[fillstyle=solid,fillcolor=black](0.2,2.0)(0.098,0.098)
\psline(0.869,1.331)(1.331,0.869)
\psellipse[fillstyle=solid,fillcolor=black](0.8,1.4)(0.098,0.098)
\psline(1.4,0.702)(1.4,0.298)
\psellipse[fillstyle=solid,fillcolor=black](1.4,0.8)(0.098,0.098)
\psellipse[fillstyle=solid,fillcolor=black](1.4,0.2)(0.098,0.098)

\end{pspicture*}

\begin{pspicture*}(0,0)(3.6,2.8)

\psline(0.651,2.451)(0.269,2.069)
\psline(0.948,2.451)(1.331,2.069)
\psellipse(0.8,2.6)(0.21,0.21)
\put(0.7,2.476){4}
\psline(1.4,2.39)(1.4,2.098)
\psellipse(1.4,2.6)(0.21,0.21)
\put(1.3,2.476){1}
\psline(2.148,2.451)(2.451,2.148)
\psline(1.851,2.451)(1.469,2.069)
\psellipse(2.0,2.6)(0.21,0.21)
\put(1.9,2.476){3}
\psline(2.451,1.851)(2.069,1.469)
\psellipse(2.6,2.0)(0.21,0.21)
\put(2.5,1.876){2}
\psline(0.269,1.931)(0.731,1.469)
\psellipse[fillstyle=solid,fillcolor=black](0.2,2.0)(0.098,0.098)
\psline(1.469,1.931)(1.931,1.469)
\psline(1.331,1.931)(0.869,1.469)
\psellipse[fillstyle=solid,fillcolor=black](1.4,2.0)(0.098,0.098)
\psline(1.931,1.331)(1.469,0.869)
\psellipse[fillstyle=solid,fillcolor=black](2.0,1.4)(0.098,0.098)
\psline(0.869,1.331)(1.331,0.869)
\psellipse[fillstyle=solid,fillcolor=black](0.8,1.4)(0.098,0.098)
\psline(1.4,0.702)(1.4,0.298)
\psellipse[fillstyle=solid,fillcolor=black](1.4,0.8)(0.098,0.098)
\psellipse[fillstyle=solid,fillcolor=black](1.4,0.2)(0.098,0.098)

\end{pspicture*}

\begin{pspicture*}(0,0)(3.6,3.4)

\psline(1.549,3.051)(1.851,2.749)
\psline(1.252,3.051)(0.948,2.749)
\psellipse(1.4,3.2)(0.21,0.21)
\put(1.3,3.076){4}
\psline(2.148,2.451)(2.451,2.148)
\psline(1.851,2.451)(1.469,2.069)
\psellipse(2.0,2.6)(0.21,0.21)
\put(1.9,2.476){3}
\psline(2.451,1.851)(2.069,1.469)
\psellipse(2.6,2.0)(0.21,0.21)
\put(2.5,1.876){2}
\psline(0.948,2.451)(1.331,2.069)
\psline(0.651,2.451)(0.269,2.069)
\psellipse(0.8,2.6)(0.21,0.21)
\put(0.7,2.476){1}
\psline(1.469,1.931)(1.931,1.469)
\psline(1.331,1.931)(0.869,1.469)
\psellipse[fillstyle=solid,fillcolor=black](1.4,2.0)(0.098,0.098)
\psline(1.931,1.331)(1.469,0.869)
\psellipse[fillstyle=solid,fillcolor=black](2.0,1.4)(0.098,0.098)
\psline(0.269,1.931)(0.731,1.469)
\psellipse[fillstyle=solid,fillcolor=black](0.2,2.0)(0.098,0.098)
\psline(0.869,1.331)(1.331,0.869)
\psellipse[fillstyle=solid,fillcolor=black](0.8,1.4)(0.098,0.098)
\psline(1.4,0.702)(1.4,0.298)
\psellipse[fillstyle=solid,fillcolor=black](1.4,0.8)(0.098,0.098)
\psellipse[fillstyle=solid,fillcolor=black](1.4,0.2)(0.098,0.098)

\end{pspicture*}
}

\noindent
These computations are enough to express all the Schubert classes as 
polynomials in
$h,s,t,q$, and the displayed relation in degree $12$ comes for free.
\end{proo}

\subsection{Adjoint and coadjoint varieties in type $E_7$}

In the cohomology of $E_7/P_1$, let $h$ denote the hyperplane class, $s$ 
resp. $t$
the classes corresponding to the elements $s_2s_4s_3s_1$ resp.
$s_7s_6s_5s_4s_3s_1$.

\begin{prop}
\label{pres-e7p1}
The quantum cohomology algebra
$QH^*(E_7/P_1,\Q)$ is the quotient of the polynomial algebra $\Q[h,s,t,q]$
by the relations
$
h^{12} - 6\ h^{8}s - 4\ h^{6}t + 9\ h^{4}s^{2} + 12\ h^{2}st - s^{3} + 
3\ t^{2}\ ,\ 
h^{14} - 6\ h^{10}s - 2\ h^{8}t + 9\ h^{6}s^{2} + 6\ h^{4}st - h^{2}s^{3} 
+ 3\ s^{2}t\ , \ $
and $\ 232\ h^{18} - 1444\ h^{14}s - 456\ h^{12}t + 2508\ h^{10}s^{2} 
+ 1520\ h^{8}st 
- 988\ h^{6}s^{3} + 133\ h^{2}s^{4} + 36\ hq.
$
\end{prop}
\begin{proo}
We know that $QH^*(E_7/P_1,\Q)$ is generated by $h,s$ and $t$ and that 
there are relations in degree $12,14, 18$.
The quantum Hasse diagram of $E_7/P_1$ begins as follows:\\
\centerline{\input{hasse-e7p1}}\\
The strategy of proof is the same as in the previous cases; therefore we 
will only sketch the arguments. We can express the Giambelli formulas up to 
degree 7 without computing any product. In degree 8 we need to compute 
$s^2$. We have
$$
s^2 = \s \poidsesept2221001 + \s \poidsesept1221101 + \s \poidsesept1122101
+ \s \poidsesept1121111.
$$
In fact, computing degrees, this follows from the existence of the following 
four tableaux:

\vskip 0.2 cm

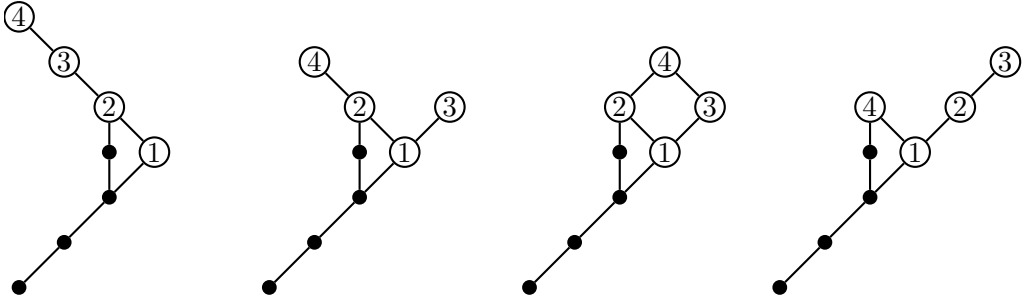
\centerline{

\begin{pspicture*}(0,0)(4.2,4.0)

\psline(0.348,3.651)(0.651,3.349)
\psellipse(0.2,3.8)(0.21,0.21)
\put(0.1,3.676){4}
\psline(0.948,3.051)(1.252,2.749)
\psellipse(0.8,3.2)(0.21,0.21)
\put(0.7,3.076){3}
\psline(1.4,2.39)(1.4,2.098)
\psline(1.549,2.451)(1.851,2.148)
\psellipse(1.4,2.6)(0.21,0.21)
\put(1.3,2.476){2}
\psline(1.4,1.902)(1.4,1.498)
\psellipse[fillstyle=solid,fillcolor=black](1.4,2.0)(0.098,0.098)
\psline(1.851,1.851)(1.469,1.469)
\psellipse(2.0,2.0)(0.21,0.21)
\put(1.9,1.876){1}
\psline(1.331,1.331)(0.869,0.869)
\psellipse[fillstyle=solid,fillcolor=black](1.4,1.4)(0.098,0.098)
\psline(0.731,0.731)(0.269,0.269)
\psellipse[fillstyle=solid,fillcolor=black](0.8,0.8)(0.098,0.098)
\psellipse[fillstyle=solid,fillcolor=black](0.2,0.2)(0.098,0.098)

\end{pspicture*}\hskip -1 cm

\begin{pspicture*}(0,0)(4.2,3.4)

\psline(0.948,3.051)(1.252,2.749)
\psellipse(0.8,3.2)(0.21,0.21)
\put(0.7,3.076){4}
\psline(1.4,2.39)(1.4,2.098)
\psline(1.549,2.451)(1.851,2.148)
\psellipse(1.4,2.6)(0.21,0.21)
\put(1.3,2.476){2}
\psline(1.4,1.902)(1.4,1.498)
\psellipse[fillstyle=solid,fillcolor=black](1.4,2.0)(0.098,0.098)
\psline(2.451,2.451)(2.148,2.148)
\psellipse(2.6,2.6)(0.21,0.21)
\put(2.5,2.476){3}
\psline(1.851,1.851)(1.469,1.469)
\psellipse(2.0,2.0)(0.21,0.21)
\put(1.9,1.876){1}
\psline(1.331,1.331)(0.869,0.869)
\psellipse[fillstyle=solid,fillcolor=black](1.4,1.4)(0.098,0.098)
\psline(0.731,0.731)(0.269,0.269)
\psellipse[fillstyle=solid,fillcolor=black](0.8,0.8)(0.098,0.098)
\psellipse[fillstyle=solid,fillcolor=black](0.2,0.2)(0.098,0.098)

\end{pspicture*}
\hskip -1 cm 

\begin{pspicture*}(0,0)(4.2,3.4)

\psline(2.148,3.051)(2.451,2.749)
\psline(1.851,3.051)(1.549,2.749)
\psellipse(2.0,3.2)(0.21,0.21)
\put(1.9,3.076){4}
\psline(2.451,2.451)(2.148,2.148)
\psellipse(2.6,2.6)(0.21,0.21)
\put(2.5,2.476){3}
\psline(1.4,2.39)(1.4,2.098)
\psline(1.549,2.451)(1.851,2.148)
\psellipse(1.4,2.6)(0.21,0.21)
\put(1.3,2.476){2}
\psline(1.4,1.902)(1.4,1.498)
\psellipse[fillstyle=solid,fillcolor=black](1.4,2.0)(0.098,0.098)
\psline(1.851,1.851)(1.469,1.469)
\psellipse(2.0,2.0)(0.21,0.21)
\put(1.9,1.876){1}
\psline(1.331,1.331)(0.869,0.869)
\psellipse[fillstyle=solid,fillcolor=black](1.4,1.4)(0.098,0.098)
\psline(0.731,0.731)(0.269,0.269)
\psellipse[fillstyle=solid,fillcolor=black](0.8,0.8)(0.098,0.098)
\psellipse[fillstyle=solid,fillcolor=black](0.2,0.2)(0.098,0.098)

\end{pspicture*}\hskip -1 cm 

\begin{pspicture*}(0,0)(4.2,3.4)

\psline(1.4,2.39)(1.4,2.098)
\psline(1.549,2.451)(1.851,2.148)
\psellipse(1.4,2.6)(0.21,0.21)
\put(1.3,2.476){4}
\psline(1.4,1.902)(1.4,1.498)
\psellipse[fillstyle=solid,fillcolor=black](1.4,2.0)(0.098,0.098)
\psline(3.051,3.051)(2.749,2.749)
\psellipse(3.2,3.2)(0.21,0.21)
\put(3.1,3.076){3}
\psline(2.451,2.451)(2.148,2.148)
\psellipse(2.6,2.6)(0.21,0.21)
\put(2.5,2.476){2}
\psline(1.851,1.851)(1.469,1.469)
\psellipse(2.0,2.0)(0.21,0.21)
\put(1.9,1.876){1}
\psline(1.331,1.331)(0.869,0.869)
\psellipse[fillstyle=solid,fillcolor=black](1.4,1.4)(0.098,0.098)
\psline(0.731,0.731)(0.269,0.269)
\psellipse[fillstyle=solid,fillcolor=black](0.8,0.8)(0.098,0.098)
\psellipse[fillstyle=solid,fillcolor=black](0.2,0.2)(0.098,0.098)

\end{pspicture*}
}

\noindent
In degree 10, we have to compute $st$: this is equal to
$\s \poidsesept1222111$, in view of:

\vskip 0.2 cm

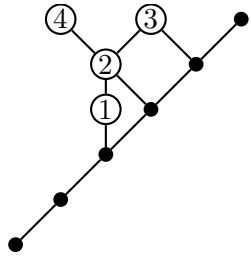
\centerline{

\begin{pspicture*}(0,0)(4.2,3.4)

\psline(0.948,3.051)(1.252,2.749)
\psellipse(0.8,3.2)(0.21,0.21)
\put(0.7,3.076){4}
\psline(1.851,3.051)(1.549,2.749)
\psline(2.148,3.051)(2.531,2.669)
\psellipse(2.0,3.2)(0.21,0.21)
\put(1.9,3.076){3}
\psline(3.131,3.131)(2.669,2.669)
\psellipse[fillstyle=solid,fillcolor=black](3.2,3.2)(0.098,0.098)
\psline(1.4,2.39)(1.4,2.21)
\psline(1.549,2.451)(1.931,2.069)
\psellipse(1.4,2.6)(0.21,0.21)
\put(1.3,2.476){2}
\psline(1.4,1.79)(1.4,1.498)
\psellipse(1.4,2.0)(0.21,0.21)
\put(1.3,1.876){1}
\psline(2.531,2.531)(2.069,2.069)
\psellipse[fillstyle=solid,fillcolor=black](2.6,2.6)(0.098,0.098)
\psline(1.931,1.931)(1.469,1.469)
\psellipse[fillstyle=solid,fillcolor=black](2.0,2.0)(0.098,0.098)
\psline(1.331,1.331)(0.869,0.869)
\psellipse[fillstyle=solid,fillcolor=black](1.4,1.4)(0.098,0.098)
\psline(0.731,0.731)(0.269,0.269)
\psellipse[fillstyle=solid,fillcolor=black](0.8,0.8)(0.098,0.098)
\psellipse[fillstyle=solid,fillcolor=black](0.2,0.2)(0.098,0.098)

\end{pspicture*}
}

\noindent
In degree 12, we have two new monomials: $t^2$ and $s^3$. We have
$$t^2 = \s \poidsesept 2222211.$$ 
To express $s^3$ as a sum of Schubert classes and thus
find the given relation of degree 12, we compute that
$$
s \cdot \s \poidsesept 2221001 = \s \poidsesept 2332101
+ \s \poidsesept 2232111 + \s \poidsesept 2222211
$$
The corresponding tableau for $t^2$ is displayed on the left and the 
three tableaux for the above product 
are on the right:

\vskip 0.2 cm

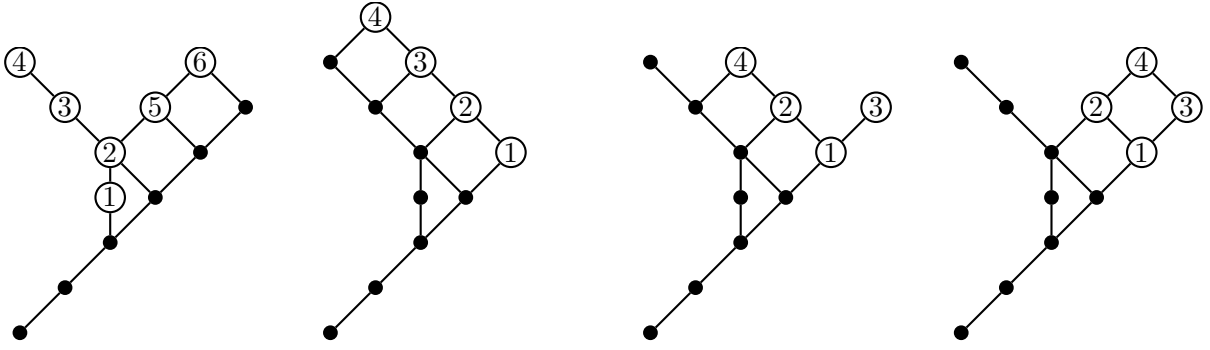
\centerline{

\begin{pspicture*}(0,0)(4.2,4.0)

\psline(0.348,3.651)(0.651,3.349)
\psellipse(0.2,3.8)(0.21,0.21)
\put(0.1,3.676){4}
\psline(0.948,3.051)(1.252,2.749)
\psellipse(0.8,3.2)(0.21,0.21)
\put(0.7,3.076){3}
\psline(2.749,3.651)(3.131,3.269)
\psline(2.451,3.651)(2.148,3.349)
\psellipse(2.6,3.8)(0.21,0.21)
\put(2.5,3.676){6}
\psline(3.131,3.131)(2.669,2.669)
\psellipse[fillstyle=solid,fillcolor=black](3.2,3.2)(0.098,0.098)
\psline(2.148,3.051)(2.531,2.669)
\psline(1.851,3.051)(1.549,2.749)
\psellipse(2.0,3.2)(0.21,0.21)
\put(1.9,3.076){5}
\psline(2.531,2.531)(2.069,2.069)
\psellipse[fillstyle=solid,fillcolor=black](2.6,2.6)(0.098,0.098)
\psline(1.4,2.39)(1.4,2.21)
\psline(1.549,2.451)(1.931,2.069)
\psellipse(1.4,2.6)(0.21,0.21)
\put(1.3,2.476){2}
\psline(1.4,1.79)(1.4,1.498)
\psellipse(1.4,2.0)(0.21,0.21)
\put(1.3,1.876){1}
\psline(1.931,1.931)(1.469,1.469)
\psellipse[fillstyle=solid,fillcolor=black](2.0,2.0)(0.098,0.098)
\psline(1.331,1.331)(0.869,0.869)
\psellipse[fillstyle=solid,fillcolor=black](1.4,1.4)(0.098,0.098)
\psline(0.731,0.731)(0.269,0.269)
\psellipse[fillstyle=solid,fillcolor=black](0.8,0.8)(0.098,0.098)
\psellipse[fillstyle=solid,fillcolor=black](0.2,0.2)(0.098,0.098)

\end{pspicture*}\hskip -0.2 cm

\begin{pspicture*}(0,0)(4.2,4.6)

\psline(0.948,4.251)(1.252,3.948)
\psline(0.651,4.251)(0.269,3.869)
\psellipse(0.8,4.4)(0.21,0.21)
\put(0.7,4.276){4}
\psline(1.549,3.651)(1.851,3.349)
\psline(1.252,3.651)(0.869,3.269)
\psellipse(1.4,3.8)(0.21,0.21)
\put(1.3,3.676){3}
\psline(2.148,3.051)(2.451,2.749)
\psline(1.851,3.051)(1.469,2.669)
\psellipse(2.0,3.2)(0.21,0.21)
\put(1.9,3.076){2}
\psline(2.451,2.451)(2.069,2.069)
\psellipse(2.6,2.6)(0.21,0.21)
\put(2.5,2.476){1}
\psline(0.269,3.731)(0.731,3.269)
\psellipse[fillstyle=solid,fillcolor=black](0.2,3.8)(0.098,0.098)
\psline(0.869,3.131)(1.331,2.669)
\psellipse[fillstyle=solid,fillcolor=black](0.8,3.2)(0.098,0.098)
\psline(1.4,2.502)(1.4,2.098)
\psline(1.469,2.531)(1.931,2.069)
\psellipse[fillstyle=solid,fillcolor=black](1.4,2.6)(0.098,0.098)
\psline(1.4,1.902)(1.4,1.498)
\psellipse[fillstyle=solid,fillcolor=black](1.4,2.0)(0.098,0.098)
\psline(1.931,1.931)(1.469,1.469)
\psellipse[fillstyle=solid,fillcolor=black](2.0,2.0)(0.098,0.098)
\psline(1.331,1.331)(0.869,0.869)
\psellipse[fillstyle=solid,fillcolor=black](1.4,1.4)(0.098,0.098)
\psline(0.731,0.731)(0.269,0.269)
\psellipse[fillstyle=solid,fillcolor=black](0.8,0.8)(0.098,0.098)
\psellipse[fillstyle=solid,fillcolor=black](0.2,0.2)(0.098,0.098)

\end{pspicture*}
\hskip -0.2 cm

\begin{pspicture*}(0,0)(4.2,4.0)

\psline(0.269,3.731)(0.731,3.269)
\psellipse[fillstyle=solid,fillcolor=black](0.2,3.8)(0.098,0.098)
\psline(1.549,3.651)(1.851,3.349)
\psline(1.252,3.651)(0.869,3.269)
\psellipse(1.4,3.8)(0.21,0.21)
\put(1.3,3.676){4}
\psline(1.851,3.051)(1.469,2.669)
\psline(2.148,3.051)(2.451,2.749)
\psellipse(2.0,3.2)(0.21,0.21)
\put(1.9,3.076){2}
\psline(3.051,3.051)(2.749,2.749)
\psellipse(3.2,3.2)(0.21,0.21)
\put(3.1,3.076){3}
\psline(0.869,3.131)(1.331,2.669)
\psellipse[fillstyle=solid,fillcolor=black](0.8,3.2)(0.098,0.098)
\psline(1.4,2.502)(1.4,2.098)
\psline(1.469,2.531)(1.931,2.069)
\psellipse[fillstyle=solid,fillcolor=black](1.4,2.6)(0.098,0.098)
\psline(1.4,1.902)(1.4,1.498)
\psellipse[fillstyle=solid,fillcolor=black](1.4,2.0)(0.098,0.098)
\psline(2.451,2.451)(2.069,2.069)
\psellipse(2.6,2.6)(0.21,0.21)
\put(2.5,2.476){1}
\psline(1.931,1.931)(1.469,1.469)
\psellipse[fillstyle=solid,fillcolor=black](2.0,2.0)(0.098,0.098)
\psline(1.331,1.331)(0.869,0.869)
\psellipse[fillstyle=solid,fillcolor=black](1.4,1.4)(0.098,0.098)
\psline(0.731,0.731)(0.269,0.269)
\psellipse[fillstyle=solid,fillcolor=black](0.8,0.8)(0.098,0.098)
\psellipse[fillstyle=solid,fillcolor=black](0.2,0.2)(0.098,0.098)

\end{pspicture*}\hskip -0.2 cm

\begin{pspicture*}(0,0)(4.2,4.0)

\psline(0.269,3.731)(0.731,3.269)
\psellipse[fillstyle=solid,fillcolor=black](0.2,3.8)(0.098,0.098)
\psline(0.869,3.131)(1.331,2.669)
\psellipse[fillstyle=solid,fillcolor=black](0.8,3.2)(0.098,0.098)
\psline(2.749,3.651)(3.051,3.349)
\psline(2.451,3.651)(2.148,3.349)
\psellipse(2.6,3.8)(0.21,0.21)
\put(2.5,3.676){4}
\psline(3.051,3.051)(2.749,2.749)
\psellipse(3.2,3.2)(0.21,0.21)
\put(3.1,3.076){3}
\psline(2.148,3.051)(2.451,2.749)
\psline(1.851,3.051)(1.469,2.669)
\psellipse(2.0,3.2)(0.21,0.21)
\put(1.9,3.076){2}
\psline(2.451,2.451)(2.069,2.069)
\psellipse(2.6,2.6)(0.21,0.21)
\put(2.5,2.476){1}
\psline(1.4,2.502)(1.4,2.098)
\psline(1.469,2.531)(1.931,2.069)
\psellipse[fillstyle=solid,fillcolor=black](1.4,2.6)(0.098,0.098)
\psline(1.4,1.902)(1.4,1.498)
\psellipse[fillstyle=solid,fillcolor=black](1.4,2.0)(0.098,0.098)
\psline(1.931,1.931)(1.469,1.469)
\psellipse[fillstyle=solid,fillcolor=black](2.0,2.0)(0.098,0.098)
\psline(1.331,1.331)(0.869,0.869)
\psellipse[fillstyle=solid,fillcolor=black](1.4,1.4)(0.098,0.098)
\psline(0.731,0.731)(0.269,0.269)
\psellipse[fillstyle=solid,fillcolor=black](0.8,0.8)(0.098,0.098)
\psellipse[fillstyle=solid,fillcolor=black](0.2,0.2)(0.098,0.098)

\end{pspicture*}
}

\noindent
To get the relation in degree 14, we have to compute the monomial $s^2t$. 
To this end we compute that 
$$t \cdot \s \poidsesept 2221001 = \s \poidsesept 2232212.$$
To express the Schubert
classes in $H^{16}(E_7/P_1)$, we need to compute $st^2$, and to this end 
we show that
$$t \cdot \s \poidsesept 1122211 = \s \poidsesept 2333212.$$ 
These relations come from the two tableaux:

\vskip 0.2 cm

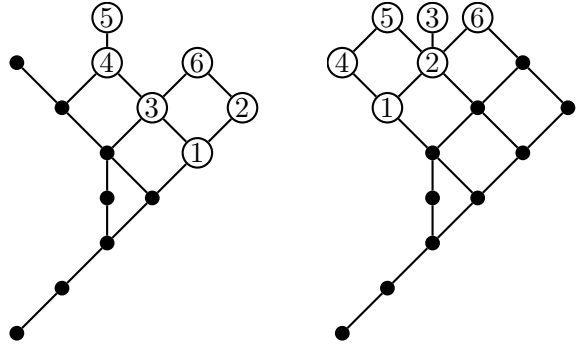
\centerline{

\begin{pspicture*}(0,0)(4.2,4.6)

\psline(0.269,3.731)(0.731,3.269)
\psellipse[fillstyle=solid,fillcolor=black](0.2,3.8)(0.098,0.098)
\psline(1.4,4.19)(1.4,4.01)
\psellipse(1.4,4.4)(0.21,0.21)
\put(1.3,4.276){5}
\psline(2.749,3.651)(3.051,3.349)
\psline(2.451,3.651)(2.148,3.349)
\psellipse(2.6,3.8)(0.21,0.21)
\put(2.5,3.676){6}
\psline(3.051,3.051)(2.749,2.749)
\psellipse(3.2,3.2)(0.21,0.21)
\put(3.1,3.076){2}
\psline(1.549,3.651)(1.851,3.349)
\psline(1.252,3.651)(0.869,3.269)
\psellipse(1.4,3.8)(0.21,0.21)
\put(1.3,3.676){4}
\psline(2.148,3.051)(2.451,2.749)
\psline(1.851,3.051)(1.469,2.669)
\psellipse(2.0,3.2)(0.21,0.21)
\put(1.9,3.076){3}
\psline(2.451,2.451)(2.069,2.069)
\psellipse(2.6,2.6)(0.21,0.21)
\put(2.5,2.476){1}
\psline(0.869,3.131)(1.331,2.669)
\psellipse[fillstyle=solid,fillcolor=black](0.8,3.2)(0.098,0.098)
\psline(1.4,2.502)(1.4,2.098)
\psline(1.469,2.531)(1.931,2.069)
\psellipse[fillstyle=solid,fillcolor=black](1.4,2.6)(0.098,0.098)
\psline(1.4,1.902)(1.4,1.498)
\psellipse[fillstyle=solid,fillcolor=black](1.4,2.0)(0.098,0.098)
\psline(1.931,1.931)(1.469,1.469)
\psellipse[fillstyle=solid,fillcolor=black](2.0,2.0)(0.098,0.098)
\psline(1.331,1.331)(0.869,0.869)
\psellipse[fillstyle=solid,fillcolor=black](1.4,1.4)(0.098,0.098)
\psline(0.731,0.731)(0.269,0.269)
\psellipse[fillstyle=solid,fillcolor=black](0.8,0.8)(0.098,0.098)
\psellipse[fillstyle=solid,fillcolor=black](0.2,0.2)(0.098,0.098)

\end{pspicture*}

\begin{pspicture*}(0,0)(4.2,4.6)

\psline(0.651,4.251)(0.348,3.948)
\psline(0.948,4.251)(1.252,3.948)
\psellipse(0.8,4.4)(0.21,0.21)
\put(0.7,4.276){5}
\psline(1.4,4.19)(1.4,4.01)
\psellipse(1.4,4.4)(0.21,0.21)
\put(1.3,4.276){3}
\psline(2.148,4.251)(2.531,3.869)
\psline(1.851,4.251)(1.549,3.948)
\psellipse(2.0,4.4)(0.21,0.21)
\put(1.9,4.276){6}
\psline(2.669,3.731)(3.131,3.269)
\psline(2.531,3.731)(2.069,3.269)
\psellipse[fillstyle=solid,fillcolor=black](2.6,3.8)(0.098,0.098)
\psline(3.131,3.131)(2.669,2.669)
\psellipse[fillstyle=solid,fillcolor=black](3.2,3.2)(0.098,0.098)
\psline(0.348,3.651)(0.651,3.349)
\psellipse(0.2,3.8)(0.21,0.21)
\put(0.1,3.676){4}
\psline(1.549,3.651)(1.931,3.269)
\psline(1.252,3.651)(0.948,3.349)
\psellipse(1.4,3.8)(0.21,0.21)
\put(1.3,3.676){2}
\psline(2.069,3.131)(2.531,2.669)
\psline(1.931,3.131)(1.469,2.669)
\psellipse[fillstyle=solid,fillcolor=black](2.0,3.2)(0.098,0.098)
\psline(2.531,2.531)(2.069,2.069)
\psellipse[fillstyle=solid,fillcolor=black](2.6,2.6)(0.098,0.098)
\psline(0.948,3.051)(1.331,2.669)
\psellipse(0.8,3.2)(0.21,0.21)
\put(0.7,3.076){1}
\psline(1.4,2.502)(1.4,2.098)
\psline(1.469,2.531)(1.931,2.069)
\psellipse[fillstyle=solid,fillcolor=black](1.4,2.6)(0.098,0.098)
\psline(1.4,1.902)(1.4,1.498)
\psellipse[fillstyle=solid,fillcolor=black](1.4,2.0)(0.098,0.098)
\psline(1.931,1.931)(1.469,1.469)
\psellipse[fillstyle=solid,fillcolor=black](2.0,2.0)(0.098,0.098)
\psline(1.331,1.331)(0.869,0.869)
\psellipse[fillstyle=solid,fillcolor=black](1.4,1.4)(0.098,0.098)
\psline(0.731,0.731)(0.269,0.269)
\psellipse[fillstyle=solid,fillcolor=black](0.8,0.8)(0.098,0.098)
\psellipse[fillstyle=solid,fillcolor=black](0.2,0.2)(0.098,0.098)

\end{pspicture*}
}

\noindent
This is enough to express all Schubert classes as polynomials in $h,s,t,q$,
and yields the last relation.
\end{proo}

\subsection{Adjoint and coadjoint varieties in type $E_8$}

In the cohomology of $E_8/P_8$, let $h$ denote the hyperplane class, $s$ 
resp. $t$
the classes corresponding to the elements $s_2s_4s_5s_6s_7s_8$ resp.
$s_6s_5s_4s_3s_2s_4s_5s_6s_7s_8$.

\begin{prop}
\label{pres-e8p8}
The quantum cohomology algebra
$QH^*(E_8/P_8,\Q)$ is the quotient of the polynomial algebra $\Q[h,s,t,q]$
by the relations
$
h^{14}s + 6\ h^{10}t - 3\ h^{8}s^{2} - 12\ h^{4}st - 10\ h^{2}s^{3} 
+ 3\ t^{2}\ , \ 
29\ h^{24} - 120\ h^{18}s + 15\ h^{14}t + 45\ h^{12}s^{2} - 30\ h^{8}st 
+ 180\ h^{6}s^{3} - 30\ h^{2}s^{2}t + 5\ s^{4}\ $,
and $-86357\ h^{30} + 368652\ h^{24}s - 44640\ h^{20}t - 189720\ h^{18}s^{2} 
+ 94860\ h^{14}st - 473680\ h^{12}s^{3} + 74400\ h^{8}s^{2}t 
- 1240\ h^{2}s^{3}t + 60\ hq.
$
\end{prop}
\begin{proo}
We know that $QH^*(E_8/P_8,\Q)$ is generated by $h,s$ and $t$ and that 
there are relations in degree $14,24, 30$.
The quantum Hasse diagram of $E_8/P_8$ begins as follows:\\
\centerline{\input{hasse-e8p8}}\\
From this we may compute the degrees of all the Schubert classes and of any 
product of two Schubert classes. Moreover we argue as before to compute 
such products, and display the
results on the companion paper \cite{giambelli}. We only indicate the 
computed products, and why they suffice to obtain the presentation. Note 
that these results were obtained with the help of a computer by S. Nikolenko 
and N. Semenov \cite{ns}.

We have
$$
s^2 = 2\s \poidsehuit 12321111 + 4 \s \poidsehuit 12222111
+ 2 \s \poidsehuit 11222211 + \s \poidsehuit 01222221 \textrm{ and }$$
$$s\cdot t = \s \poidsehuit 12432112 + 5 \s \poidsehuit 12332212 
+ 6 \s \poidsehuit 12322222
+ 2 \s \poidsehuit 12333211 + 7 \s \poidsehuit 12332221.
$$
This allows to compute some Giambelli
formulas up to degree 17. To compute $s^3$ we compute
$$
s \cdot \s \poidsehuit 01222221 = \s \poidsehuit 12432222 
+ \s \poidsehuit 12333222
+ \s \poidsehuit 12333321,
$$
which allows to get the Giambelli formulas up to degree 21.
To get the relation in degree 20 we compute
$$
\begin{array}{rcl}
t^2 & = & 4 \s \poidsehuit 23433212 + 2 \s \poidsehuit 13443212 
+ 7 \s \poidsehuit 23432222
+ 14 \s \poidsehuit 13433222\\
& & + 8 \s \poidsehuit 12443222 + 16 \s \poidsehuit 12433322.
\end{array}
$$
\noindent
To express the Schubert classes of degree 22 we need to compute $s^2t$: 
to this end we
compute
$$
t \cdot \s \poidsehuit 01222221 = \s \poidsehuit 23443222 
+ 3 \s \poidsehuit 13543323
+ 2 \s \poidsehuit 23444322 + \s \poidsehuit 12433322.
$$
To get the relation in degree 24 we need to know $s^4$, and for this we show
$$
\s \poidsehuit 01222221 ^ 2 = \s  \poidsehuit 23444322.
$$
Finally we need to express the Schubert classes of degree 28 in terms of the generators; for this
we need to compute $s^3t$ which follows from the equality
$$
\s \poidsehuit 01222221 \cdot \s \poidsehuit 12332221
=
\s \poidsehuit 24644323 + \s \poidsehuit 24654322 + 4 \s \poidsehuit 24554323
+ \s \poidsehuit 23654323.
$$
The relation in degree 30 follows from the Chevalley formula.
\end{proo}

\section{Semi-simplicity of the quantum cohomology and strange duality}
\label{section-ss}

Recall that in \cite{cmp} and \cite{cmp3} we proved that the localisation 
$$QH^*(X,\C)_{loc}:=QH^*(X,\C)[q^{-1}]$$ 
of the quantum parameter in 
the quantum cohomology algebra of any minuscule or cominuscule variety $X$ 
is semi-simple. A general argument using complex conjugation (see 
Theorem 2.1 in \cite{cmp3}), implies that in that case there exists an 
algebra involution on $QH^*(X,\C)_{loc}$ sending $q$ on its inverse and more 
generally a quantum cohomology class of degree $d$ to a quantum cohomology 
class of degree $-d$. We conjectured in \cite{cmp3} that the only rational 
homogeneous spaces $X$ with Picard number one such that $QH^*(X,\C)_{loc}$
is semi-simple were the minuscule and cominuscule ones.

\subsection{Coadjoint varieties}

\begin{theo}
\label{theo-non-semisimple}
  Assume that there exists a 
presentation of the classical cohomology algebra $H^*(X,\C)$ satisfying:
  \begin{itemize}
  \item there is no equation of degree multiple of $c_1(X)$, 
\item there is at least one equation of degree smaller than $c_1(X)$.
  \end{itemize}
Then the localised quantum cohomology algebra $QH^*(X,\C)$ is not
semi-simple.
\end{theo}

\begin{proo}
  First remark that to prove that the localised quantum cohomology algebra 
$QH^*(X,\C)$  is not semi-simple, it suffices to prove that the 
specialisation of the quantum cohomology algebra at $q=1$ is not semi-simple. 

Let us fix some notation for the presentation of the cohomology ring 
$H^*(X,\C)$. This ring can always be presented as a homogeneous complete 
intersection ring: the generators and the equations are homogeneous and 
there are as many equations as there are generators. Let us denote by 
$(X_i)_{i\in[1,s]}$ the generators and with $(F_i)_{i\in[1,s]}$ the 
equations. The generators will be called classical generators and the 
equations classical equations.

A presentation for the quantum cohomology algebra can be derived from 
a presentation of the classical cohomology algebra by taking the same 
classical generators together with $q$ as generators and by deforming
the classical  relations $F_i$ in relations $G_i$ of the same degree
but eventually with $q$ terms (see \cite{ST}). 
By hypothesis we get homogeneous equations of degree different from 
multiples of $c_1$. In particular any monomial with a $q$ term is also 
multiple of a classical generator. This implies that $X_i=0$ for all 
$i\in[1,s]$ and $q=1$ is a solution of the equations 
$((G_i=0)_{i\in[1,s]},q=1)$. Denote by $0$ this point. We want to prove 
that 0 is a multiple solution of the system $((G_i=0)_{i\in[1,s]},q=1)$.

For this, because the specialisation at $q=1$ of the quantum algebra 
$QH^*(X,\C)$ is a complete intersection algebra of dimension 0, we only 
need to prove that the familly $(d_0G_i)_{i\in[1,s]}$ of differentials is 
not linearly independant. But an equation $F_i$ is not deformed if 
$\deg(F_i)<\deg(q)=c_1(X)$. In that case $G_i=F_i$ and by minimality $F_i$ 
is at least quadratic in the variables $(X_i)_{i\in[1,s]}$ therefore 
its differential $d_0F_i$ vanishes. In particular  the familly 
$(d_0G_i)_{i\in[1,s]}$ of differentials is not linearly independant and 0
is a multiple solution.
\end{proo}

\begin{coro}
  Let $X$ be a rational homogeneous space with Picard group $\Z$. Then the 
localised quantum algebra $QH^*(X,\C)_{loc}$ is not semi-simple except maybe 
in the following cases:

\vskip 0.2 cm

\begin{tabular}{c|l}
  Type&Variety\\
\hline
$A_n$ & $\G(p,n+1)$ for all $p\in[1,n]$.\\
$B_n$ & $\G_Q(p,2n+1)$ for all $p\in[1,n]$ with $p$ even, for $p=1$ and 
$p=n$.\\
$C_n$ & $\G_\omega(p,2n)$ for all $p\in[1,n]$ with $p$ odd and for $p=n$.\\
$D_n$ & $\G_Q(p,2n)$ for all $p\in[1,n-1]$ with $p$ odd and for $p=n$.\\
$E_6$ & $E_6/P_1$ and $E_6/P_3$ (and $E_6/P_6$ and $E_6/P_5$ which are 
isomorphic to the previous ones).\\
$E_7$ & $E_7/P_2$, $E_7/P_4$, $E_7/P_5$ and $E_7/P_7$.\\
$E_8$ & $E_8/P_4$ and $E_8/P_6$.\\
$F_4$ & $F_4/P_1$ and $F_4/P_2$.\\
$G_2$ & $G_2/P_1$ and $G_2/P_2$.\\
\hline
\end{tabular}
\end{coro}

\begin{rema}
  We already proved in \cite{cmp3} that for $X$ minuscule or cominuscule the
localised quantum cohomology ring $QH^*(X,\C)_{loc}$ is semi-simple. In the 
next subsection we shall prove that the same result holds for adjoint non 
coadjoint varieties. The only cases left are therefore the following ones 
(for $G_2$ the coadjoint variety is a quadric):

\vskip 0.2 cm 

\begin{tabular}{c|l}
  Type&Variety\\
\hline
$B_n$ & $\G_Q(p,2n+1)$ for all $p\in[1,n-1]$ with $p$ even.\\
$C_n$ & $\G_\omega(p,2n)$ for all $p\in[2,n-1]$ with $p$ odd.\\
$D_n$ & $\G_Q(p,2n)$ for all $p\in[2,n-1]$ with $p$ odd.\\
$E_6$ & $E_6/P_3$ (and $E_6/P_5$ which is isomorphic to the previous one).\\
$E_7$ & $E_7/P_2$, $E_7/P_4$ and $E_7/P_5$.\\
$E_8$ & $E_8/P_4$ and $E_8/P_6$.\\
$F_4$ & $F_4/P_2$.\\
\hline
\end{tabular}
\end{rema}

\begin{coro}
  For $X$ a coadjoint variety for a semi-simple group different from $G_2$, 
the localised quantum cohomology algebra $QH^*(X,\C)$ is not semi-simple.
\end{coro}

\subsection{Adjoint non coadjoint varieties}

There are three adjoint non coadjoint varieties: $\G_Q(2,2n+1)$, $F_4/P_1$ 
and $G_2/P_2$. In this subsection we prove the following

\begin{theo}
\label{ss-adjoint}
  For $X$ an adjoint non coadjoint rational homogeneous space, the 
localised quantum cohomology algebra $QH^*(X,\Q)_{loc}$ is semi-simple. 
It therefore admits an algebra involution $\iota$ sending $q$ to its inverse. 
  \end{theo}

We prove this result by a case by case analysis. We start with type $B_n$:

\begin{prop}
\label{ss-Bn}
  The localised quantum cohomology algebra $QH^*(\G_Q(2,2n+1),\Q)_{loc}$ 
is semi-simple. It therefore admits an algebra involution sending $q$ to 
its inverse.
\end{prop}

\begin{proo}
  A presentation of the ring $QH^*(\G_Q(2,2n+1),\Z)$ (over the integers) has 
been given in \cite{BKTnouveau}. We will use a simplified presentation 
which is only valid over $\Q$ but this is enough for our purpose. 

Let us denote by $K$ resp. $Q$ the tautological subbundle resp. quotient 
bundle on $\G_Q(2,2n+1)$. The rank of $K$ is 2 while the rank of $Q$ is 
$2n-1$. Pulling these vector bundles to $G/T$ we get decompositions into line 
bundles of the pull-backs. We denote $x_1,x_2$ the Chern classes of the two
line bundles corresponding to $K$.
Because the pull-back is injective in cohomology, 
we may therefore write
$$c_t(K)=(1+x_1t)(1+x_2t)\ \ {\rm and}\ \ 
c_t(Q)=(1-x_1t)(1-x_2t)\prod_{i=3}^n(1-x_i^2t^2),$$
where $c_t$ denotes the Chern polynomial. The tautological exact sequence 
gives $c_t(K)c_t(Q)=1.$
We can compute explicit formulas for the Chern classes of $K$ and $Q$. Let us 
write $c_i=c_i(K)$ for $i\in[0,2]$ and $\psi_j=c_j(Q)$ for $j\in[0,n-1]$. 
We have the equalities:
$$
\begin{array}{l}
c_0=\psi_0=1,\ c_1=x_1+x_2,\ c_2=x_1x_2,\\ 
\psi_{2j}=(-1)^{j}  
\s_{j}(x_3^2,\cdots,x_n^2) + (-1)^{j-1} x_1x_2 
\s_{j-1}(x_3^2,\cdots,x_n^2),\\
\psi_{2j+1}=(-1)^{j+1} (x_1+x_2) 
\s_{j}(x_3^2,\cdots,x_n^2),\\
\end{array}$$
where $\s_j$ is the $j$-th elementary symmetric function in $n-2$ variables. 
Remark that $\s_{n-1}(x_3^2,\cdots,x_n^2)=0$ thus we have 
$\psi_{2n-2}=(-1)^n\s_{n-2}(x_3^2,
\cdots,x_n^2)$ and $\psi_{2n-1}=0$. This yields a presentation of the 
cohomology algebra
\begin{equation}
  \label{eq:presentation}
  H^*(\G_Q(2,2n+1),\Q)=\Q[c_1,c_2,(\psi_j)_{j\in[1,2n-2]}]/
(\Sigma_k(x_l^2)_{k\in[1,n]}),
\end{equation}
where $\Sigma_j$ is the $j$-th elementary symmetric function in $n$ variables.
This is the Borel presentation, see \cite{borel} or \cite{BGG}. Using the 
first $n-2$ relations we can express the Chern classes 
$(\psi_j)_{j\in[1,2n-2]}$ in terms of $c_1$ and $c_2$. We get for example 
the following formulas for $j\in[1,n-2]$:
$$ \s_j(x_3^2,\cdots,x_n^2)=(-1)^j\sum_{k=0}^jx_1^{2k}x_2^{2(j-k)}.$$
In particular, over $\Q$, only the last two equations are worth defining a 
presentation of the cohomology ring. These equations are:
\begin{equation}
  \label{eq:equations}
\psi_{2n-2}+c_1\psi_{2n-3}+c_2\psi_{2n-4}=0\ \ {\rm and}\ \ 
c_2\psi_{2n-2}=0.
\end{equation}
In these expressions, the $\psi$ terms can be removed but we shall keep them 
for the moment as they are easier to express in terms of Schubert varieties. 
This will be usefull to deform these equations in the quantum cohomology ring. 
Let us describe the classes $\psi_{2n-4}$, $\psi_{2n-3}$ and $\psi_{2n-2}$ 
as Schubert varieties. This is classical (see for example \cite{BKTnouveau}). 
Let $(W_p)_{p\in[1,n]}$ be a complete isotropic flag, we have the equalities:
$$
\begin{array}{l}
\psi_{2n-4}=2[\{V_2\in\G_Q(2,2n+1)\ /\ \dim(W_3\cap V_2)\geq 1\}]\\
\psi_{2n-3}=2[\{V_2\in\G_Q(2,2n+1)\ /\ \dim(W_2\cap V_2)\geq1\}]\\
\psi_{2n-2}=2[\{V_2\in\G_Q(2,2n+1)\ /\ \dim(W_1\cap V_2)\geq1\}]\\
\end{array}$$
where $[Z]$ denotes the cohomology class of $Z$. In terms of roots as 
described in Section \ref{section-hasse} (with notation as in \cite{bou}), 
we have $\psi_{2n-4}=2\s_{\a_1+\a_2}$, $\psi_{2n-3}=2\s_{\a_1}$ and 
$\psi_{2n-2}=2\s_{-\a_1}$.
%
%
The Chern clases $c_1$ and $c_2$ are also related to Schubert classes:
$$
\begin{array}{l}
c_1=[\{V_2\in\G_Q(2,2n+1)\ /\ \dim(W_2^\perp\cap V_2)\geq 1\}]\\
c_2=[\{V_2\in\G_Q(2,2n+1)\ /\ V_2\subset W_1^\perp\}]\\
\end{array}$$
or with the notation of Section \ref{section-hasse}: $c_1=\s_{\T-\a_{2}}$ 
and $c_2=\s_{\T-\a_2-\a_1}$.

To compute quantum cohomology, we need to deform the preceding relations 
using the quantum parameter $q$. Remark that the degree of $q$ is $
2n-2$ therefore only the last two equations in the presentation 
(\ref{eq:presentation}) can be deformed. We thus need to compute the quantum 
terms in the equations (\ref{eq:equations}).
Now using Example \ref{calcul-de-q}, we can easily compute the $q$ terms
in the first equation of (\ref{eq:equations}). We get
\begin{equation}
  \label{eq:equation-quant-1}
\psi_{2n-2}+c_1*\psi_{2n-3}+c_2*\psi_{2n-4}=2q.
\end{equation}
To deform the last equation, we only need to use the affine symmetries (see 
Subsection \ref{subsection-affine}). Indeed $\frac{1}{2}\psi_{2n-2}$ 
corresponds to the opposite of the cominuscule simple root. We get
  \begin{equation}
    \label{eq:equation-quant-2}
c_2*\psi_{2n-2}=2qc_2.
  \end{equation}
To study the semi-simplicity of the quantum cohomology ring, we only need to prove that when we set $q=1$, 
the scheme ${\rm Spec}(QH^*(\G_Q(2,2n+1),\Q)_{q=1})$ is reduced.
In turn, since $QH^*(\G_Q(2,2n+1),\Q)$ is the ring of $S_2$-invariants of
the quotient
$\Q[x_1,x_2]/(\psi_{2n-2}+c_1*\psi_{2n-3}+c_2*\psi_{2n-4}-2q,
c_2*\psi_{2n-2}-2qc_2)$, it is enough to check that this quotient is
reduced.

For this we solve equations (\ref{eq:equation-quant-1}) and 
(\ref{eq:equation-quant-2}).
We get
$$-\sum_{k=0}^{n-1}x_1^{2k}x_2^{2n-2-2k}=2q\ \ {\rm and}\ \ 
x_1^2x_2^2\sum_{k=0}^{n-2}x_1^{2k}x_2^{2n-4-2k}=2qx_1x_2.$$
To solve these equations for $q=1$, we first remove the solutions 
$(x_1,x_2)$ with $x_i=0$ and $x_{3-i}$ an $2(n-1)$-th root of $-2$. There are 
$4(n-1)$ distinct such solutions. Then we assume $x_1x_2\neq0$ and write 
$x_2=\lt x_1$. We get
$$-x_1^{2(n-1)}(1+\lt^2+\cdots+\lt^{2(n-1)})
=2\ \ {\rm and}\ \ x_1^{2(n-1)}
\lt(1+\lt^2+
\cdots+\lt^{2(n-2)})
=2.$$
Equating these relations and dividing by $x_1^{2(n-1)}$ we get:
$$1+\lt+\cdots+\lt^{2n-2}=0$$
therefore $\lt$ is a $(2n-1)$-th root of unity different from 1. Remark 
that therefore $\lt\neq-1$ and we have 
$$1+\lt^2+\cdots+\lt^{2(n-1)}=\frac{\lt^{2n}-1}{\lt^2-1}=
\frac{\lt-1}{\lt^2-1}=\frac{1}{\lt+1}.$$
There are $2n-2$ different solutions for $\lt$ and then $x_1$ is a 
$(2n-2)$-th root of $-2(\lt+1)$ and there are $2n-2$ solutions for $x_1$. 
We described $4(n-1)+(2n-2)(2n-2)=2n(2n-2)$ solutions therefore the above 
affine scheme is reduced.
%
%
%
%
%
\end{proo}


\begin{prop}
  The localised quantum cohomology algebra $QH^*(F_4/P_1,\Q)_{loc}$ 
is semi-simple. It therefore admits an algebra involution sending $q$ to 
its inverse.
\end{prop}

\begin{proo}
  Recall that we have a presentation of the quantum cohomology ring given by
$$QH^*(F_4/P_1,\C)=\C[h,s,q]/(-h^8 + 12s^2 + 16q, 
3h^{12} - 18h^8s + 24h^4s^2 + 8s^3)$$
where $h$ is of degree 1 and $s$ of degree 4. Recall also that $c_1=8$. We
solve these equations in $h$ and $s$ when setting $q=1$. First remark that 
there is no solution on the closed subset $s=0$ neither on the closed subset 
$h=0$. Eliminating $h^8$ 
in the second equation using the first one we get
$$3h^{4}(5s^2+4)=4s(13s^2+18).$$
Taking the square of this equation and eliminating $h^8$ again we get
\begin{equation}
\label{equa-f4}
s^6-108s^4-576s^2-576=0.
\end{equation}
For $P(S)=S^3-108S^2-576S-576$ we have $P(-2)=136$ and
$P(0)=-576$, thus (\ref{equa-f4})
is is a degree 3 equation in $s^2$ whose solutions are real and 
distinct. The result follows.
\end{proo}

\begin{prop}
  The localised quantum cohomology algebra $QH^*(G_2/P_2,\Q)_{loc}$ 
is semi-simple. It therefore admits an algebra involution sending $q$ to 
its inverse.
\end{prop}

\begin{proo}
   Recall that we have a presentation of the quantum cohomology ring given by
$$QH^*(G_2/P_2,\C)=\C[h,q]/(h^6-18h^3q-27q^2)$$
where $h$ is of degree 1. Recall also that $c_1=3$. We solve this equation 
in $h$ when setting $q=1$. This is a degree 2 equation in $h^3$ whose 
solutions are real non vanishing and distinct.
\end{proo}

\subsection{Some properties of the involution for adjoint non coadjoint 
varieties}

Recall that in the minuscule and cominuscule cases, it was proved in 
\cite{cmp} and \cite{cmp3} that the involution $\iota$ given by complex 
conjugation has the nice property of sending a Schubert class to a 
multiple of such a class. In this section we shall see that contrary to 
these cases this property is not satisfied.

We shall first prove a nice property of the involution $\iota$, namely it 
sends a class to a multiple of itself and the quantum parameter $q$ for 
classes of degree $c_1$.

\begin{prop}
\label{prop-inv-c1}
For $X$ an adjoint non coadjoint rational homogeneous space, let $\s$ be a 
Schubert class of degree $dc_1$ for some integer $d$, then we 
have
$$\iota(\s)=\frac{\s}{q^{2d}}$$
\end{prop}

\begin{proo}
  To prove this result, we may restrict ourselves to the specialisation of the 
quantum cohomology obtained by setting $q=1$. The classes of this 
specialisation can be seen as functions with values in $\C$ on the variety 
${\rm Spec}(QH^*(X,\C)_{q=1})$. The involution $\iota$ is given by complex 
conjugation on these functions. It is thus enough to prove that the values of 
these functions for classes of degree multiple by $c_1$ are real.

For this we prove that $\s$ annihilates a polynomial equation whose solutions 
are real. This polynomial equation will be the minimal polynomial of the 
endomorphism given by the multiplication by $\s$. Let us first remark that it 
is enough to compute the minimal polynomial of the multiplication by $\s$ on 
the subspace of degree $c_1$ classes. Indeed, assume that $P(\s)$ vanishes on 
that subspace, then $P(\s)q=0$ and therefore because we inverted $q$ we have
$P(\s)=0$.

We therefore restrict ourselves on the subspace $H$ of degree $c_1$ classes 
and prove that the endomorphism given by multipication by $\s$ is adjoint for
a non degenerate real quadratic form. This will prove the result because the
eigenvalues of such endomorphisms are real.

Let us define a pairing on $H$ by 
$$(a,b)=\textrm{coefficient of $q^2$ in the product $a*b$}.$$
This pairing is clearly symmetric and we shall prove that the quantum Schubert 
classes of degree $c_1$ form an orthogonal basis for this pairing. Indeed, 
denote by ${\rm pt}$ the class of a point in $X$, this orthogonality result 
comes from a formula for the multiplication by the class ${\rm pt}$. This 
formula is a consequence of a more general formula for the multiplication by 
the Schubert class corresponding to the reflection with respect to the highest 
root (see \cite[Proposition 11.2]{LS}) which is ${\rm pt}$ in our situation: 

\begin{fact}
  We have, for $\a$ a long simple root and $\varpi_\a$ the associated 
fundamental weight, the formula:
$${\rm pt}*\s_{-\a}=\sca{\varpi_\a}{\Theta^\vee}q^2\s_{\a}$$
\end{fact}

In particular the quantum Schubert classes in $H$ form an orthogonal basis 
and the form is positive definite. But now we may compute
$$(\s*a,b)=\textrm{coefficient of $q^2$ in the product $\s*a*b$}=(a,\s*b)$$
and $\s$ is adjoint form the form $(\ ,\ )$. This concludes the proof.
\end{proo}


However, we will now prove that the involution $\iota$ is not as simple as in 
the minuscule or cominuscule cases. In particular, we prove the following

\begin{theo}
\label{pas-classe}
  For $X$ an adjoint non coadjoint rational homogeneous space, there exists 
a Schubert class $\s$ such that $\iota(\s)$ is not a multiple of a Schubert
class.
\end{theo}

\begin{proo}  
 Remark that to prove this result, we may set $q=1$ so that the 
involution is the identity on degree $c_1$ classes. The proof will be a
case by case analysis even if the idea of the proof is always the same:
take $\s$ a class of degree a divisor of $c_1$, say $d\deg(\s)=c_1$. Then
we have $\deg(\s^d)=c_1$ therefore $\iota(\s)^d=\iota(\s^d)=\s^d$ by 
Proposition \ref{prop-inv-c1}. We are then reduced to prove that $\s$ is
the 
only Schubert class $\tau$ such that $\deg(\tau)\equiv-\deg(\s)\ ({\rm
  mod}\ c_1)$ and $\tau^d$ is a scalar multiple of $\s^d$. When 
$\deg(\s)\equiv-\deg(\s)\ ({\rm mod}\ c_1)$ then we will also need to
prove that $\iota(\s)$ is not a scalar multiple of $\s$.

Let us start with $X=G_2/P_2$. Recall that in this situation
$c_1=3$. Take $\s=h$. We have by Chevalley formula the equality
$\s^3=6\s_3+3q$ where $\s_3$ is the only Schubert class of degree
3. If the class $\iota(\s)$ is a multiple of a Schubert class, then it
has to be a multiple of $\s_2$ the unique degree 2 Schubert class or a
multiple of ${\rm pt}$ the class of a point. But $3\s_2=h^2$ therefore
we may compute $\s_2^3$ using the equation of the presentation of the
cohomology ring \ref{pres-g2p2} we get $3\s_2^3=2h^3q+3q^2$ and by
Chevalley formula $\s_2^3=4\s_3q+3q^2$ which is not a scalar multiple
of $\iota(\s^3)=6\s_3+3q$. We can also compute ${\rm
  pt}^3=2q^2\ell*{\rm pt}$ where $\ell$ is the class of a line in $X$
(see Corollary \ref{gw-point-au-carre}). But by Theorem 
\ref{coro_dmax_adjoint} we have ${\rm pt}*\ell=\lt q^2\s_3$ because the other
term would be $q^3$ (by \cite[Proposition 11.2]{LS} or several
applications of the Chevalley formula we even have ${\rm
  pt}*\ell=q^2\s_3$) therefore ${\rm pt}^3=2\lt q^4\s_3$ which is not
a scalar multiple of $\iota(\s^3)=6\s_3+3q$.

We now deal with $X=F_4/P_1$. Recall that in that situation
$c_1=8$. Take $\s$ to be the Schubert class associated to the element
$s_4s_3s_2s_1$ in the Weyl group. We have $\s^2=\s_{-\a_1}+\s_{-\a_2}$
where $\a_1$ and $\a_2$ are simple root.
This formula can be proved using the presentation of
Proposition \ref{pres-f4p1}, we proved it using our program
\cite{programme}. If
the class $\iota(\s)$ is a multiple of a Schubert class, then it has
to be a multiple of $\s$, of $\tau$ the other degree 4 Schubert class
or of $\s_{12}$ the unique degree 12 Schubert class. Using our
program, we have $\tau^2=6\s_{-\a_1}+8\s_{-\a_2}$ and
$\s_{12}^2=3q^2\s_{-\a_1}+5q^2\s_{-\a_2}+q^3$. Therefore none of them
is a multiple of $\iota(\s)^3=\s_{-\a_1}+\s_{-\a_2}$. To conclude, we
need to prove that $\iota(\s)$ is not equal to a multiple of $\s$. For
this we consider $\iota(\s)$ as a function on the scheme ${\rm
  Spec}(QH^*(X,\C)_{q=1})$.
In the basis $(q,\s_{-\a_1},\s_{-\a_2})$ of $H$ the matrix of
the multiplication by $\iota(\s)^2=\s^2$ is given by (once more we
use our program):
$$
\left(\begin{array}{ccc}
  0&2&3\\
1&5&9\\
1&6&10
\end{array}\right)$$
whose minimal polynomial is $P(x)=x^3+15x^2+9x+1$. It has two real positive
roots and one real negative root, and from $P(\s^2)=0$ we deduce that $\s$
has four real values and two complex
conjugate values. Its complex conjugate is therefore not a multiple of
itself.

We finish with $X=\G_Q(2,2n+1)$. Recall that in that situation
$c_1=2n-2$. We first prove that there is a Schubert class $\s$ such that its
square $\s^2$ is a multiple of $\s_{-\a_{n-1}}$. For this we need to discuss 
on the parity of $n$. We take for $\s$ the Schubert class associated to the 
Schubert variety (here we fix a complete isotropic flag $(W_i)_{i\in[1,n]}$):
$$\begin{array}{ll}
  \{V_2\in\G_Q(2,2n+1)\ /\ V_2\subset W_p^\perp\}&\textrm{if $n-1=2p$ 
is even,}\\
  \{V_2\in\G_Q(2,2n+1)\ /\ V_2\subset W_{p+2}\} & \textrm{if $n-1=2p+1$ 
is odd.}\\
\end{array}$$
Remark that the degree of $\s$ is given by
$$\deg(\s)=\left\{
  \begin{array}{ll}
    n-1 & \textrm{if $n-1$ is even,}\\
3(n-1) &  \textrm{if $n-1$ is odd.}\\
  \end{array}\right.$$

\begin{fact}
\label{fait-s^2}
  The square $\s^2$ is a scalar (and $q$) multiple of $\s_{\-\a_{n-1}}$.
\end{fact}

\begin{proo}
Let us deal with the case $n-1$ even first. By Example \ref{calcul-de-q}, 
there is no $q$ term in $\s^2$. This is of degree $2(n-1)$ and we need to 
compute the classical products $\s\cup\s\cup\s_{\a_k}$ for $\a_k$ a simple 
root and $k\in[1,n-1]$.
The Schubert variety associated to the class $\s_{\a_k}$ is 
$$\{V_2\in\G_Q(2,2n+1)\ /\ V_2\subset W_{k-1}^\perp\ \textrm{and}\ 
\dim(V_2\cap W_{k+1})\geq1\}$$
where $(W_i)_{i\in[1,n]}$ is a complete isotropic flag. We therefore want to 
count the number of elements $V_2\in\G_Q(2,2n+1)$ such that
$$V_2\subset W_p^\perp,\ V_2\subset T_p^\perp,\ V_2\subset U_{k-1}^\perp\ 
\textrm{and } \dim(V_2\cap U_{k+1})\geq1$$
where $(U_i)_{i\in[1,n]}$, $(T_i)_{i\in[1,n]}$ and $(W_i)_{i\in[1,n]}$ 
are three generic complete isotropic flags. We see that these conditions imply 
that $\dim(U_{k+1}\cap W_p^\perp\cap T_p^\perp)\geq1$ and because these flags 
are in general position that $k+1-2p\geq1$ therefore $k\geq n-1$. In particular
the only Schubert class appearing in $\s^2$ is $\s_{-\a_{n-1}}$ and the 
result follows. We even know (from \cite{FW}) that this product is not zero 
and easy manipulations with isotropic subspaces show that 
$\s^2=2\s_{-\a_{n-1}}$ (we will not need this fact here).

Now assume $n-1$ to be odd. By Corollary \ref{coro_dmax_adjoint}, there is no 
$q^3$ term in $\s^2$, and by Example \ref{calcul-de-q}, there is no
$q$ term,
thus there are only terms of the form $q^2\s_{-\a_k}$ for $\a_k$ a simple 
root and $k\in[1,n-1]$. We 
therefore need to compute $I_2(\s,\s,\s_{\a_k})$. But if 
$f:\pu\to\G_Q(2,2n+1)$ is a degree 2 morphism then all the subspaces in the 
image are contained in a 4-dimensional, non isotropic in general, subspace 
$V_4$. We shall use this subspace to get some information. Indeed, if $f$ 
meets three general Schubert varieties representing the Schubert classes $\s$, 
$\s$ and $\s_{\a_k}$ then the subspace $V_4$ satisfies:
$$
\begin{array}{l}
  \dim(V_4\cap W_{p+2})\geq2,\\ 
\dim(V_4\cap T_{p+2})\geq2,\\ 
V_4\cap U_{k-1}^\perp \textrm{ contains a 2-dimensional isotropic subspace 
and }\\ 
\dim(V_4\cap U_{k+1})\geq1,
\end{array}$$
where $(U_i)_{i\in[1,n]}$, $(T_i)_{i\in[1,n]}$ and $(W_i)_{i\in[1,n]}$ 
are three generic complete isotropic flags. We see that these conditions imply 
that $\dim(U_{k+1}\cap (W_{p+2}+ T_{p+2}))\geq1$ and because these flags 
are in general position that $k+1+p+2+p+2-(2n+1)\geq1$ therefore $k\geq n-1$.
We conclude as before. Again one can prove the equality 
$\s^2=2q^2\s_{-\a_{n-1}}$.
\end{proo}

We now prove that there is no Schubert class $\tau$ different from $\s$ 
with $\deg(\tau)\equiv-\deg(\s)\equiv n-1\ ({\rm mod}\ c_1)$ and $\tau^2$ 
equal to a scalar multiple of $\s_{-\a_{n-1}}$. For this we describe the 
Schubert varieties associated to such classes $\tau$ of degree equal to $n-1$
modulo $c_1$:
$$
\begin{array}{l}
\{V_2\in\G_Q(2,2n+1)\ /\ V_2\subset W_{a}^\perp\ \textrm{and}\ 
\dim(V_2\cap W_{a+b+1})\geq1\}\ \textrm{with $a\in[0,\frac{n-1}{2}[$ 
and $2a+b=n-1$ 
}\\
\{V_2\in\G_Q(2,2n+1)\ /\ V_2\subset W_{a+b+2}\ \textrm{and}\ 
\dim(V_2\cap W_{a+1})\geq1\}\ \textrm{with $a\in[0,\frac{n-2}{2}[$ and 
$2a+b=n-2$}
\\
\end{array}
$$
where $(W_i)_{i\in[1,n]}$ is a complete isotropic flag. We want to compute 
the coefficient of $\s_{-\a_{n-1}}$ (resp. $q^2\s_{-\a_{n-1}}$) is the square 
$\tau^2$ in the first (resp. second) example. 

\begin{fact}
  The only Schubert class $\tau$ of degree $n-1$ or $3(n-1)$ such that 
$\tau^2$ is a multiple of $\s^2$ by a scalar product of a power of $q$ are the classes $\pm\s$.
\end{fact}

\begin{proo}
  For $\tau$ a degree $n-1$ class, we proceed as follows: first we prove 
that with notation 
as in Proposition \ref{ss-Bn}, we have $\s_2^{n-1}=\s^2$ (recall that $\s_2$ 
is the degree 2 cohomology class of the Schubert variety 
$\{V_2\in\G_Q(2,2n+1)\ /\ 
V_2\subset W_1^\perp\}$ with $(W_i)_{i\in[1,n]}$ is a complete isotropic flag).
Indeed, we can easily compute using jeu de taquin (cf. \cite{littlewood}) that 
$\s_2^k=\s_{2k}$ where $\s_{2k}$ is the cohomology class of the Schubert 
variety $\{V_2\in\G_Q(2,2n+1)\ /\ V_2\subset W_k^\perp\}$ with 
$(W_i)_{i\in[1,n]}$ is a complete isotropic flag and $k\leq n-2$. Then we 
ealisy see (use for example the heaps like in \cite{littlewood}) that the 
only degree $2n-2$ class containing $\s_{n-1}$ is $\frac{1}{2}\s^2$ and from 
Example \ref{calcul-de-q} there is no $q$ term in the product $\s_2*\s_{n-2}$. 
Therefore $\s_2^{n-1}$ is a scalar multiple of $\frac{1}{2}\s^2$. From this 
one easily get $\s_2^{n-1}=\s^2$ (we will only need that it is a scalar 
multiple of $\s^2$).

Recall from our presentation of 
the quantum cohomology in Proposition \ref{ss-Bn} that any degree $n-1$ 
class $\tau$ can be expressed as a symmetric degree $n-1$ polynomial in the 
variables $x_1$ and $x_2$. Moreover, the first equation being in degree $2n-2$ 
and of the form $q=P(x_1,x_2)$ with $P$ a degree $2n-2$ polynomial, we have 
that, in the subring generated by $x_1$ and $x_2$ in $QH^*(\G_Q(2,2n+1),\C)$,
the multiplication is injective from degree $n-1$ classes to degree $2n-2$ 
classes. In particular, the only classes $\tau$ such that $\tau^2=\s^2=
\s_2^{n-1}$ is $\pm\s_2^{p}$ with $2p=n-1$ for $n-1$ even and no class for 
$n-1$ odd. Now $\s_2^p=\s_{2p}$ as above and is the class $\s$.

For the second case, we were not able to do the same trick and we need to do 
geometry on subspaces as in Fact \ref{fait-s^2}. We need to compute the square 
of the class $\tau$ of a Schubert variety 
$\{V_2\in\G_Q(2,2n+1)\ /\ V_2\subset W_{a+b+2}\ \textrm{and}\ 
\dim(V_2\cap W_{a+1})\geq1\}$ with $a\in[0,\frac{n-2}{2}[$ and 
$2a+b=n-2$ and where $(W_i)_{i\in[1,n]}$ is a complete isotropic flag. To 
prove that its square is not a scalar multiple of $q^2\s^2$, we compute its 
coefficient on the class $\s_{-\a_{n-b-1}}$.
We therefore have to compute the Gromov-Witten invariant 
$I_2(\tau,\tau,\s_{\a_{n-b-1}})$. For this we want to count the number of 
4-dimensional (non isotropic) subspaces $V_4$ such that
$$\begin{array}{ll}
\dim(V_4\cap W_{a+1})\geq1, &  \dim(V_4\cap W_{a+b+2})\geq2,\\ 
\dim(V_4\cap T_{a+1})\geq1,  &\dim(V_4\cap T_{a+b+2})\geq2,\\ 
\dim(V_4\cap U_{n-b})\geq1, 
& V_4\cap U_{n-b-2}^\perp 
\textrm{ contains a 2-dimensional isotropic subspace,}\\ 
\end{array}$$
where $(U_i)_{i\in[1,n]}$, $(T_i)_{i\in[1,n]}$ and $(W_i)_{i\in[1,n]}$ 
are three generic complete isotropic flags.

\vskip .2cm

Suppose we are given such $V_4$, and let us analyse some properties of this
linear subspace. Since
$(U_i)_{i\in[1,n]}$, $(T_i)_{i\in[1,n]}$ and $(W_i)_{i\in[1,n]}$ are
generic, the above
inequalities are equalities. We denote $P$ the intersection point of
$\p V_4$ and $\p U_{n-b}$. We denote
$L_1=V_4\cap W_{a+b+2}$ and $L_2=V_4\cap T_{a+b+2}$.
In the plane $\p (V_4 \cap P^\perp)$ there are two isotropic lines
through $P$, 
the line $\p L$ not meeting $\p L_1$ and $\p L_2$ and the line $\p
L_3$. 
%
%
In the plane $\p \pi_P(V_4)$ we denote $p = \p \pi_P(L)$.
In $\p \pi_L(L_4)$ we denote $\Lambda_3$ the projection of $\p L_3$.
Note that $\pi_L(V_4) = \pi_L(V_4 \cap W_{a+1}) \oplus 
\pi_L(V_4 \cap T_{a+1})$ and therefore
$\Lambda_3 \in \p ( \pi_L(W_{a+1}) \oplus \pi_L(T_{a+1} )) $.
Finally
we denote $\pi_3$ the projection with center generated by $L$ and $L_3$
and we observe that
$\Lambda := \p \pi_3 (V_4) \in \pi_3(W_{a+1}) \cap \pi_3(T_{a+1})$.

\vskip .2cm

We now explain how to find exactly two subspaces $V_4$ satisfying the above
conditions given the isotropic flags.

First of all
we see that 
$U_{n-b}\cap (W_{a+b+2} + T_{a+b+2})$ has dimension $n-b+a+b+2+a+b+2-(2n+1)=1$.
Thus the point $P$ is determined by
$P = \p (U_{n-b}\cap (W_{a+b+2} + T_{a+b+2}))$. We then have
$\dim \pi_P ( W_{a+b+2} ) = \dim \pi_P ( T_{a+b+2} ) = a+b+2$ whereas
$\dim \pi_P ( W_{a+b+2} + T_{a+b+2} ) = 2(a+b+2) - 1$, thus
the spaces $\p \pi_P ( W_{a+b+2} )$ and $\p \pi_P ( T_{a+b+2} ) $
intersect in one
point. Note moreover that $\pi_P(W_{a+b+2}) \supset \pi_P(L_1) \supset \pi_P(L)$,
and similarly for $\pi_P(T_{a+b+2})$, and thus $p$ is detemined by
$p = \p (\pi_P ( W_{a+b+2} ) \cap \pi_P ( T_{a+b+2} )) $. This defines
uniquely the line $\p L$ in $\p \C^{2n+1}$.

Now we observe that $\p (U_{n-b-2}^\perp \cap P^\perp)$ meets $\p L$
at $P$ (except for $(a,b)=(0,n-2)$), so that $U_{n+b+1} := \pi_L(
U_{n-b-2}^\perp \cap P^\perp) $ is a quadratic space 
of dimension $n+b+1$ (for $(a,b)=(0,n-2)$ the space $\pi_L(
U_{n-b-2}^\perp\cap P^\perp)$ is of dimension $2n-2$). The intersection 
$U_{n+b+1} \cap (\pi_L(W_{a+1}) + \pi_L(T_{a+1}))$ has dimension 2 and is
non-degenerate for this quadratic form (for $(a,b)=(0,n-2)$, the
space $\pi_L(W_{a+1}) + \pi_L(T_{a+1})$ is already of dimension 2 with
a non degenerate quadratic form). Since we have seen that
$\Lambda_3 \in \p ( \pi_L(W_{a+1}) \oplus \pi_L(T_{a+1} ) )$, it follows
that we have exactly two choices for $\Lambda_3$: the two isotropic
points of $\p ( U_{n+b+1} \cap (\pi_L(W_{a+1}) + \pi_L(T_{a+1})))$.
Let us choose one isotropic point $\Lambda_3$. This defines the projection
$\pi_3$. Observe that $\pi_3 ( W_{a+1} )$ and
$\pi_3(T_{a+1})$ meet in dimension 1 in $\pi_3 (  W_{a+1} \oplus T_{a+1} )$
of dimension $2a+1$, and therefore $\Lambda$ is defined
by $\Lambda = \p ( \pi_3 ( W_{a+1} ) \cap \pi_3(T_{a+1}))$. This defines
$V_4$ as the inverse image of $\Lambda$ by $\pi_3$.
\end{proo}

To conclude the proof we need to prove the following fact:

\begin{fact}
  We have $\iota(\s)\neq\pm\s$.
\end{fact}

\begin{proo}
  For this we first specialise to $q=1$ and we will prove that 
$\s^2=2\s_{-\a_{n-1}}$ takes positive and negative real values. This will 
imply that for its square root the complex 
conjugation (i.e. the involution $\iota$) is never multiplication by $\pm1$. 

Consider the multiplication by $\s_{-\a_{n-1}}$ as an endomorphism $u$
of the space $H$ (the vector space generated by the degree $2n-2$ classes) 
and adjoint for the scalar product $(\ ,\ )$ defined in the proof of 
Proposition \ref{prop-inv-c1}. First of all because the quantum multiplication 
of Schubert classes has non negative integer coefficients, we see that 
${\rm tr}(u)>0$ and therefore $u$ has positive eigenvalues. If all 
eigenvalues of $u$ where non negative, then we would have $(u(x),x)\geq0$ 
for all $x\in H$ with equality if and only if $u(x)=0$. But take 
$x=1+\s_{-\a_1}$ (recall that multipication by $\s_{-\a_1}$ induces the 
affine symmetry described in Subsection \ref{subsection-affine}). We have 
$u(x)=2\s_{-\a_{n-1}}$ which is non trivial but $(u(x),x)=0$. Therefore $U$ 
has negative eigenvalues and the result follows.
\end{proo}

This concludes the proof: $\iota(\s)$ is never a Schubert class.
\end{proo}

\begin{rema}
 Let us compute explicitely in an example the element $\iota(\s)$ with $\s$ 
as in the previous proposition. Assume that $X=\G_Q(2,7)$. We have $c_1=4$. 
In that case, we have $\s=\s_2$ with the notation of Proposition \ref{ss-Bn}. 
Recall also that there is a unique other Schubert class of degree 2 denoted 
$\tau_2$. We have $\s^2=2\s_{-\a_{ad}}$ and easy computations give the 
following equation:
$$s^2(\s^4-4\s^2-16)=0.$$
Therefore the values of $\s^2$ are 0, $2(1+\sqrt{5})$ and $2(1-\sqrt{5})$. 
In particular we see that the polynomial 
$$Q(x)=\frac{x(x^2-2)}{2\sqrt{5}}$$
maps the square roots of the values of $\s^2$ to their complex conjugate 
therefore 
$$\iota(\s)=\frac{s^3-2sq}{2q^2\sqrt{5}}=\frac{1}{q^2\sqrt{5}}
(2\ell+2\tau_2-\s_2)$$
where $\ell$ is the class of a line.
\end{rema}

\section{The incidence variety}
\label{section-incidence}

In this section we deal with the adjoint variety $X$ of type $A_n$,
namely the incidence variety between points and hyperplanes in
$\p^n$. For many aspects the results proved in the previous sections
are still valid. However, as the Picard group of $X$ is $\Z^2$, some
of the results change. We will explain in this section which of the
results proved remain valid and which changes are required for the
others. 

\begin{fact}
  (\i) The parametrisation of Schubert classes by roots in Section
  \ref{section-hasse} remains valid. In particular Proposition \ref{chev} is
  valid.

(\i\i) As the Picard group is $\Z^2$, there are two classes of degree
1 that we shall denote by $h_1$ and $h_2$ and there are two quantum
parameters $q_1$ and $q_2$. Therefore, the parametrisation of quantum
monomials by roots of the affine root system fails.

(\i\i\i) Chevalley and Quantum Chevalley formulas also fail. However,
setting $h=h_1+h_2$ and $q=q_1+q_2$, then the formula in Theorem
\ref{theo-chevalley} is true.

(\i v) The formula for affine symmetries fails by lack of a
description of the quantum monomials.

(v) The description of $\varpi$-minuscule elements is valid as well
as the fact that $\dmax=2$.
\end{fact}

It is easy to derive an explicit formula for the quantum Chevalley
formula. This formula completly determines the quantum cohomology ring
as the degree 1 classes are generators. We get the following:

%

\begin{prop}
%
%
The quantum cohomology ring $QH^*(X,\Q)$ is generated by $h_1$, $h_2$
and the quantum parameters $q_1$ and $q_2$. It is the quotient of
$\Q[h_1,h_2,q_1,q_2]$ by the relations 
$$\sum_{k=0}^nh_1^k(-h_2)^{n-k}=q_1+(-1)^nq_2\ {\rm and}\
h_1^{n+1}=q_1(h_1+h_2).$$ 

(\i\i\i) The quantum cohomology ring is semi-simple for quantum parameters 
$q_1$ and $q_2$ satisfying
$$q_1\neq0,\ q_2\neq0\  \textrm{and}\ q_1+(-1)^nq_2\neq0.$$
\end{prop}

\begin{proo}
  The first two formulas are direct interpretations of the quantum
  Chevaley formula given in \cite{FW}. 

Let us prove the semi-simplicity. By multiplication 
of the quantum parameters $q_1$ and $q_2$ by the same scalar, we may assume 
$q_1+(-1)^nq_2=1$. Let us first remark that because $q_1\neq0$, if $h_1=0$
then from the second equation we have $h_2=0$. This is not compatible with the
first equation. We may thus assume $h_1\neq0$ and set $h_2=\lt h_1$. We obtain
$$h_1^n\sum_{k=0}^n(-\lt)^k=1\ {\rm and}\
h_1^n=q_1(1+\lt).$$
Remark that because $h_1\neq0$ we have $1+\lt\neq0$ therefore, multiplying the
first equation by $q_1(1+\lt)$ and equating, we get
$h_1^nq_1(1+(-1)^n\lt^{n+1})=h_1^n$ and the equation
$$q_1\lt^{n+1}=(-1)^{n+1}(q_1-1).$$
Because $q_1\neq0$ (otherwise $h_1=0$) and $q_1\neq1$ (otherwise $q_2=0$) we
get $n+1$ distinct solutions for $\lt$ and $n(n+1)$ solutions for
$(h_1,h_2)$.
\end{proo}


\bigskip\noindent
Pierre-Emmanuel {\sc Chaput}, \\
{\it Laboratoire de Math{\'e}matiques Jean Leray,} 
UMR 6629 du CNRS, UFR Sciences et Techniques,  2 rue de la
Houssini{\`e}re, BP 92208, 44322 Nantes cedex 03, France and \\
{\it Hausdorff Center for Mathematics,}
Universit{\"a}t Bonn, Villa Maria, Endenicher
Allee 62, 
53115 Bonn, Germany

\noindent {\it email}: \texttt{pierre-emmanuel.chaput@math.univ-nantes.fr}.

\medskip\noindent
Nicolas {\sc Perrin}, \\
{\it Hausdorff Center for Mathematics,}
Universit{\"a}t Bonn, Villa Maria, Endenicher
Allee 62, 
53115 Bonn, Germany and \\
{\it Institut de Math{\'e}matiques de Jussieu,} 
Universit{\'e} Pierre et Marie Curie, Case 247, 4 place
Jussieu, 75252 Paris Cedex 05, France.

\noindent {\it email}: \texttt{nicolas.perrin@hcm.uni-bonn.de}.

\end{document}